\documentclass[12pt, a4paper]{article}
\usepackage{setspace}
\usepackage{geometry}
\usepackage[OT1]{fontenc}
\usepackage{amsthm,amsmath}
\usepackage[authoryear,nonamebreak,bibstyle]{natbib}
\usepackage[colorlinks,citecolor=blue,urlcolor=blue]{hyperref}
\usepackage{amssymb,amsmath,bm,mathrsfs,makeidx,amsfonts,graphicx,amsthm,multirow}
\usepackage{booktabs}
\def\be{\begin{equation}}
\def\ee{\end{equation}}
\def\bea{\begin{eqnarray}}
\def\eea{\end{eqnarray}}

\numberwithin{equation}{section}
\theoremstyle{plain}
\newtheorem{thm}{Theorem}[section]
\newtheorem{prop}{Proposition}
\newtheorem{lem}{Lemma}

\newtheorem{rmk}{Remark}

\newcommand{\non}{\nonumber \\}
\newcommand{\bbA}{{\bf A}}

\newcommand{\bbB}{{\bf B}}
\newcommand{\bbC}{{\bf C}}

\newcommand{\bbD}{{\bf D}}
\newcommand{\bbd}{{\bf d}}
\newcommand{\bbe}{{\bf e}}
\newcommand{\bbE}{{\bf E}}

\newcommand{\bbF}{{\bf F}}

\newcommand{\bbG}{{\bf G}}

\newcommand{\bbH}{{\bf H}}

\newcommand{\bbI}{{\bf I}}

\newcommand{\bbM}{{\bf M}}

\newcommand{\bbN}{{\bf N}}

\newcommand{\bbQ}{{\bf Q}}

\newcommand{\bbr}{{\bf r}}

\newcommand{\bbS}{{\bf S}}

\newcommand{\bbT}{{\bf T}}

\newcommand{\bbu}{{\bf u}}

\newcommand{\bbv}{{\bf v}}

\newcommand{\bbX}{{\bf X}}

\newcommand{\bbx}{{\bf x}}

\newcommand{\bbY}{{\bf Y}}
\newcommand{\bby}{{\bf y}}

\newcommand{\bbz}{{\bf z}}

\newcommand{\ep}{\ensuremath{\epsilon}}

\begin{document}

\begin{center}
\Large{High Dimensional Correlation Matrices: \\CLT and Its Applications}
\end{center}

\begin{center}
Jiti Gao\\
Department of Econometrics and Business Statistics\\
Monash University, Caulfield East, VIC 3145\\
email:\texttt{jiti.gao@monash.edu}\\ \vspace{0.2in}
Xiao Han\\
Division of Mathematical Sciences, Nanyang Technological University\\
 Singapore, 637371\\ email: \texttt{xhan011@e.ntu.edu.sg}\\  \vspace{0.2in}
Guangming Pan\\
Division of Mathematical Sciences, Nanyang Technological University\\  Singapore, 637371 \\ email: \texttt{gmpan@ntu.edu.sg}\\ \vspace{0.2in}
Yanrong Yang\\
Department of Econometrics and Business Statistics\\
Monash University, Caulfield East, VIC 3145\\ email:\texttt{yanrong.yang@monash.edu}
\end{center}

\newpage
%
\begin{center}
\textbf{Abstract}
\end{center}
Statistical inferences for sample correlation matrices are important in high dimensional data analysis. Motivated by this, this paper establishes a new central limit theorem (CLT)  for a linear spectral statistic (LSS) of high dimensional sample correlation matrices for the case where the dimension $p$ and the sample size $n$ are comparable. This result is of independent interest in large dimensional random matrix theory. Meanwhile, we apply the linear spectral statistic to an independence test for $p$ random variables, and then an equivalence test for $p$ factor loadings and $n$ factors in a factor model. The finite sample performance of the proposed test shows its applicability and effectiveness in practice. An empirical application to test the independence of household incomes from different cities in China is also conducted.
\vspace*{.3in}


\newpage

\section{Introduction}

Big data issues arising in various fields bring great challenges to classical statistical inferences. High dimensionality and large sample size are two critical features of big data. In statistical inferences, there are serious problems, such as, noise accumulation, spurious correlations, and incidental homogeneity, arisen by high dimensionality. In view of this, the development of new statistical models and methods is necessary for big data research. Thus, our task in this paper is to analyze the correlation matrix of a $p$-dimensional random vector $\bbx=(X_1, X_2, \ldots, X_p)^{*}$, with available samples $\bbx_1, \bbx_2, \ldots, \bbx_n$, where $\bbx_i=(X_{1i}, X_{2i}, \ldots, X_{pi})^{*}$, where $*$ denotes the conventional conjugate transpose. We consider the setting of the dimensionality $p$ and the sample size $n$ being in the same order.

Correlation matrices are commonly used in statistics to investigate relationships among different variables in a group. It is well known that the sample correlation matrix is not a `good' estimator of its corresponding population version when the number $p$ of random variables under investigation is comparable to the sample size $n$. Thus, it is of great interest to understand and investigate the asymptotic behaviour of the sample correlation matrices of high dimensional data. Sample correlation matrices have appeared in some classical statistics for hypothesis tests. \cite{S2005} utilized sample correlation matrices to test independence for a large number of random variables having a multivariate normal distribution. Concerning statistical inference for high dimensional data, furthermore, there are many available research methods based on sample covariance matrices such as \cite{J2001}. 
 As the population mean and variance of the original data are usually unknown, sample covariance matrices cannot provide us with sufficient and correct information about the data. To illustrate this point, a simple example is that we will make an incorrect conclusion in an independence test if the variance of the data under investigation is not identical to one while the statistics based on sample covariance matrices require the variance to be one. Moreover, the main advantage of using sample correlation matrices over sample covariance matrices is that it does not require the first two population moments of the elements of $\bbx$ to be known. This point makes the linear spectral statistics based on sample correlation matrices more practical in applications. By contrast, linear spectral statistics for sample covariances involve unknown moments, and are therefore practically infeasible.

Large dimensional random matrix theory provides us with a powerful tool to establish asymptotic theory for high dimensional sample covariance matrices. \cite{BS2004} contributed to the establishment of asymptotic theory for linear spectral statistics based on high dimensional sample covariance matrices. Meanwhile, there are few results available in the literature for investigating high dimensional sample correlation matrices. \cite{Jiang2004}, among one of the first, established a limiting spectral distribution for sample correlation matrices. \cite{CJ2011} developed some limiting laws of coherence for sample correlation matrices. In addition, both \cite{BPZ2012} and \cite{PY2012} established asymptotic distributions for the extreme eigenvalues of the sample correlation matrices under study. By moving one step further, this paper develops a new central limit theorem for a linear spectral statistic (LSS), which is based on the empirical spectral distribution (ESD) of the sample correlation matrix of $\bbx$. LSS is a general class of statistics in the sense of being able to cover a lot of commonly used statistics. This new CLT is also of independent interest in large dimensional random matrix theory.

In addition to the establishment of a new CLT, we discuss two relevant statistical applications of both the linear spectral statistic of the sample correlation matrix and the resulting asymptotic theory. The first one is an independence test for $p$ random variables included in the vector $\bbx$. A related study is \cite{S2005}, who discussed this kind of independence test for $p$ normal random variables. The second application is to test the equivalence of factor loadings or factors in a factor model. As we discuss in Section 3 below, sample correlation matrices can be used directly for testing purposes without estimating factor loadings and factors first.

The rest of the paper is organized as follows. Section 2 introduces a class of linear spectral statistics. An asymptotic theory is established in Section 3.1 and its applications are established in Section 3.2. The finite sample performance of the proposed test is reported and discussed in Section 4. An empirical application to test independence for household incomes from different cities in China is provided in Section 5. Section 6 concludes the main discussion of this paper. The proofs of the main theory stated in Section 3.1 is given in an appendix. The proofs of some necessary lemmas are provided in Section 8.

\section{Linear Spectral Statistics}

Given a $p$-dimensional random vector $\bbx=(X_1, X_2, \ldots, X_p)^{*}$ with $n$ random samples $\bbx_1, \bbx_2, \ldots, \bbx_n$, where $\bbx_i=(X_{1i}, X_{2i}, \ldots, X_{pi})^{*}$, $i=1, 2, \ldots, n$.
Let $\bbX_n=(\bby_1-\bar{\bby}_1, \bby_2-\bar{\bby}_2, \ldots, \bby_p-\bar{\bby}_p)^{*}$, where $\bby_i=(X_{i1}, X_{i2}, \ldots, X_{in})^{T}$ for $i=1,2,\ldots,p$ and $\bar{\bby}_i=\frac{1}{n}\sum^{n}_{j=1}X_{ij}\bbe$ with $\bbe$ being a $p$-dimensional vector whose elements are all $1$, in which $T$ denotes the transpose of a matrix or a vector.

Consider the sample correlation matrix $\bbB_n=(\rho_{ik})_{p\times p}$ with
\begin{eqnarray*}
\rho_{ik}=\frac{(\bby_i-\bar{\bby}_i)^{*}(\bby_k-\bar{\bby}_k)}{||\bby_i-\bar{\bby}_i||\cdot||\bby_k-\bar{\bby}_k||},
\end{eqnarray*}
where $||\cdot||$ is the usual Euclidean norm. $\bbB_n$ can also be written as
\begin{eqnarray*}
\bbB_n=\bbY^{*}_n\bbY_n=\bbD_n\bbX_n^*\bbX_n\bbD_n,
\end{eqnarray*}
with
\begin{eqnarray*}
\bbY_n=\Big(\frac{\bby_1-\bar{\bby}_1}{||\bby_1-\bar{\bby}_1||}, \frac{\bby_2-\bar{\bby}_2}{||\bby_2-\bar{\bby}_2||}, \ldots, \frac{\bby_p-\bar{\bby}_p}{||\bby_p-\bar{\bby}_p||}\Big)
\end{eqnarray*}
and $\bbD_n={\rm diag}\left(\frac{1}{\|\bby_i-\bar \bby_i\|}\right)_{p\times p}$ is a diagonal matrix.

Let us consider a class of statistics related to the eigenvalues of $\bbB_n$. To this end, define the empirical spectral distribution (ESD) of the sample correlation matrix $\bbB_n$ by $F^{\bbB_n}(x)=\frac{1}{p}\sum^{p}_{i=1}I(\lambda_i\leq x)$, where $\lambda_1\leq\lambda_2\leq\ldots\leq\lambda_p$ are the eigenvalues of $\bbB_n$ and $I(\cdot)$ is an indicator function.

 If $X_1, X_2, \ldots, X_p$ are independent, $F^{\bbB_n}(x)$ converges with probability one to the Marcenko-Pastur (simply called M-P) law $F_{c}(x)$ with $c=\lim_{n\rightarrow\infty}p/n$ (see \cite{Jiang2004}), whose density has an explicit expression of the form
\begin{eqnarray*}
f_c(x)=\left\{\begin{array}{cc}
\frac{1}{2\pi xc}\sqrt{(b-x)(x-a)}, & a\leq x\leq b;\\
0, & otherwise;
\end{array}
\right.
\end{eqnarray*}
and a point mass $1-1/c$ at the origin if $c>1$, where $a=(1-\sqrt{c})^2$ and $b=(1+\sqrt{c})^2$.


Linear spectral statistics of the sample correlation matrix are of the form:
\begin{eqnarray*}
\frac{1}{p}\sum^{p}_{j=1}f(\lambda_j)=\int f(x)dF^{\bbB_n}(x),
\end{eqnarray*}
where $f$ is an analytic function on $[0,\infty)$.

We then consider a normalized and scaled linear spectral statistic of the form:
\begin{eqnarray}\label{0520(2)}
T_{n}(f)=\int f(x)dG_n(x),
\end{eqnarray}
where $G_n(x)=p[F^{\bbB_n}(x)-F_{c_n}(x)]$.

The test statistic $T_n(f)$ is a general statistic in the sense that it covers many classical statistics as special cases. For example,
\begin{enumerate}
\item
Schott's Statistic (\cite{S2005}):
\begin{eqnarray*}
f_1(x)=x^2-x: \ \  T_n(f_1)=tr(\bbB_n^2)-p-p\int (x^2+x)dF_{c_n}(x).
\end{eqnarray*}
\item
The Likelihood Ratio Test Statistic (\cite{M2005}):
\begin{eqnarray*}
f_2(x)=\log(x): \ \ T_n(f_2)=\sum^{p}_{i=1}\log(\lambda_i)-p\int\log(x)dF_{c_n}(x),
\end{eqnarray*}
where $\lambda_i: i=1,2,\ldots,p$ are eigenvalues of $\bbB_n$.
\end{enumerate}

One important tool used in developing an asymptotic distribution for $T_n(f)$ is the Stieltjes transform. The Stieltjes transform $m_G$ for any c.d.f $G$ is defined by
\begin{eqnarray*}
m_{G}(z)=\int\frac{1}{\lambda-z}dG(\lambda), \ \ \Im(z)>0.
\end{eqnarray*}
The Stieltjes transform $m_G(z)$ and the corresponding distribution $G(x)$ satisfy the following relation:
\begin{eqnarray*}
G([x_1,x_2])=\frac{1}{\pi}\lim_{\varepsilon\rightarrow 0}\int^{x_2}_{x_1}\Im\big(m_G(x+i\varepsilon)\big)dx,
\end{eqnarray*}
where $x_1$ and $x_2$ are continuity points of $G$.

Furthermore, the linear spectral statistic can be expressed via the Stieltjes transform of ESD of $\bbB_n$ as follows:
\begin{eqnarray}
\int f(x)dF^{\bbB_n}(x)=-\frac{1}{2\pi i}\oint_{\mathcal{C}} f(z)m_{F^{\bbB_n}}(z)dz,
\end{eqnarray}
where the contour $\mathcal{C}$ contains the support of $F^{\bbB_n}$ with probability one.

\section{Asymptotic Theory and Two Applications}
First, we establish a new central lint theorem for the linear statistic (\ref{0520(2)}) in Theorem \ref{thm1}. Second, we show how to apply the linear statistic and its resulting limiting distribution for an independence test for $p$ random variables and then an equivalence test for factor loadings or factors respectively.

\subsection{Asymptotic Theory}

Before we establish our main theorem, we introduce some notion. Let $\underline{\bbB}_n=\bbY_n\bbY_n^{*}$. The Stieltjes transforms of ESD and LSD for $\underline{\bbB}_n$ are denoted by $\underline{m}_{n}(z)$ and $\underline{m}_{c}(z)$, respectively. Their analogues for $\bbB_n$ are denoted by $m_n(z)$ and $m_c(z)$, respectively. Moreover, $\underline{m}_{c_n}(z)$ and $m_{c_n}(z)$ become $\underline{m}_{c}(z)$ and $m_c(z)$, respectively, when $c$ is replaced by $c_n$. For ease of notation, we denote $m_c(z)$ and $\underline{m}_c(z)$ by $m(z)$ and $\underline{m}(z)$, respectively with omitting the subscript $c$. Moreover, let $\kappa=\lim_{p\rightarrow \infty}\frac{1}{p}\sum_{i=1}^p\frac{\mathbb{E}|X_{i1}-\mathbb{E}X_{i1}|^4}{(\mathbb{E}|X_{i1}-\mathbb{E}X_{i1}|^2)^2}$, and $m^{'}(z)$ denote the first derivative of $m(z)$ with respect to $z$, throughout the rest of this paper.

The following theorem is to establish a joint central limit theorem for the linear spectral statistic of the correlation matrix $\bbB_n$.

\begin{thm}\label{thm1}
Assume that $\{X_{ij}: i=1,2,\ldots,p; j=1,2,\ldots,n\}$ are independent with $\sup_{ 1 \le i\le p}\mathbb{E}|X_{i1}|^4<\infty$. Let $p/n\rightarrow c\in(0,+\infty)$ as $n\rightarrow \infty$. Let $f_1, f_2, \ldots, f_r$ be functions on $\mathbb{R}$ and analytic on an open interval containing
\begin{eqnarray*}
\big[(1-\sqrt{c})^2, (1+\sqrt{c})^2\big].
\end{eqnarray*}

Then, the random vector $\Big(\int f_1(x)dG_n(x), \ldots, \int f_r(x)dG_n(x)\Big)$ converges weakly to a Gaussian vector $(X_{f_1}, \ldots, X_{f_r})$.

When $X_{ij}$ are real random variables, the asymptotic mean is
{\small
\bea
&& E_r\left[X_{f_j}\right]= \frac{\kappa-1}{2\pi i}\oint_{\mathcal{C}}f(z)\frac{c\underline{m}(z)\big(z(1+\underline{m}(z))+1-c\big)}{\big(\big(z(1+\underline{m}(z))-c\big)^2-c\big)\big(z(1+\underline{m}(z))-c\big)}dz
\nonumber\\
&&-\frac{\kappa-|\psi|^2-2}{2\pi i}\oint_{\mathcal{C}}f(z)\frac{cz\underline{m}(z)m^2(z)\big(1+\underline{m}(z)\big)\big(z(1+\underline{m}(z))+1-c\big)}{\big((z(1+\underline{m}(z))-c)^2-c\big)\big(1+cm(z)\big)}dz\nonumber\\
&&-\frac{1}{2\pi i}\oint_{\mathcal{C}}f(z)\frac{c\underline{m}^{'}(z)\Big(z(1+\underline{m}(z))+1-c\Big)}{\underline{m}(z)\big(z+z\underline{m}(z)-c\big)\Big(\big(z(1+\underline{m}(z))-c\big)^2-c\Big)}dz\nonumber\\
&& +\frac{1}{2\pi i} \oint_{\mathcal{C}}f(z)\frac{c\Big(1+z\underline{m}(z)-zm(z)\underline{m}(z)-z^2m(z)\underline{m}^2(z)\Big)\Big(1+\underline{m}(z)\Big)\Big(z(1+\underline{m}(z))+1-c\Big)}{z(1+cm(z))\big(z(1+\underline{m}(z))-c)^2-c\big)}dz\nonumber\\
&&+\frac{1}{2\pi i}\oint_{\mathcal{C}}f(z)\Big(\frac{cm(z)}{z}-czm(z)\underline{m}^{'}(z)\Big)dz
\nonumber
\eea
and the asymptotic covariance function
\bea
&& Cov_r(X_{f_j}, X_{f_r})
\nonumber\\
&& = \  -\frac{1}{2\pi^2}\oint_{\mathcal{C}_1}\oint_{\mathcal{C}_2}f_j(z_1)f_r(z_2)\frac{cm^{'}(z_1)m^{'}(z_2)}{\big(1+c(m(z_1)+m(z_2))+c(c-1)m(z_1)m(z_2)\big)^2}dz_1dz_2
\nonumber\\
&&+\ \frac{\kappa-1}{4\pi^2}\oint_{\mathcal{C}_1}\oint_{\mathcal{C}_2}f_j(z_1)f_r(z_2)\frac{c\underline{m}^{'}(z_1)\underline{m}^{'}(z_2)}{(1+\underline{m}(z_1))^2(1+\underline{m}(z_2))^2}dz_1dz_2
\nonumber\\
&&- \ \frac{\kappa-|\psi|^2-2}{4\pi^2}\oint_{\mathcal{C}_1}\oint_{\mathcal{C}_2}f_j(z_1)f_r(z_2)V(c,m(z_1),m(z_2))dz_1dz_2,
\nonumber
\eea
in which
$\psi=\frac{\mathbb{E}(X_{i1}-\mathbb{E}X_{i1})^2}{\mathbb{E}|X_{i1}-\mathbb{E}X_{i1}|^2}\equiv 1$ under the real case,
\bea
&&V(c,m(z_1),m(z_2)) = c\Big(m(z_1)\underline{m}(z_1)+z_1m(z_1)\underline{m}^{'}(z_1)+z_1m'(z_1)\underline{m}(z_1)\Big)
\nonumber\\
&& \times \ \Big(m(z_2)\underline{m}(z_2)+z_2m(z_2)\underline{m}^{'}(z_2)+z_2m'(z_2)\underline{m}(z_2)\Big)
\nonumber
\eea
for $j,k=1,2,\ldots,r$, and the contour $\oint_{\mathcal{C}}$ is closed and taken in the positive direction in the complex plane, each enclosing the support of $F_{c}(\cdot)$.

When $\left\{X_{ij}\right\}$ are complex variables, assuming that $\psi=\frac{\mathbb{E}(X_{i1}-\mathbb{E}X_{i1})^2}{\mathbb{E}|X_{i1}-\mathbb{E}X_{i1}|^2}$ are the same for i=1,2,...,p, the asymptotic mean is
\bea
&&E_c\left[X_{f_j}\right] = E_r\left[X_{f_j}\right]
\nonumber\\
&& - \ \frac{1}{2\pi i}\oint_{\mathcal{C}}f(z)\Big(\frac{z\underline{m}^{'}(z)}{(1+\underline{m}(z))(z+z\underline{m}(z)-c)}-\frac{c|\psi|^2 m^2(z)}{(1+cm(z))[(1+cm(z))^2-c|\psi|^2m^2(z)]}\Big)
\nonumber\\
&& \times \ \Big(-\frac{c(1+\underline{m}(z))\big(z(1+\underline{m}(z))+1-c\big)}{z\underline{m}(z)\Big(\big(z(1+\underline{m}(z))-c\big)^2-c\Big)}\Big)dz;
\nonumber
\eea
and the asymptotic variance is
\bea
&& Cov_c(X_{f_j}, X_{f_r})=Cov_r(X_{f_j}, X_{f_r})
\nonumber\\
&& + \  \frac{1}{4\pi^2}\oint_{\mathcal{C}_1}\oint_{\mathcal{C}_2}\frac{f_j(z_1)f_r(z_2)cm^{'}(z_1)m^{'}(z_2)dz_1dz_2}{\big(1+c(m(z_1)+m(z_2))+c(c-1)m(z_1)m(z_2)\big)^2}\nonumber\\
&&- \ \frac{|\psi|^2}{4\pi^2}\oint_{\mathcal{C}_1}\oint_{\mathcal{C}_2}\frac{f_j(z_1)f_r(z_2)cm^{'}(z_1)m^{'}(z_2)dz_1dz_2}{[(1+cm(z_1))(1+cm(z_2))-c|\psi|^2m(z_1)m(z_2)]^2}.
\nonumber
\eea
}
\end{thm}

\begin{rmk}
Especially, when $X_{ij}\sim \mathcal{N}(\mu_i,\sigma_i^2)$, i=1,2,...,p; j=1,2,...,n, we have $\kappa\equiv 3$. Although the asymptotic means and variances given above look complicated, they are easy to calculate in practice. In fact, the LSD's $m(z)$ and $\underline{m}(z)$ can be estimated by $\frac{1}{p}tr(\bbB_n-z\bbI_p)^{-1}$ and $\frac{1}{n}tr(\underline{\bbB}_n-z\bbI_n)^{-1}$ respectively. Moreover, asymptotic distributions are still the same after plugging in such estimators due to Slutsky's theorem. The integrals involved in Theorem  3.1 may be calculated by the function `quad' or `dblquad' in MATLAB.
\end{rmk}

\subsection{Two Applications}
In this section, we provide two statistical applications of linear spectral statistics for sample correlation matrices. They are an independence test for high dimensional random vector and an equivalence test for factor loadings or factors in a factor model.

\subsubsection{Independence Test}

For the $p$ random variables grouped in the vector $\bby$, our goal is to test the following hypotheses:
\begin{eqnarray}\label{0520(1)}
&&\mathbb{H}_{10}: X_1, \ldots, X_p \ are \ independent; \ vs \
\mathbb{H}_{1a}: X_1, \ldots, X_p \ are \ dependent.\non
\end{eqnarray}

For this independence test, we make the best use of the linear spectral statistic (\ref{0520(2)}) based on the sample correlation matrix of $\bbx$ with the available $n$ samples $\bbx_1, \bbx_2, \ldots, \bbx_n$. As stated in the last section, under the null hypothesis, the limit spectral distribution of $\bbB_n$ is the M-P law. We use this point to imply independence when applying linear spectral statistics. For simplicity, we choose $f(x)=x^2$ in (\ref{0520(2)}).

\subsubsection{Test for Equivalence of Factor Loadings or Factors}
Since it is difficult to find consistent estimators for unknown factors and loadings, this section proposes to use the proposed linear spectral statistic of the sample correlation matrix for directly testing equivalence for either the factor or the loading without requiring consistent estimators.

Consider the factor model
\begin{eqnarray}\label{3(1)}
X_{it}=\boldsymbol{\lambda}^{T}_i\bbF_t+\varepsilon_{it}, \ \ i=1, 2, \ldots, p;\  t=1, 2, \ldots, n,
\end{eqnarray}
where $\boldsymbol{\lambda}_i$ is an $r$-dimensional factor loading, $\bbF_t$ is the corresponding $r$-dimensional common factor, $\{\varepsilon_{it}: i=1,2,\ldots,p; t=1,2,\ldots,n\}$ are the idiosyncratic components and they are independent for $i=1,2,\ldots,p$ and $t=1,2,\ldots,n$.

One goal is to test
\begin{eqnarray}\label{0812(4)}
\mathbb{H}_{20}: \boldsymbol{\lambda}_1=\boldsymbol{\lambda}_2=\ldots=\boldsymbol{\lambda}_p.
\end{eqnarray}

The proposed statistic is the linear spectral statistic based on the sample correlation matrix $\bbB_n$.
Under $\mathbb{H}_{20}$, model (\ref{3(1)}) reduces to
\begin{eqnarray}\label{3(2)}
X_{it}=\boldsymbol{\lambda}^{T}\bbF_t+\varepsilon_{it}.
\end{eqnarray}

From (\ref{3(2)}), we have
\begin{eqnarray*}
X_{it}-\bar{X}_t=\varepsilon_{it}-\bar{\varepsilon}_t,
\end{eqnarray*}
where $\bar{X}_t=\frac{1}{N}\sum^{N}_{i=1}X_{it}$ and $\bar{\varepsilon}_t=\frac{1}{N}\sum^{N}_{i=1}\varepsilon_{it}$.

In view of this, under the null hypothesis $\mathbb{H}_{20}$, the sample correlation matrix of $\bbx=(X_{i1}, X_{i2}, \ldots, X_{in})^{T}$ is the same as that of $\boldsymbol{\varepsilon}=(\varepsilon_{i1}, \varepsilon_{i2}, \ldots, \varepsilon_{in})^{T}$. Since the components of $\boldsymbol{\varepsilon}$ are independent, the linear spectral statistic (\ref{0520(2)}) follows the asymptotic distribution in Theorem \ref{thm1}. This is the reason why the proposed statistic works in this case.

Another goal is to test
\begin{eqnarray}\label{0813(1)}
\mathbb{H}_{30}: \bbF_1=\bbF_2=\ldots=\bbF_n.
\end{eqnarray}

Similarly, we also propose the linear spectral statistic based on the sample correlation matrix $\bbB_n$.
Under $\mathbb{H}_{30}$, model (\ref{3(1)}) reduces to
\begin{eqnarray}\label{4(2)}
X_{it}=\boldsymbol{\lambda}^{T}_i\bbF+\varepsilon_{it},
\end{eqnarray}

From (\ref{4(2)}), we have
\begin{eqnarray*}
X_{it}-\bar{X}_i=\varepsilon_{it}-\bar{\varepsilon}_i,
\end{eqnarray*}
where $\bar{X}_i=\frac{1}{n}\sum^{n}_{t=1}X_{it}$ and $\bar{\varepsilon}_i=\frac{1}{n}\sum^{n}_{i=1}\varepsilon_{it}$.

Then under the null hypothesis $\mathbb{H}_{30}$, the sample correlation matrix of $\widetilde{\bbx}=(X_{1t}, X_{2t}, \ldots, X_{pt})^{T}$ is the same as that of $\widetilde{\boldsymbol{\varepsilon}}=(\varepsilon_{1t}, \varepsilon_{2t}, \ldots, \varepsilon_{pt})^{T}$. This point makes the proposed statistic (\ref{0520(2)}) applicable and useful in this situation.

\begin{rmk}
We consider a special example of interactive factor model (\ref{3(1)}) of the form:
\begin{eqnarray}\label{0812(1)}
X_{it}=\alpha_i+f_t+\varepsilon_{it}, \ \ i=1,2,\ldots,p; \ t=1,2,\ldots,n,
\end{eqnarray}
where $\alpha_i$ is the specific fixed effects corresponding to section $i$ for $i=1,2,\ldots,n$, $f_t=f(\frac{t}{T})$ is a trend function, $\{\varepsilon_{it}: i=1,2,\ldots,p; t=1,2,\ldots,n\}$ are the idiosyncratic components and they are independent for $i=1,2,\ldots,p$ and $t=1,2,\ldots,n$.

For model (\ref{0812(1)}), we consider the null hypothesis test
\begin{eqnarray}\label{0812(2)}
\mathbb{H}_{40}: \alpha_1=\alpha_2=\cdots=\alpha_p.
\end{eqnarray}
We may propose the same statistic as that for (\ref{0812(4)}).
\end{rmk}

\section{Finite sample analysis}

The finite sample performance of the proposed linear spectral statistic in the two applications are being investigated. We present the empirical sizes and powers of the proposed test.

\subsection{Empirical sizes and powers}

First, we introduce the method of calculating the empirical sizes and powers. Since the asymptotic distribution of the proposed test statistic $R_n$ is a standard normal distribution, it is not difficult to compute the empirical sizes and powers. Let $z_{1-\frac{1}{2}\alpha}$ and $z_{\frac{1}{2}\alpha}$ be the $100(1-\frac{1}{2}\alpha)\%$ and $\frac{1}{2}\alpha$ quantiles of the standard normal distribution. With $K$ replications of the data set simulated under the null hypothesis, we calculate the empirical size as
\begin{equation}
\hat \alpha=\frac{\{\sharp \ of \ R_n^H\geq z_{1-\frac{1}{2}\alpha} or \ R_n^H\leq z_{\frac{1}{2}\alpha}\}}{K},
\end{equation}
where $R_n^H$ represents the value of the test statistic $R_n$ based on the data simulated under the null hypothesis.

In our simulation, we choose $K=1000$ as the number of the replications. The significance level is $\alpha=0.05$. Similarly, the empirical power is calculated as
\begin{equation}
\hat \beta=\frac{\{\sharp \ of \ R_n^A\geq z_{1-\frac{1}{2}\alpha} or \ R_n^A\leq z_{\frac{1}{2}\alpha}\}}{K},
\end{equation}
where $R_n^A$ represents the value of the test statistic $R_n$ based on the data simulated under the alternative hypothesis.

\subsection{Independence Test}

First, we generate the data $\bbx=(X_1, X_2, \ldots, X_p)$ with $n$ random samples $\bbx_1, \bbx_2, \ldots, \bbx_n$ in the following data generating process. Let $\bbx_i=\bbT\bbz_i$, where $\bbz_i=(Z_{1i}, Z_{2i}, \ldots, Z_{pi})^{T}$ with the first $[p/2]$ components $(Z_{1i}, Z_{2i}, \ldots, Z_{[p/2]i})$ being generated from the standard normal distribution and the rest of the components $(Z_{[p/2]+1,i}, Z_{[p/2]+2,i}, \ldots, Z_{pi})$ being generated from Gamma(1,1), in which $[m]\leq m$ denotes the largest integer of $m$. The $p\times p$ deterministic matrix $\bbT$ is generated in the following scenarios:

\begin{enumerate}
\item Independent case: $\bbT=\bbI_p$, where $\bbI_p$ is an identity matrix;
\item Dependent case(1): $\bbT=\bbI_p+\frac{1}{\sqrt{n}}\bbu\bbv^{T}$, where $\bbu$ and $\bbv$ are $p\times 1$ random vectors whose elements are generated from the standard normal distribution;
\item Dependent case(2): $\bbT=\bbI_p+\bbd\bbe^{T}+\bbe\bbd^{T}$, where $\bbd=(0.5, 0, 0, \ldots, 0)^{T}$ is $p\times 1$ vector with the first element being $0.5$ and the rest of the elements being $0$, and $\bbe$ is a $p\times 1$ vector whose elements are all $1$.

\end{enumerate}

The empirical sizes corresponding to the independent case are listed in Table \ref{tb1}. The table shows that, as the pair $(n, p)$ increases jointly, the sizes are close to the true value $0.05$. The empirical powers under the two dependent cases above are presented in Table \ref{tb4} and Table \ref{tb5} respectively. The tendency of the powers going to $1$, as $(n,p)$ increases, illustrates both the finite--sample applicability and the effectiveness of the proposed test statistic.

\begin{table}[htbp]
  \centering
  \caption{Independent test: size(half gamma)}
  \label{tb1}
  {\scriptsize
    \begin{tabular}{cccccc}
   \hline
   n\verb|\|c & 0.2   & 0.4   & 0.6   & 0.8   & 1 \\
   \hline
    20    & 0.0248 & 0.0310 & 0.0376 & 0.0366 & 0.0374 \\
    30    & 0.0360 & 0.0376 & 0.0440 & 0.0400  & 0.0416 \\
    40    & 0.0360 & 0.0424 & 0.0446 & 0.0452 & 0.0436 \\
    50    & 0.0410 & 0.0482 & 0.0484 & 0.0512 & 0.0440 \\
    60    & 0.0428 & 0.0486 & 0.0448 & 0.0482 & 0.0516 \\
    \hline
    \end{tabular}%
    }
  \label{tab:addlabel}%
\end{table}

\begin{table}[htbp]
  \centering
  \caption{Independent test: power(I+$\frac{1}{\sqrt n}u_p v^*_p$)}
  \label{tb4}
  {\scriptsize
    \begin{tabular}{cccccc}
   \hline
  n\verb|\|c& 0.2   & 0.4   & 0.6   & 0.8   & 1.0 \\
   \hline
    10    & 0.1640 & 0.2902 & 0.4704 & 0.6404 & 0.7682 \\
    20    & 0.4092 & 0.7342 & 0.9114 & 0.9816 & 0.9952 \\
    30    & 0.6244 & 0.9384 & 0.9942 & 0.9998 & 1.0000 \\
    40    & 0.8076 & 0.9890 & 0.9994 & 1.0000     & 1.0000 \\
    50    & 0.9022 & 0.9986 & 1.0000     & 1.0000     & 1.0000 \\
    \hline
    \end{tabular}%
    }
  \label{tab:addlabel}%
\end{table}%
\begin{table}[htbp]
  \centering
  \caption{Independent test: power(a=0.5)}
  \label{tb5}
  {\scriptsize
    \begin{tabular}{ccccccccc}
   \hline
      (n,c)\verb|\|d      & 0.1   & 0.2   & 0.3   & 0.4   & 0.5   & 0.6   & 0.7   & 0.8 \\
    \hline
    (20,0.4) & 0.2916 & 0.6368 & 0.6310 & 0.8082 & 0.8930 & 0.9318 & 0.9506 & 0.9534 \\
    (20,0.8) & 0.2416 & 0.3700  & 0.5806 & 0.6662 & 0.7808 & 0.8452 & 0.8692 & 0.9066 \\
    (30,0.4) & 0.3102 & 0.6916 & 0.9326 & 0.9668 & 0.9784 & 0.9884 & 0.9892 & 0.9928 \\
    (30,0.8) & 0.2580 & 0.6384 & 0.7828 & 0.9048 & 0.9444 & 0.9622 & 0.9836 & 0.9902 \\
    (40,0.4) & 0.7000 & 0.8826 & 0.9762 & 0.9874 & 0.9974 & 0.9976 & 0.9988 & 0.9996 \\
    (40,0.8) & 0.4080 & 0.7628 & 0.9284 & 0.9730 & 0.9870 & 0.9944 & 0.9984 & 0.9994 \\
    \hline
    \end{tabular}%
    }
  \label{tab:addlabel}%
\end{table}%

\subsection{Equivalence Tests for Factor Loadings or Factors}

As for the equivalence test (\ref{0812(4)}) for factor loadings, we generate data for factors and idiosyncratic components as follows. The idiosyncratic components $\{\varepsilon_{it}: i=1,2,\ldots,p; t=1,2,\ldots,n\}$ are generated from the standard normal distribution and the factors $\bbF_t$ is $AR(1)$, i.e.
\begin{eqnarray*}
\bbF_t=a\bbF_{t-1}+\boldsymbol{\eta}_t, \ \ t=1,2,\ldots,n,
\end{eqnarray*}
where $a=0.2$ and $\{\boldsymbol{\eta}_t\}$ is generated independently from the standard normal distribution. The initial value $\bbF_0=\textbf{0}$.  The number of factors takes values of $2$ and $3$, respectively, in the simulation.

Factor loadings are generated in the following two scenarios.

\begin{enumerate}
\item DGP(1): $\boldsymbol{\lambda}_i=\boldsymbol{\lambda}$ for $i=1,2,\ldots,p$, where $\boldsymbol{\lambda}$ is generated from the standard normal distribution.
\item DGP(2): $\boldsymbol{\lambda}_i=\boldsymbol{\lambda}$ for $i=1,2,\ldots,[d\cdot p]$, where $d=0.1$; $\boldsymbol{\lambda}_j$ is generated independently from the standard normal distribution for each $j=[d\cdot p], [d\cdot p]+1, \ldots, p$.
\end{enumerate}

For this test, the empirical sizes under DGP(1) are shown in Table \ref{tb2} while the empirical powers under DGP(2) are given in Table \ref{tb6} and Table \ref{tb7}. As $(n,p)$ increases jointly, the empirical sizes tend to the nominal level of $5\%$. The powers show that our proposed test statistic can capture some local alternatives effectively. As $p=30$, there are $3$ different factor loadings under the alternative hypothesis which can be distinguished by the proposed test statistic.

\begin{table}[htbp]
  \centering
  \caption{Factor loading test: size}
  \label{tb2}
  {\scriptsize
    \begin{tabular}{cccccc}
    \hline
     n\verb|\|c& 0.2   & 0.4   & 0.6   & 0.8   & 1 \\
    \hline
    20    & 0.0234 & 0.0320 & 0.0348 & 0.0324 & 0.0346 \\
    30    & 0.0328 & 0.0374 & 0.0376 & 0.0386 & 0.0404 \\
    40    & 0.0338 & 0.0386 & 0.0462 & 0.0444 & 0.0454 \\
    50    & 0.0348 & 0.0440 & 0.0456 & 0.0460 & 0.0424 \\
   \hline
    \end{tabular}%
    }
  \label{tab:addlabel}%
\end{table}

\begin{table}[htbp]
  \centering
  \caption{Factor loading test: power(r=2, different factor loadings are at n-direction)}
  \label{tb6}
  {\scriptsize
    \begin{tabular}{cccccccccc}
   \hline
     (n,p)\verb|\|d   & 0.1   & 0.2   & 0.3   & 0.4   & 0.5   & 0.6   & 0.7   & 0.8   & 0.9 \\
    \hline
    (10,10) & 0.0690  & 0.1096  & 0.1446  & 0.1812  & 0.2070  & 0.2256  & 0.2486  & 0.2526  & 0.2394  \\
    (20,10) & 0.0726  & 0.1100  & 0.1536  & 0.1886  & 0.2180  & 0.2392  & 0.2646  & 0.2682  & 0.2700  \\
    (30,10) & 0.0742  & 0.1134  & 0.1624  & 0.1964  & 0.2214  & 0.2432  & 0.2586  & 0.2634  & 0.2782  \\
    (20,20) & 0.1100  & 0.2070  & 0.3068  & 0.3964  & 0.4616  & 0.5216  & 0.5578  & 0.6092  & 0.6264  \\
    (30,20) & 0.1010  & 0.1830  & 0.2884  & 0.3744  & 0.4464  & 0.4954  & 0.5486  & 0.6062  & 0.6126  \\
    (30,30) & 0.1412  & 0.2624  & 0.4088  & 0.5266  & 0.6172  & 0.7004  & 0.7464  & 0.8050  & 0.8368  \\
    \hline
    \end{tabular}%
    }
  \label{tab:addlabel}%
\end{table}%

\begin{table}[htbp]
  \centering
  \caption{Factor loading test: power(r=3, different factor loadings are at n-direction)}
  \label{tb7}
  {\scriptsize
    \begin{tabular}{cccccccccc}
   \hline
     (n,p)\verb|\|d   & 0.1   & 0.2   & 0.3   & 0.4   & 0.5   & 0.6   & 0.7   & 0.8   & 0.9 \\
    \hline
    (10,10) & 0.0942  & 0.1648  & 0.2074  & 0.2418  & 0.2854  & 0.2956  & 0.2834  & 0.2964  & 0.2920  \\
    (20,10) & 0.1618  & 0.2970  & 0.4078  & 0.4862  & 0.5642  & 0.6044  & 0.6466  & 0.6854  & 0.6826  \\
    (30,10) & 0.2202  & 0.4230  & 0.5646  & 0.6578  & 0.7318  & 0.8028  & 0.8342  & 0.8612  & 0.8766  \\
    (20,20) & 0.1692  & 0.2816  & 0.4252  & 0.5226  & 0.5998  & 0.6518  & 0.7026  & 0.7406  & 0.7438  \\
    (30,20) & 0.2068  & 0.4228  & 0.5774  & 0.7024  & 0.7808  & 0.8478  & 0.8812  & 0.9074  & 0.9348  \\
    (30,30) & 0.1954  & 0.4052  & 0.5770  & 0.6918  & 0.7768  & 0.8372  & 0.8848  & 0.9092  & 0.9320  \\
    \hline
    \end{tabular}%
    }
  \label{tab:addlabel}%
\end{table}%

Similarly, for the equivalence test (\ref{0813(1)}) for factors, the idiosyncratic components are generated in the same way as the test above. The factor loading $\{\boldsymbol{\lambda}_i\}$ is generated independently from the standard normal distribution.

Factors are generated in the following two scenarios.
\begin{enumerate}
\item DGP(3): $\bbF_t=\bbF$ for $t=1,2,\ldots,n$, where $\bbF$ is generated independently from the standard normal distribution.
\item DGP(4): $\bbF_t=\bbF$ for $i=1,2,\ldots,[d\cdot n]$, where $d=0.1$; $\bbF_t$ is generated independently from the standard normal distribution for $t=[d\cdot n], [d\cdot n]+1, \ldots, n$.
\end{enumerate}

The empirical sizes under DGP(3) are shown in Table \ref{tb3} while the empirical powers under DGP(4) are given in Table \ref{tb8} and Table \ref{tb9}. The behaviours of the sizes and powers are similar to those discussed in the factor loading test.

\begin{table}[htbp]
  \centering
  \caption{Factor test: size}
  \label{tb3}
  {\scriptsize
    \begin{tabular}{cccccc}
    \hline
      n\verb|\|c     & 0.2   & 0.4   & 0.6   & 0.8   & 1 \\
    \hline
    20    & 0.0286 & 0.0330 & 0.0348 & 0.0384 & 0.0390 \\
    30    & 0.0322 & 0.0352 & 0.0396 & 0.0398 & 0.0412 \\
    40    & 0.0322 & 0.0362 & 0.0410 & 0.0420 & 0.0414 \\
    50    & 0.0360 & 0.0442 & 0.0462 & 0.0456 & 0.0440 \\
     \hline
    \end{tabular}%
    }
  \label{tab:addlabel}%
\end{table}%

\begin{table}[htbp]
  \centering
  \caption{Factors test: power(r=2, different factors are at n-direction)}
  \label{tb8}
  {\scriptsize
    \begin{tabular}{cccccccccc}
    \hline
      (n,p)\verb|\|d    & 0.1   & 0.2   & 0.3   & 0.4   & 0.5   & 0.6   & 0.7   & 0.8   & 0.9 \\
    \hline
    (10,10) & 0.0696  & 0.1170  & 0.1528  & 0.1822  & 0.1994  & 0.2248  & 0.2272  & 0.2530  & 0.2474  \\
    (20,10) & 0.1146  & 0.2016  & 0.3024  & 0.3684  & 0.4386  & 0.4850  & 0.5316  & 0.5606  & 0.5710  \\
    (30,10) & 0.1582  & 0.2970  & 0.4260  & 0.5338  & 0.6088  & 0.6850  & 0.7192  & 0.7564  & 0.7734  \\
    (20,20) & 0.1024  & 0.2038  & 0.3002  & 0.3918  & 0.4612  & 0.5214  & 0.5548  & 0.5988  & 0.6158  \\
    (30,20) & 0.1354  & 0.2896  & 0.4130  & 0.5492  & 0.6340  & 0.7116  & 0.7574  & 0.8096  & 0.8310  \\
    (30,30) & 0.1358  & 0.2810  & 0.4058  & 0.5304  & 0.6268  & 0.6988  & 0.7594  & 0.8094  & 0.8302  \\
   \hline
    \end{tabular}%
    }
  \label{tab:addlabel}%
\end{table}%

\begin{table}[htbp]
  \centering
  \caption{Factors test: power(r=3, different factors are at n-direction)}
  \label{tb9}
  {\scriptsize
    \begin{tabular}{cccccccccc}
    \hline
     (n,p)\verb|\| d  & 0.1   & 0.2   & 0.3   & 0.4   & 0.5   & 0.6   & 0.7   & 0.8   & 0.9 \\
    \hline
    (10,10) & 0.0996  & 0.1590  & 0.2088  & 0.2566  & 0.2688  & 0.2910  & 0.3004  & 0.2960  & 0.2930  \\
    (20,10) & 0.1606  & 0.2996  & 0.3968  & 0.4984  & 0.5556  & 0.6016  & 0.6298  & 0.6632  & 0.6784  \\
    (30,10) & 0.2272  & 0.4100  & 0.5502  & 0.6568  & 0.7334  & 0.7912  & 0.8298  & 0.8620  & 0.8748  \\
    (20,20) & 0.1554  & 0.2988  & 0.4358  & 0.5252  & 0.5906  & 0.6592  & 0.7010  & 0.7336  & 0.7584  \\
    (30,20) & 0.2138  & 0.4166  & 0.5762  & 0.7024  & 0.7880  & 0.8506  & 0.8826  & 0.9120  & 0.9256  \\
    (30,30) & 0.2074  & 0.4028  & 0.5660  & 0.6960  & 0.7850  & 0.8362  & 0.8842  & 0.9210  & 0.9304  \\
    \hline
    \end{tabular}%
    }
  \label{tab:addlabel}%
\end{table}%

Another equivalence test (\ref{0812(2)}) is also analyzed. The idiosyncratic components $\{\varepsilon_{it}: i=1,2,\ldots,p; t=1,2,\ldots,n\}$ are generated independently from the standard normal distribution, and the trend function $f_t=t/n$.

The specific character $\alpha_i$  for each section $i=1,2,\ldots,p$ is generated in the following two scenarios.
\begin{enumerate}
\item DGP(1): $\alpha_i=\alpha$ with $i=1,2,\ldots,p$ where $\alpha$ is generated from standard normal distribution.
\item DGP(2): $\alpha_i=\alpha$ with $i=1,2,\ldots,[d\cdot p]$ where $d=0.1$; $\alpha_j$ is generated from standard normal distribution independently for each $j=[d\cdot p], [d\cdot p]+1, \ldots, p$.
\end{enumerate}

The empirical sizes and powers are illustrated in Table \ref{tb33} and Table \ref{tb44} respectively. In contrast with the powers in the factor loading test, the powers are relatively lower. It is reasonable because the specific characteristic $\alpha_i$ is not affected by the common factors. In summary, the proposed statistic still works well numerically in this case.

\begin{table}[htbp]
  \centering
  \caption{Specific characteristic test: size}
  \label{tb33}
  {\scriptsize
    \begin{tabular}{cccccccccc}
    \hline
      (n,p)\verb|\|d   & 0.1   & 0.2   & 0.3   & 0.4   & 0.5   & 0.6   & 0.7   & 0.8   & 0.9 \\
    \hline
   (10,10) & 0.0286  &  0.0302 &  0.0292 &   0.0318 &   0.0276 &   0.0348 &   0.0324 &   0.0344 &   0.0328\\
   (20,10) & 0.0364  &  0.0334  &  0.0350 &   0.0392 &   0.0366 &   0.0400 &   0.0360 &   0.0350 &   0.0334\\
   (30,10) & 0.0360   & 0.0424  &  0.0334  &  0.0338 &   0.0386 &   0.0400 &   0.0398 &   0.0360 &   0.0360\\
   (20,20) & 0.0372    &0.0344  &  0.0392  &  0.0388 &   0.0402 &   0.0378 &   0.0386 &   0.0414 &   0.0392\\
   (30,20) & 0.0390   & 0.0408   & 0.0388  &  0.0356 &   0.0432 &   0.0418 &   0.0418 &   0.0390 &   0.0382\\
   (30,30) & 0.0440   & 0.0420 &   0.0434  &  0.0412 &   0.0432 &   0.0396 &   0.0434 &   0.0442 &   0.0436\\
    \hline
    \end{tabular}%
    }
  \label{tab:addlabel}%
\end{table}%

\begin{table}[htbp]
  \centering
  \caption{Specific characteristic test: power}
  \label{tb44}
  {\scriptsize
    \begin{tabular}{cccccccccc}
    \hline
      (n,p)\verb|\|d   & 0.1   & 0.2   & 0.3   & 0.4   & 0.5   & 0.6   & 0.7   & 0.8   & 0.9 \\
    \hline
    (10,10) & 0.0560 &   0.0892 &   0.1298  &  0.1644 &   0.1998 &   0.2292 &   0.2544  &  0.2726 &   0.2822\\
    (20,10) & 0.0638 &   0.0940 &   0.1266  &  0.1670 &   0.1970 &   0.2256 &   0.2368  &  0.2664 &   0.2562\\
    (30,10) & 0.0532 &   0.0864 &   0.1158  &  0.1572 &   0.1768 &   0.2058 &   0.2150  &  0.2480 &   0.2392\\
    (20,20) & 0.0758 &   0.1534 &   0.2258  &  0.3076 &   0.3756 &   0.4428 &   0.5078  &  0.5430 &   0.5826\\
    (30,20) & 0.0644 &   0.1434 &   0.2056  &  0.2738 &   0.3404 &   0.4106 &   0.4608  &  0.5094 &   0.5396\\
    (30,30) & 0.0912 &   0.1852 &   0.2924  &  0.3946 &   0.4972 &   0.5766 &   0.6544  &  0.7102 &   0.7638\\
    \hline
    \end{tabular}%
    }
  \label{tab:addlabel}%
\end{table}%

\section{Empirical Application}
In this section, we analyze the relationship of the household incomes among different cities for rural China. The main goal is to test whether they are independent or not.

The data set is drawn from the `Rural Household Income and Expenditure Survey' conducted by the State Statistics Bureau of China (SSB) and the Chinese Academy of Social Science (CASS). The data set was collected in 1995 and provides useful information about 7998 households in rural areas of 19 Chinese provinces.

In this study, we focus on testing independence of the household incomes among different cities. After deleting observations with missing or implausible values of the household income variables, a sample of $96$ households is retained for $69$ different cities.

The proposed linear spectral statistic is applied to this independence test. Different number of cities and various number of households are considered. The p-values of the proposed test are reported in Table \ref{tb55}. The $p$-values decrease as the number of cities increases. This phenomenon makes sense since the possibility of the dependence becomes larger as the number of cities becomes bigger. Since the $p$-values are all greater than $0.01$, we conclude that the household incomes from different cities are independent.

\begin{table}[h]
{\small \caption{P-values of independence test for household incomes from different cities}\label{tb55}}
\begin{center}
{\footnotesize
\begin{tabular}{cccccccccccccccccc}
&\multicolumn{7}{c}{}\\
\hline
(p,n) & (5, 10)  &  (15, 20)  &  (40, 50)   &  (50, 60)   &  (60, 70)  &  (69,80)  &  (69, 96) \\
\hline
$p-values$        & 0.5260 & 0.4430 & 0.5620 & 0.5290 & 0.0890 &	0.0680 & 0.0540 \\
\hline
\end{tabular}}
\\ \end{center}\medskip
\end{table}

\section{Conclusions}
In this paper, we have established a new central limit theorem for a linear spectral statistic of sample correlation matrices for the case where the dimensionality $p$ and the sample size $n$ are comparable. Two useful statistical applications are considered. The first one is an independence test for $p$ random variables while the second one is an equivalence test in factor models. The advantage of using the linear spectral statistic based on sample correlation matrices over sample covariance matrices is that we do not require the knowledge of the first two moments or the underlying distribution of the $p$ random variables under investigation. The finite sample performance of the proposed test is evaluated. An empirical application to test cross-section independence for the household income in different cities of China is discussed.

{\small

\section{Appendix: Proof of the main theorem}

We start by listing some necessary lemmas.
\subsection{Lemmas}
\begin{lem}[\cite{Jiang2004}; \cite{HW2010}]\label{lem1}
Suppose $p/n\rightarrow c\in(0, +\infty)$. If $\mathbb{E}|X_{11}|^{4}<\infty$ and $\mathbb{E}X_{11}=0$, then $\lambda_{\max}(\bbB_n)\stackrel{a.s.}{\rightarrow}(1+\sqrt{c})^2$ and $\lambda_{\min}(\bbB_n)\stackrel{a.s.}{\rightarrow}(1-\sqrt{c})^2$.
\end{lem}

\begin{lem}[Corollary 7.38 of \cite{HJ1999}]\label{lem0525}
Let A and B be two complex $p \times n$ matrices. Define r=$\min\left\{p,n\right\}$. If $\sigma_1 \ge \sigma_2 \ge ...\ge \sigma_r$ are the first r largest eigenvalues of $A^{*}A$ and $\lambda_1 \ge \lambda_2 \ge...\ge \lambda_r$ are the first r largest eigenvalues of $B^{*}B$, then
$$\max_{1 \le i \le r}|\sqrt{\sigma_i}-\sqrt{\lambda_i}|\le \|A-B\|,$$
where $\|A-B\|$ denotes the largest eigenvalue of $(A-B)^{*}(A-B).$
\end{lem}
\begin{lem}[\cite{B1973}]\label{lem0527}
Let $\left\{X_k\right\}$ be a complex martingale difference sequence with respect to the increasing $\sigma-field \left\{\mathcal F_k\right\}$. Then for q$> 1$,
$$\mathbb{E}|\sum X_k|^q \le K_q \left(\mathbb{E}\left(\sum \mathbb{E}_{k-1}|X_k|^2\right)^{q/2}+\mathbb{E}\left[\sum|X_k|^q\right]\right).
$$
\end{lem}
\begin{lem}[Theorem 35.12 of \cite{B1995}]\label{lem0521}
Suppose for each $n$, $Y_{n1}, Y_{n2}, \ldots, Y_{nr_n}$ is a real martingale difference sequence with respect to the increasing $\sigma$-field $\{\mathcal{F}_{nj}\}$ having second moments. If as $n\rightarrow\infty$,
\begin{eqnarray}
\sum^{r_n}_{j=1}\mathbb{E}(Y_{nj}^2|\mathcal{F}_{n,j-1})\stackrel{i.p.}{\rightarrow}\sigma^2,
\end{eqnarray}
where $\sigma^2$ is positive constant, and for each $\varepsilon>0$,
\begin{eqnarray}
\sum^{r_n}_{j=1}\mathbb{E}(Y_{nj}^2I_{|Y_{nj}|\geq\varepsilon})\rightarrow 0,
\end{eqnarray}
then $\sum^{r_n}_{j=1}Y_{nj}\stackrel{D}{\rightarrow}\mathcal{N}(0,\sigma^2)$.

\end{lem}

The proofs of Lemmas 5-7 below are given in the supplementary.

\begin{lem}\label{lem han}
Suppose that $\left\{X_i\right\}_{i=1}^n$ are i.i.d. random variables with $EX_{1}=0$ and $E|X_{1}|^2 =1$. Let $\bby=(X_1 ,..., X_n)^{T}$ and $\bar{\bby}=\frac{\sum_{i=1}^{n}X_{i}}{n}\bbe$,  where $\bbe=(1,1,...,1)^T $ is an n-dimensional vector.  Assuming that $\bbA$ is a deterministic complex matrix, then for any given $ q\ge 2$ , there is a positive constant $K_q$ depending on q such that
\begin{equation}
E\left| \pmb\alpha^{*}\bbA\pmb\alpha-\frac{1}{n}tr\bbA \right|^q \leq K_q  \left\{n^{-q} (v_{2q}tr(\bbA\bbA^*)^{q}+(v_4 tr(\bbA\bbA^*))^{q/2})+\mathbb{P}(B^{c}_{n}(\ep))\|\bbA\|^q\right\},
\end{equation}
where $B_n(\ep) =\left\{\bby: |\frac{ \|\bby-\bar{\bby}\|^2}{n}-1| \le \ep \right\}$ and $\pmb\alpha=\frac{(\bby-\bar{\bby})^T}{\|\bby-\bar{\bby}\|}$, in which $\ep >0$ is a constant.

\end{lem}

\begin{rmk}
Note that $\mathbb{P}(B^{c}_{n}(\ep))=O(n^{-q/2}v^{q/2}_4+ n^{-q+1}v_{2q})$. 
If  $\|\bbA\|\leq K$ and $|X_i|\le \sqrt n \delta_n$, we have
\begin{equation}\label{1019.2}
\mathbb{E}\left|\pmb\alpha^{*}\bbA\pmb\alpha-\frac{1}{n}tr\bbA \right|^q  \leq K_{q}n^{-1}\delta^{2q-4}_{n}.
\end{equation}
\end{rmk}
\begin{rmk}\label{remark1}
Similar to Lemma \ref{lem han}, one can prove that under the same conditions of Lemma \ref{lem han} (replacing $\pmb\alpha*$ by $\pmb\alpha^{T}$), we have
\begin{equation}
\mathbb{E}\left|\pmb\alpha^{T}\bbA\pmb\alpha-\frac{\mathbb{E}X_{1}^2}{n}tr\bbA \right|^q \leq K_q  \left\{n^{-q} (v_{2q}tr(\bbA\bbA^*)^{q}+(v_4 tr(\bbA\bbA^*))^{q/2})+\mathbb{P}(B^{c}_{n}(\ep))\|\bbA\|^q\right\}.
\end{equation}
\end{rmk}

 \begin{lem}\label{lem yang}
In addition to the condition of Lemma \ref{lem han}, if $\mathbb{E}|X_1|^4<\infty$, $\|\bbA\|\le K$ and $\|\bbB\|\le K$, then
 \begin{eqnarray}
  &&\mathbb{E}(\pmb\alpha^{*}\bbA\pmb\alpha-\frac{1}{n}tr\bbA)(\pmb\alpha^{*}\bbB\pmb\alpha-\frac{1}{n}tr\bbB) =  \sum_{i=1}^{n}\frac{1}{n^2}(\mathbb{E}|X_{1}|^{4}-|\mathbb{E}(X_{1}^{2})|^2 -2)\bbA_{ii}\bbB_{ii}\non
    &&+ \ \frac{|\mathbb{E}X_{1}^{2}|^2}{n^2}tr(\bbA\bbB^{T}) + \frac{1}{n^2}tr(\bbA\bbB)+\frac{1-\mathbb{E}|X_{1}|^{4}}{n^{3}}tr \bbA tr \bbB+o\left(\frac{1}{n}\right).\nonumber
  \end{eqnarray}

    \end{lem}

    In the sequel, we assume that $\{X_{ij}\}$ satisfies
\begin{eqnarray}\label{*}
&&\\
|X_{ij}|<\delta_n\sqrt{n}, \  \mathbb{E}X_{ij}=0, \ \mathbb{E}|X_{ij}|^2=1, \ \mathbb{E}|X_{ij}|^4<\infty, \ \mbox{and} \ \ \ \kappa=\lim_{p \rightarrow \infty}\frac{1}{p}\sum_{i=1}^p\mathbb{E}|X_{i1}|^4.\nonumber
\end{eqnarray}

\begin{lem}\label{lem eig}
For any $l \in \bbN^+$, $\mu_1 >(1+\sqrt{c})^{2}$ and $0<\mu_2 <\mathbf{I}_{(0,1)}(c)(1-\sqrt{c})^2$, under condition (\ref{*}), we have
\begin{eqnarray}\label{gg1}
P(\|\bbB_n\|\geq \mu_1)=o(n^{-l})
\end{eqnarray}
and
\begin{eqnarray}\label{gg2}
P(\lambda_{\min}^{\bbB_n}\leq \mu_2)=o(n^{-l}).
\end{eqnarray}
\end{lem}

\subsection{Proof of Theorem \ref{thm1}}
The overall strategy of our proof is similar to that in \cite{BS2004}. Since many tools proposed in \cite{BS2004} can not be utilized for the sample correlation matrix case, we therefore develop a number of new techniques for the proof of Theorem \ref{thm1}. Among them, to apply the Cauchy integral formula in (\ref{1019.1}) below and prove tightness, we develop Lemma \ref{lem eig} to make sure that the extreme eigenvalues of $\bbB_n$ are highly concentrated around two edges of the support. To convert random quadratic forms into the corresponding traces, we establish a moment inequality for random quadratic forms in Lemma \ref{lem han}. Lemma \ref{lem yang} also establishes a precise estimator for the expectation of the product of two random quadratic forms before we may apply central limit theorems for martingale differences.  Moreover, we find out the limit of the quadratic form $\mathbf{1}^T\mathbb{E}(\underline{\bbB}_n-z\bbI)^{-1}\mathbf{1}/n$ is independent of $\underline{m}(z)$, which is quite different from what may be obtained in the case of covariance matrices (here all entries of the vector $\mathbf{1}$ are one). One can refer to Lemma 8 in the supplementary  for detail.

By the Cauchy integral formula, we have
\begin{equation}\label{1019.1}
\int f(x)dG(x)=-\frac{1}{2\pi i}\int_\mathcal{C}f(z)m_G(z)dz
\end{equation}
valid for any c.d.f $G$ and any function $f$ analytic on an open set containing the support of $G$. In our case, $G(x):=G_n(x)=p(F^{\bbB_n}(x)-F_{c_n}(x))$.

Note that the support of $G_n(x)$ is random. Fortunately, it is well known that the extreme eigenvalues of $\bbB_n$ are highly concentrated around two edges of the support of the limiting M-P law $F_c(x)$ (see Lemma \ref{lem eig}). Then the contour $\mathcal{C}$ can be appropriately chosen. Moreover,  as in \cite{BS2004}, by Lemma \ref{lem eig}, we can replace the process $\{M_n(z),\mathcal{C}\}$ by a slightly modified process $\{\widehat{M}_n(z), \mathcal{C}\}$. Below we present the definitions of the contour $\mathcal{C}$ and the modified process $\widehat{M}_n(z)$. Let $x_r$ be any number greater than $(1+\sqrt{c})^2$. Let $x_l$ be any negative number if $c\geq 1$. Otherwise we choose $x_l\in(0,(1-\sqrt{c})^2)$. Now let $\mathcal{C}_u=\{x+iv_0:x\in[x_l,x_r]\}$.

Then we define $\mathcal{C}^+\equiv\{x_l+iv:v\in[0,v_0]\}\cup\mathcal{C}_u\cup\{x_r+iv: v\in[0,v_0]\}$, and $\mathcal{C}=\mathcal{C}^+\cup \overline{\mathcal{C}^+}$. Now we define the subsets $\mathcal{C}_n$ of $\mathcal{C}$ on which $M_n(\cdot)$ equals to $\widehat{M}_n(\cdot)$. Let $\{\varepsilon_n\}$ be a sequence decreasing to zero satisfying for some $\alpha\in(0,1)$, $\varepsilon_n\geq n^{-\alpha}$.

Now we set
\begin{equation*}
\mathcal{C}_l=\left\{
\begin{array}{ccc}
\{x_l+iv:v\in[n^{-1}\varepsilon_n,v_0]\}~~~~ \mathrm{if}~x_l>0,\\
\{x_l+iv:v\in[0,v_0]\} ~~~~~~~~~~\mathrm{if}~x_l<0,
\end{array}
\right.
\end{equation*}
and $\mathcal{C}_r=\{x_r+iv: v\in[n^{-1}\varepsilon,v_0]\}$.

Then we define $\mathcal{C}_n=\mathcal{C}_l\cup\mathcal{C}_u\cup\mathcal{C}_r$. The process $\widehat{M}_n(z)$ is defined as
\begin{equation*}
\widehat{M}_n(z)=\left\{
\begin{array}{ccc}
M_n(z)&~\mathrm{for}~ z\in \mathcal{C}_n,\\
M_n(x_r+in^{-1}\varepsilon_n)&~~~~~~~~~\mathrm{for}~ x=x_r,v\in[0,n^{-1}\varepsilon_n],\\
M_n(x_l+in^{-1}\varepsilon_n)&~~~~~~~~~\mathrm{for}~
x=x_l,v\in[0,n^{-1}\varepsilon_n].
\end{array}\right.
\end{equation*}

To prove Theorem \ref{thm1}, as in \cite{BS2004}, it suffices to prove the CLT for $\widehat{M}_n(z)$ with $z\in \mathcal{C}$. We state the result in the following proposition and then prove it.
\begin{prop}\label{prop1}
Under the conditions of Theorem \ref{thm1}, $\{\widehat{M}_n(\cdot)\}$ forms a tight sequence on $\mathcal{C}^+$. And  $\{\widehat{M}_n(\cdot)\}$ converges weakly to a two-dimensional Gaussian process $\{M(\cdot)\}$ satisfying for $z\in\mathcal{C}^+$. Under the real random variable case,
\begin{eqnarray}\label{yyy1}
\mathbb{E}M(z)&=&-\Big(\kappa-1\Big)\frac{c\underline{m}(z)\big(z(1+\underline{m}(z))+1-c\big)}{\big(\big(z(1+\underline{m}(z))-c\big)^2-c\big)\big(z(1+\underline{m}(z))-c\big)}\non
&&+\Big(\kappa-|\psi|^2-2\Big)\frac{cz\underline{m}(z)m^2(z)\big(1+\underline{m}(z)\big)\big(z(1+\underline{m}(z))+1-c\big)}{\big((z(1+\underline{m}(z))-c)^2-c\big)\big(1+cm(z)\big)}\non
&&+\frac{c\underline{m}^{'}(z)\Big(z(1+\underline{m}(z))+1-c\Big)}{\underline{m}(z)\big(z+z\underline{m}(z)-c\big)\Big(\big(z(1+\underline{m}(z))-c\big)^2-c\Big)}\non
&&-\frac{c\Big(1+z\underline{m}(z)-zm(z)\underline{m}(z)-z^2m(z)\underline{m}^2(z)\Big)\Big(1+\underline{m}(z)\Big)\Big(z(1+\underline{m}(z))+1-c\Big)}{z(1+cm(z))\big(z(1+\underline{m}(z))-c)^2-c\big)}\non
&&-\frac{cm(z)}{z}+czm(z)\underline{m}^{'}(z)\non
\end{eqnarray}
and for $z_i,z_j\in\mathcal{C}$
\begin{eqnarray}\label{yyy2}
\mathrm{Cov}(M(z_i),M(z_j))&=&2\frac{cm^{'}(z_1)m^{'}(z_2)}{\big(1+c(m(z_1)+m(z_2))+(c^2-c)m(z_1)m(z_2)\big)^2}\non
&&-\Big(\kappa-1\Big)\frac{c\underline{m}^{'}(z_1)\underline{m}^{'}(z_2)}{(1+\underline{m}(z_1))^2(1+\underline{m}(z_2))^2}\non
&&+\Big(\kappa-|\psi|^2-2\Big)V(c,m(z_1),m(z_2)),\non
\end{eqnarray}
where $V(c,m(z_1),m(z_2))$ is defined in Theorem 3.1.

When $\{X_{ij}\}$ are complex variables, assuming that $\psi=\frac{\mathbb{E}(X_{i1}-\mathbb{E}X_{i1})^2}{\mathbb{E}|X_{i1}-\mathbb{E}X_{i1}|^2}$ are the same for i=1,2,...,p, the asymptotic mean is
\begin{eqnarray}
&&(\ref{yyy1})+\Big(\frac{z\underline{m}^{'}(z)}{(1+\underline{m}(z))(z+z\underline{m}(z)-c)}-\frac{c|\psi|^2m^2(z)}{(1+cm(z))[(1+cm(z))^2-c|\psi|^2m^2(z)]}\Big)\non
&&\ \ \ \ \ \ \ \ \ \ \ \ \ \ \ \ \ \ \ \ \cdot\Big(-\frac{c(1+\underline{m}(z))\big(z(1+\underline{m}(z))+1-c\big)}{z\underline{m}(z)\Big(\big(z(1+\underline{m}(z))-c\big)^2-c\Big)}\Big);\non
\end{eqnarray}
and the asymptotic variance is
\begin{eqnarray}
(\ref{yyy2})&-&\frac{cm^{'}(z_1)m^{'}(z_2)}{\big(1+c(m(z_1)+m(z_2))+c(c-1)m(z_1)m(z_2)\big)^2}\\
&+&|\psi|^2\frac{cm^{'}(z_1)m^{'}(z_2)}{[(1+cm(z_1))(1+cm(z_2))-c|\psi|^2m(z_1)m(z_2)]^2}\nonumber.
\end{eqnarray}
\end{prop}

By the discussions in \cite{BS2004}, we see that Theorem \ref{thm1} holds if Proposition \ref{prop1} is proved. Thus the rest of the work will be devoted to the proof of Proposition \ref{prop1}.

Before proving Proposition \ref{prop1}, we need to truncate the elements of $\bbX_n$ as follows.

\subsubsection{Truncation, Centralization and Rescaling}

By the same method as that in page 559 of \cite{BS2004}, we can choose a positive sequence of $\{\delta_n\}$ such that
$$\delta_n\rightarrow0, \ \delta_nn^{1/4}\rightarrow\infty,\  \delta_n^{-4}EX_{11}^4I(|X_{11}|\geq\delta_n\sqrt{n})\rightarrow 0.$$

 Let $\hat\bbB_n=\hat\bbD_n\hat\bbX_n\hat\bbX_n^{*}\hat\bbD_n$, where $\hat\bbX_n$ is $p\times n$ matrix having $(i,j)th$ entry $\hat{X}_{ij}-\frac{1}{n}\sum_{k=1}^{n}\hat{X}_{ik}$, $\hat{X}_{ij}=X_{ij}I_{\{|X_{ij}|<\delta_n\sqrt{n}\}}$ and $\hat\bbD_n$ is $\bbD_n$ with $\bbX_n$ replaced by $\hat\bbX_n$. We then have
\begin{eqnarray*}
P(\bbB_n\neq\hat\bbB_n)\leq np\cdot P(|X_{11}|\geq\delta_n\sqrt{n})\leq K\delta_n^{-4}\int_{\{|X_{11}|\geq\delta_n\sqrt{n}\}}|X_{11}|^4=o(1).
\end{eqnarray*}

Define $\widetilde{\bbB}_n=\frac{1}{n}\widetilde{\bbD}_n\widetilde{\bbX}_n\widetilde{\bbX}_n^{*}\widetilde{\bbD}_n$, where $\widetilde{\bbX}_n$ is $p\times n$ matrix having $(i,j)th$ entry $\widetilde{X}_{ij}-\frac{1}{n}\sum_{k=1}^{n}\tilde{X}_{ik}$, $\widetilde{X}_{ij}=(\hat X_{ij}-\mathbb{E}\hat{X}_{ij})/\sigma_n$ with $\sigma_n^2=\mathbb{E}|\hat X_{ij}-\mathbb{E}\hat X_{ij}|^2$; and $\widetilde{\bbD}_n$ is $\bbD_n$ with $X_{ij}$ replaced by $\widetilde{X}_{ij}$. Throughout this paper, we use M and K to denote a constant which can represent different constants at difference appearance.

From Lemma \ref{lem1} and \cite{YBK1988}, we see that
\begin{eqnarray*}
|\lim\sup_{n}\lambda_{\max}^{\hat\bbB_n}|\leq M(1+\sqrt{c})^2,\ \
|\lim\sup_{n}\lambda_{\max}^{\widetilde{\bbB}_n}|\leq M(1+\sqrt{c})^2.
\end{eqnarray*}

Let $\hat G_n(x)$ and $\widetilde{G}_n(x)$ be $G_n(x)$ with $\bbB_n$ replaced by $\hat\bbB_n$ and $\widetilde{\bbB}_n$ respectively. Then for each $j=1,2,\ldots,k$,
\begin{eqnarray*}
&&|\int f_j(x)d\hat G_n(x)-\int f_j(x)d\widetilde{G}_n(x)|\leq M\sum^{p}_{k=1}|\lambda_k^{\hat\bbB_n}-\lambda_k^{\widetilde{\bbB}_n}|\non
&\leq&M\big(\sum^{p}_{k=1}(\sqrt{\lambda_k^{\hat\bbB_n}}-\sqrt{\lambda}_k^{\widetilde{\bbB}_n})^2\big)^{1/2}
\big(\sum^{p}_{k=1}(\lambda_k^{\hat\bbB_n}+\lambda_k^{\widetilde{\bbB}_n})\big)^{1/2}.
\end{eqnarray*}

Moreover, similar to \cite{BS2004} (page 560), we have
\begin{eqnarray*}
&&\sum^{p}_{k=1}(\sqrt{\lambda_k^{\hat\bbB_n}}-\sqrt{\lambda_k^{\widetilde{\bbB}_n}})^2\non
&\leq&M\big(\frac{1}{n}tr(\hat\bbD_n\hat\bbX_n-\widetilde{\bbD}_n\widetilde{\bbX}_n)(\hat\bbD_n\hat\bbX_n-\widetilde{\bbD}_n\widetilde{\bbX}_n)^{*}\big)^{1/2}\big(p(\lambda_{\max}^{\hat\bbB_n}+\lambda_{\max}^{\widetilde{\bbB}_n})\big)^{1/2},
\end{eqnarray*}
where $K_j$ is a bound in $|f_j^{'}(z)|$.

Meanwhile, we have
\begin{eqnarray}\label{1002.1}
&&|\sigma_n^2-1|\leq M\int_{\{|X_{11}|\geq 2\delta_n\sqrt{n}\}}|X_{11}|^2\\
&\leq& 2\delta_n^{-2}n^{-1}\int_{\{|X_{11}|\geq\delta_n\sqrt{n}\}}|X_{11}|^4=o(\delta_n^2n^{-1})\nonumber
\end{eqnarray}
and $|\mathbb{E}\hat X_{11}|\leq 2\big|\int_{\{|X_{11}|\geq\delta_n\sqrt{n}\}}X_{11}\big|=o(\delta_nn^{-3/2})$.

We therefore obtain
\begin{eqnarray}\label{f1}\\
&& \big(\frac{1}{n}tr\big[(\hat\bbD_n\hat\bbX_n-\widetilde{\bbD}_n\widetilde{\bbX}_n)(\hat\bbD_n\hat\bbX_n-\widetilde{\bbD}_n\widetilde{\bbX}_n)^{*}\big]\big)^{1/2}
\nonumber\\
&= & \big(\frac{1}{n}tr\big[\big(\hat\bbD_n(\hat\bbX_n-\widetilde{\bbX}_n)+(\hat\bbD_n-\widetilde{\bbD}_n)\widetilde{\bbX}_n\big)
\big(\hat\bbD_n(\hat\bbX_n-\widetilde{\bbX}_n)+(\hat\bbD_n-\widetilde{\bbD}_n)\widetilde{\bbX}_n\big)^{*}\big]\big)^{1/2}
\nonumber\\
&= & \big(\frac{1}{n}tr\big[\big(\hat\bbD_n(\hat\bbX_n-\widetilde{\bbX}_n)(\hat\bbX_n-\widetilde{\bbX}^*_n)\hat\bbD_n+\hat\bbD_n(\hat\bbX_n-\widetilde{\bbX}_n)\widetilde{\bbX}^*_n(\hat\bbD_n-\widetilde{\bbD}_n)
\nonumber\\
&+ & (\hat\bbD_n-\widetilde{\bbD}_n)\widetilde{\bbX}_n(\hat\bbX_n-\widetilde{\bbX}^*_n)\hat\bbD_n+(\hat\bbD_n-\widetilde{\bbD}_n)\widetilde{\bbX}_n\widetilde{\bbX}^*_n(\hat\bbD_n-\widetilde{\bbD}_n)\big]\big)^{1/2}.
\nonumber
\end{eqnarray}

For the first term on the right hand side above, under the condition of Theorem 1, 
By (\ref{1002.1}), we can see
\begin{eqnarray}
&&\frac{1}{n}\mathbb{E}tr\big[\hat\bbD_n(\hat\bbX_n-\widetilde{\bbX}_n)(\hat\bbX_n-\widetilde{\bbX}_n)^{*}\hat\bbD_n^{*}\big]\\
&=&\frac{1}{n}\sum_{i,j}\mathbb{E}\frac{[\hat{X}_{ij}-\tilde{X}_{ij}-\frac{1}{n}\sum_{k=1}^{n}(\hat{X}_{ik}-\tilde{X}_{ik})]^2}{\sum_{k=1}^{n}(\hat{X}_{ik}-\frac{1}{n}\sum_{l=1}^{n}\hat{X}_{il})^2}
\nonumber\\
&=&\frac{(1-\frac{1}{\sigma_n})^2}{n}\sum_{i,j}\mathbb{E}\frac{(\hat{X}_{ij}-\frac{1}{n}\sum_{l=1}^{n}\hat{X}_{il})^2}{\sum_{k=1}^{n}(\hat{X}_{ik}-\frac{1}{n}\sum_{l=1}^{n}\hat{X}_{il})^2} = \frac{p}{n}(1-\frac{1}{\sigma_n})^2=o(\delta_n^4n^{-2}).
\nonumber
\end{eqnarray}

The remaining terms of (\ref{f1}) can be similarly verified to have an order of o(1/n) and so (\ref{f1}) =$o(n^{-1/2})$. In view of above, we obtain $\int f_j(x)dG_n(x)=\int f_j(x)d\widetilde{G}_n(x)+o_p(1)$. Since $\mathbb{E}|\tilde X_{ij}|^4=\mathbb{E}|X_{ij}|^4+O(n^{-1})$, it will not affect $\kappa=\lim_{p \rightarrow \infty}\frac{1}{p}\sum_{i=1}^p\mathbb{E}|X_{i1}|^4$. To simplify notation we below still use $X_{ij}$ instead of $\tilde{X}_{ij}$ and $\{X_{ij}$\} satisfy (\ref{*}).

\subsubsection{Convergence of $M_n(z)$}
Let $\underline{\boldsymbol{\bbB}}_n=\bbY_n\bbY_n^{*}$. The Stieltjes transforms of ESD and LSD for $\underline{\boldsymbol{\bbB}}_n$ are denoted by $\underline{m}_{n}(z)$ and $\underline{m}_{c}(z)$ respectively. Their analogues for $\boldsymbol{\bbB}_n$ are denoted by $m_n(z)$ and $m_c(z)$ respectively. Moreover, $\underline{m}_{c_n}(z)$ and $m_{c_n}(z)$ are $\underline{m}_{c}(z)$ and $m_c(z)$ respectively with $c$ replaced by $c_n$. For ease of notation, we also denote $m_c(z)$ and $\underline{m}_c(z)$ by $m(z)$ and $\underline{m}(z)$ respectively with omitting the subscript $c$.

Since $M_n(z)=p\big[m_{n}(z)-m_{c_n}(z)\big]=n\big[\underline{m}_{n}(z)-\underline{m}_{c_n}(z)\big]$, we write for $z\in \mathcal{C}_n$, $M_n(z)=M_n^{(1)}(z)+M_n^{(2)}(z)$, where $M_n^{(1)}(z)=n\big[\underline{m}_{n}(z)-\mathbb{E}\underline{m}_{n}(z)\big]$ and $M^{(2)}_n(z)=n\big[\mathbb{E}\underline{m}_{n}(z)-\underline{m}_{c_n}(z)\big]$.

Following the steps in \cite{BS2004}, it suffices to show the following four statements:
\begin{enumerate}
\item \emph{Finite-dimensional convergence of $M_n^{(1)}(z)$ in distribution on $\mathcal{C}_n$};

\item \emph{$M_n^{(1)}(z)$ is tight on $\mathcal{C}_n$};

\item \emph{$\{M_n^{(2)}(z)\}$ for $z\in \mathcal{C}_n$ is bounded and equicontinuous};

\item \emph{$M_n^{(2)}(z)$ converges to a constant and find its limit}.
\end{enumerate}
\textbf{Step 1:}

First, we introduce some notations. In the following proof, we assume that $v=\Im z\ge v_0>0$. Moreover,
\begin{eqnarray*}
&&\bbr_j=\frac{\bby_j-\bar{\bby}_j}{||\bby_j-\bar{\bby}_j||}, \forall j=1,2,\ldots,p; \ \ \underline\bbB_j^{(n)}=\underline\bbB_n-\bbr_j\bbr_j^{*}, \ \ \bbD(z)=\underline\bbB_n-z\bbI_n, \non
&&\bbD_j(z)=\bbD(z)-\bbr_j\bbr_j^{*}, \ \ \beta_j(z)=\frac{1}{1+\bbr_j^{*}\bbD_j^{-1}(z)\bbr_j^{*}}, \ \ \widetilde{\beta}_j(z)=\frac{1}{1+\frac{1}{n}tr\bbD^{-1}_j(z)},\non
&&b_n(z)=\frac{1}{1+\frac{1}{n}\mathbb{E}tr\bbD_1^{-1}(z)}, \ \ \varepsilon_j(z)=\bbr_j^{*}\bbD_j^{-1}(z)\bbr_j-\frac{1}{n}tr\bbD_j^{-1}(z),
\end{eqnarray*}
and $\delta_j(z)=\bbr_j^{*}\bbD_j^{-2}(z)\bbr_j-\frac{1}{n}tr\bbD_j^{-2}(z)$. By Lemma \ref{lem han}, we have for any $r\ge 2$
\begin{eqnarray}\label{1002.2}
\mathbb{E}|\varepsilon_j(z)|^r\le \frac{M}{v^{2r}}n^{-1}\delta_n^{2r-4}
\end{eqnarray}
and
\begin{eqnarray}\label{1002.3}
\mathbb{E}|\delta_j(z)|^r\le \frac{M}{v^{2r}}n^{-1}\delta_n^{2r-4}.
\end{eqnarray}
It is easy to see that
\begin{eqnarray}\label{0626(2)}
\bbD^{-1}(z)-\bbD_{j}^{-1}(z)=-\bbD_{j}^{-1}(z)\bbr_j\bbr_j^{*}\bbD_{j}^{-1}(z)\beta_{j}(z),
\end{eqnarray}
where we use the formula that $\bbA_1^{-1}-\bbA_2^{-1}=\bbA_2^{-1}(\bbA_2-\bbA_1)\bbA_1^{-1}$ holds for any two invertible matrices $\bbA_1$ and $\bbA_2$. Note that $|\beta_j(z)|$, $|\widetilde{\beta}_j(z)|$ and $|b_n(z)|$ are bounded by $\frac{|z|}{v}$.

Let $\mathbb{E}_0(\cdot)$ denote expectation and $\mathbb{E}_j(\cdot)$ denote conditional expectation with respect to the $\sigma$-field generated by $\bbr_1, \bbr_2, \ldots, \bbr_j$, where $j=1,2,\ldots,p$. Next, we write $M_n^{(1)}(z)$ as a sum of martingale difference sequences(MDS), and then utilize the CLT of MDS which is listed in Lemma \ref{lem0521} to derive the asymptotic distribution of $M_n^{(1)}(z)$, which can be written as
\begin{eqnarray}\label{0521(6)}
&& M_n^{(1)}(z) = n[\underline{m}_{n}(z)-\mathbb{E}\underline{m}_{n}(z)]
=tr[\bbD^{-1}(z)-\mathbb{E}\bbD^{-1}(z)]\\
&&=\ \sum^{p}_{j=1}[tr\mathbb{E}_j\bbD^{-1}(z)-tr\mathbb{E}_{j-1}\bbD^{-1}(z)]\non
&&= \ \sum^{p}_{j=1}\Big(tr\mathbb{E}_j[\bbD^{-1}(z)-\bbD^{-1}_j(z)]-tr\mathbb{E}_{j-1}[\bbD^{-1}(z)-\bbD_j^{-1}(z)]\Big)\non
&&= \ -\sum^{p}_{j=1}(\mathbb{E}_j-\mathbb{E}_{j-1})\beta_j(z)\bbr_j^{*}\bbD_j^{-2}(z)\bbr_j.\nonumber
\end{eqnarray}
Evidently, $\beta_j(z)$ can be written as
\begin{eqnarray*}
\beta_j(z)=\widetilde{\beta}_j(z)-\beta_j(z)\widetilde{\beta}_j(z)\varepsilon_j(z)
=\widetilde{\beta}_j(z)-\widetilde{\beta}_j^2(z)\varepsilon_j(z)+\widetilde{\beta}_j^2(z)\beta_j(z)\varepsilon_j^2(z).
\end{eqnarray*}

From this and the definition of $\delta_j(z)$, (\ref{0521(6)}) has the following expression
\begin{eqnarray}\label{0521(7)}
&&\\
&&(\mathbb{E}_j-\mathbb{E}_{j-1})\beta_j(z)\bbr_j^{*}\bbD_j^{-2}(z)\bbr_j = (\mathbb{E}_j-\mathbb{E}_{j-1})\big[\big(\widetilde{\beta}_j(z)-\widetilde{\beta}_j^2(z)\varepsilon_j(z)\non
&&  +\widetilde{\beta}_j^2(z)\beta_j(z)\varepsilon_j^2(z)\big)
\big(\delta_j(z)+\frac{1}{n}tr\bbD_j^{-2}(z)\big)\big] =  (\mathbb{E}_j-\mathbb{E}_{j-1})\big[\widetilde{\beta}_j(z)\delta_j(z)
\nonumber\\
&& - \ \widetilde{\beta}_j^2(z)\varepsilon_j(z)\delta_j(z)-\widetilde{\beta}_j^2(z)\varepsilon_j(z)\frac{1}{n}tr\bbD_j^{-2}(z)+\widetilde{\beta}_j^2(z)\beta_j(z)\varepsilon_j^2(z)\bbr_j^{*}\bbD_j^{-2}(z)\bbr_j\big]\non
&&= \ \mathbb{E}_j\big[\widetilde{\beta}_j(z)\delta_j(z)-\widetilde{\beta}_j^2(z)\varepsilon_j(z)\frac{1}{n}tr\bbD_j^{-2}(z)\big]
\nonumber\\
&& - \ (\mathbb{E}_j-\mathbb{E}_{j-1})\big[\widetilde{\beta}_j^2(z)\big(\varepsilon_j(z)\delta_j(z)-\beta_j(z)\bbr_j^{*}\bbD_j^{-2}(z)\bbr_j\varepsilon_j^2(z)\big)\big] -  \mathbb{E}_{j-1}[\widetilde{\beta}_j(z)\delta_j(z)],\nonumber
\end{eqnarray}
where the second equality uses the fact that $(\mathbb{E}_j-\mathbb{E}_{j-1})\widetilde{\beta}_j(z)\frac{1}{n}tr\bbD_j^{-2}(z)=0$.

By making a minor change to Lemma 8(i.e. Replace $\bbD^{-1}(z)$ by $\bbD_j^{-1}(z)$ ), we have
$\mathbb{E}\Big|\frac{1}{n}\textbf{1}^{T}\bbD_j^{-1}(z)\textbf{1}+\frac{1}{z}\Big|^2\rightarrow 0$.
Thus
\begin{eqnarray}\label{0524(51)}
&&\\
&& -\sum^{p}_{j=1}\mathbb{E}_{j-1}\delta_j(z) = \frac{1}{n(n-1)}\sum^{p}_{j=1}\sum_{k\neq\ell}\mathbb{E}_{j-1}\big(\bbD_j^{-2}(z)\big)_{k\ell}
\nonumber\\
&& = \ \frac{1}{n(n-1)}\sum^{p}_{j=1}\mathbb{E}_{j-1}\big(\textbf{1}^{T}\bbD_j^{-2}(z)\textbf{1}-tr\bbD_j^{-2}(z)\big) \stackrel{i.p.}{\rightarrow} \ \frac{c}{z^2}-c\underline{m}^{'}(z),\nonumber
\end{eqnarray}
where the last step uses the fact that
$
\textbf{1}^{T}\bbD_j^{-2}(z)\textbf{1}=(\textbf{1}^{T}\bbD_j^{-1}(z)\textbf{1})^{'}\rightarrow\frac{1}{z^2} , \ \ \frac{1}{n}tr\bbD_j^{-2}(z)\rightarrow \underline{m}^{'}(z)
$  in $L_2$  by  Lemma 2.3 of \cite{BS2004}.

It follows from (\ref{0524(51)}) that
\begin{eqnarray}\label{0525(80)}
\sum_{j=1}^p\mathbb{E}_{j-1}\widetilde{\beta}_j(z)\delta_j(z)\stackrel{i.p.}{\rightarrow}\frac{cm(z)}{z}-czm(z)\underline{m}^{'}(z),
\end{eqnarray}
where we use the fact that
\begin{eqnarray*}
&& \mathbb{E}\Big|\sum^{p}_{j=1}\mathbb{E}_{j-1}\Big(\hat{\beta}_j(z)-b_1(z)\Big)\delta_j(z)\Big|
\nonumber\\
&& \leq\sum^{p}_{j=1}\frac{1}{n(n-1)}\mathbb{E}\Big|\hat{\beta}_j(z)-b_1(z)\Big|\cdot\Big|\textbf{1}^{T}\bbD_j^{-1}(z)\textbf{1}-tr\bbD_j^{-1}(z)\Big|\rightarrow 0.
\end{eqnarray*}

By (\ref{1002.2}) and (\ref{1002.3}), we have
\begin{eqnarray}\label{0525(81)}
E\big|\sum^{p}_{j=1}(\mathbb{E}_j-\mathbb{E}_{j-1})\widetilde{\beta}_j^2(z)\varepsilon_j(z)\delta_j(z)\big|^2
&=&\sum^{p}_{j=1}\mathbb{E}\big|(\mathbb{E}_j-\mathbb{E}_{j-1})\widetilde{\beta}_j^2(z)\varepsilon_j(z)\delta_j(z)\big|^2\non
&\leq& 4\sum^{p}_{j=1}\mathbb{E}\big|\widetilde{\beta}_j^2(z)\varepsilon_j(z)\delta_j(z)\big|^2=o(1),
\end{eqnarray}
where the first equality uses the fact that $(\mathbb{E}_j-\mathbb{E}_{j-1})\widetilde{\beta}_j^2(z)\varepsilon_j(z)\delta_j(z)$ is a martingale difference sequence.
Therefore, $\sum^{p}_{j=1}(\mathbb{E}_j-\mathbb{E}_{j-1})\widetilde{\beta}_j^2(z)\varepsilon_j(z)\delta_j(z)$ converges to $0$ in probability. By the same argument, we have
\begin{eqnarray}\label{0525(82)}
\sum^{p}_{j=1}(\mathbb{E}_j-\mathbb{E}_{j-1})\widetilde{\beta}_j^2(z)\beta_j(z)\bbr_j^{*}\bbD_j^{-2}(z)\bbr_j\varepsilon_j^2(z)\stackrel{i.p.}{\rightarrow}0.
\end{eqnarray}

With (\ref{0521(6)})--(\ref{0525(82)}), we need to consider the limit of the following term:
\begin{eqnarray}\label{0626(1)}
\sum^{r}_{i=1}\alpha_i\sum^{p}_{j=1}Y_j(z_i)=\sum^{p}_{j=1}\sum^{r}_{i=1}\alpha_iY_j(z_i),
\end{eqnarray}
where $\Im(z_i)\neq 0$, $\{\alpha_i: i=1,2,\ldots,r\}$ are constants, and
\begin{eqnarray*}
Y_j(z)=-\mathbb{E}_j\Big(\widetilde{\beta}_j(z)\delta_j(z)-\widetilde{\beta}_j^2(z)\varepsilon_j(z)\frac{1}{n}tr\bbD_j^{-2}(z)\Big)=-\mathbb{E}_j\frac{d}{dz}\big(\widetilde{\beta}_j(z)\varepsilon_j(z)\big).
\end{eqnarray*}
By Lemma \ref{lem han}, we obtain,
\begin{eqnarray}\label{0521(15)}
\mathbb{E}|Y_j(z)|^4\leq K
\mathbb{E}|\varepsilon_j(z)|^4=o\left(\frac{1}{p}\right).
\end{eqnarray}
It follows from (\ref{0521(15)}) that
\begin{eqnarray*}
\sum^{p}_{j=1}\mathbb{E}\left(\left|\sum^{r}_{i=1}\alpha_iY_j(z_i)\right|^2I_{(|\sum^{r}_{i=1}\alpha_iY_j(z_i)|\geq\varepsilon)}\right)\leq\frac{1}{\varepsilon^2}\sum^{p}_{j=1}\mathbb{E}\left|\sum^{r}_{i=1}\alpha_iY_j(z_i)\right|^4\rightarrow 0.
\end{eqnarray*}
From Lemma \ref{lem0521}, it suffices to prove that
\begin{eqnarray}\label{yang1}
\sum^{p}_{j=1}\mathbb{E}_{j-1}[Y_j(z_1)Y_j(z_2)]
\end{eqnarray}
converges in probability to a constant.
Once it is proved, we can conclude that $M_n^{(1)}(z)$ converges in finite dimension to a normal distribution. 

Since
\begin{eqnarray*}
\frac{\partial^2}{\partial z_1\partial z_2}\Big(\sum^{p}_{j=1}\mathbb{E}_{j-1}\big[\mathbb{E}_j\big(\widetilde{\beta}_j(z_1)\varepsilon_j(z_1)\big)\mathbb{E}_j\big(\widetilde{\beta}_j(z_2)\varepsilon_j(z_2)\big)\big]\Big)=(\ref{yang1}),
\end{eqnarray*}
and by the same arguments as those on page 571 of \cite{BS2004},
it is enough to consider the limit of
\begin{eqnarray}\label{yang2}
\sum^{p}_{j=1}\mathbb{E}_{j-1}\big[\mathbb{E}_j\big(\widetilde{\beta}_j(z_1)\varepsilon_j(z_1)\big)\mathbb{E}_j\big(\widetilde{\beta}_j(z_2)\varepsilon_j(z_2)\big)\big].
\end{eqnarray}

By the fact that
\begin{eqnarray}\label{0524(6)}
\lim_{n\rightarrow}\mathbb{E}\Big|\frac{1}{n}tr\mathbb{E}_j\Big(\bbD_j^{-1}(z)\Big)-\underline{m}(z)\Big|^2=0,
\end{eqnarray}
$|\beta_j(z)\bbr_j^*\bbD_j^{-2}(z)\bbr_j|\le \frac{|z|}{v^2}$  and Lemma \ref{lem0527}, we obtain
\begin{eqnarray}\label{0521(20)}
\mathbb{E}\big|\widetilde{\beta}_j(z_i)-b_n(z_i)\big|^2\leq \frac{K}{n}.
\end{eqnarray}
By (\ref{0521(20)}), we have
\begin{eqnarray*}
\mathbb{E}\big|
\mathbb{E}_{j-1}[\mathbb{E}_j\big(\widetilde{\beta}_j(z_1)\varepsilon_j(z_1)\big)\mathbb{E}_j\big(\widetilde{\beta}_j(z_2)\varepsilon_j(z_2)\big)]
-\mathbb{E}_{j-1}[\mathbb{E}_j\big(b_n(z_1)\varepsilon_j(z_1)\big)\mathbb{E}_j\big(b_n(z_2)\varepsilon_j(z_2)\big)]\big|
=O(n^{-3/2}).
\end{eqnarray*}
From this, it follows that
\begin{eqnarray*}
\sum^{p}_{j=1}\mathbb{E}_{j-1}\big[\mathbb{E}_j\big(\widetilde{\beta}_j(z_1)\varepsilon_j(z_1)\big)\mathbb{E}_j\big(\widetilde{\beta}_j(z_2)\varepsilon_j(z_2)\big)\big]
-b_n(z_1)b_n(z_2)\sum^{p}_{j=1}\mathbb{E}_{j-1}\big[\mathbb{E}_j\big(\varepsilon_j(z_1)\big)\mathbb{E}_j\big(\varepsilon_j(z_2)\big)\big]
\stackrel{i.p.}{\rightarrow}0.
\end{eqnarray*}
Then it is enough to prove that
\begin{eqnarray}\label{yang3}
b_n(z_1)b_n(z_2)\sum^{p}_{j=1}\mathbb{E}_{j-1}\big[\mathbb{E}_j\big(\varepsilon_j(z_1)\big)\mathbb{E}_j\big(\varepsilon_j(z_2)\big)\big]
\end{eqnarray}
converges to a constant in probability, which further gives the limit of (\ref{0626(1)}).

By Lemma \ref{lem yang}, (\ref{yang3}) becomes
\begin{eqnarray}\label{yang4}
&&\non
&&(\ref{yang3})=\left\{\begin{array}{cc}
  J_1+2J_2+J_3+o_P(1), & under\  the \ real \  case; \\
  J_1+J_2+J_3+J_4+o_P(1), & under \ the \ complex \ case,
\end{array}
\right.
\end{eqnarray}
where
\begin{eqnarray*}
&&J_1=\frac{1}{n^3}b_n(z_1)b_n(z_2)\Big[\sum^{p}_{j=1}(1-\mathbb{E}|X_{j1}|^4)tr\mathbb{E}_j\big(\bbD_j^{-1}(z_1)\big)tr\mathbb{E}_j\big(\bbD_j^{-1}(z_2)\big)\Big];\non
&&J_2=\frac{1}{n^2}b_n(z_1)b_n(z_2)\Big[\mathbb{E}\sum^{p}_{j=1}tr\big[\mathbb{E}_j\big(\bbD_j^{-1}(z_1)\big)\mathbb{E}_j\big(\bbD_j^{-1}(z_2)\big)\big];\non
&&J_3=\frac{1}{n^2}b_n(z_1)b_n(z_2)\Big[\sum^{p}_{j=1}(\mathbb{E}|X_{j1}|^4-2-|\mathbb{E}X_{j1}^2|^2)\sum^{n}_{k=1}\mathbb{E}_j\big(\bbD_j^{-1}(z_1)\big)_{kk}\mathbb{E}_j\big(\bbD_j^{-1}(z_2)\big)_{kk}\Big];
\nonumber\\
&&J_4=b_n(z_1)b_n(z_2)\frac{1}{n^2}\Big[\sum^{p}_{j=1}|\mathbb{E}X_{j1}^2|^2\big[\mathbb{E}_j\big(\bbD_j^{-1}(z_1)\big)\mathbb{E}_j\big(\bbD_j^{-1}(z_2)^{T}\big)\big]\Big].
\end{eqnarray*}

Next, we study the limit of the term $J_2$. Let $\bbD_{ij}(z)=\bbD(z)-\bbr_i\bbr_i^{*}-\bbr_j\bbr_j^{*}$, $b_1(z)=\frac{1}{1+\frac{1}{n}\mathbb{E}tr\bbD_{12}^{-1}(z)}$ and $\beta_{ij}(z)=\frac{1}{1+\bbr_i^{*}\bbD_{ij}^{-1}(z)\bbr_i}$.

We have the equality $\bbD_j(z_1)+z_1\bbI_n-\frac{p-1}{n}b_1(z_1)\bbI_n=\sum^{p}_{i\neq j}\bbr_i\bbr_i^{*}-\frac{p-1}{n}b_1(z_1)\bbI_n$. Multiplying by $(z_1\bbI_n-\frac{p-1}{n}b_1(z)\bbI_n)^{-1}$ on the left-hand side and $\bbD_j^{-1}(z_1)$ on the right-hand side, and using
\begin{eqnarray}\label{0521(40)}
\bbr_i^{*}\bbD_j^{-1}(z_1)=\beta_{ij}(z_1)\bbr_i^{*}\bbD_{ij}^{-1}(z_1),
\end{eqnarray}
we get
\begin{eqnarray}\label{yang10}
\bbD_j^{-1}(z_1)&=&-\bbH_n(z_1)+\sum_{i\neq j}^{p}\beta_{ij}(z_1)\bbH_n(z_1)\bbr_i\bbr_i^{*}\bbD_{ij}^{-1}(z_1)
-\frac{p-1}{n}b_1(z_1)\bbH_n(z_1)\bbD_j^{-1}(z_1)\non
&=&-\bbH_n(z_1)+b_1(z_1)A(z_1)+B(z_1)+C(z_1),
\end{eqnarray}
where $\bbH_n(z_1)=\big(z_1\bbI_n-\frac{p-1}{n}b_1(z_1)\bbI_n\big)^{-1}$, $A(z_1)=\sum_{i\neq j}^{p}\bbH_n(z_1)\big(\bbr_i\bbr_i^{*}-\frac{1}{n}\bbI_n\big)\bbD_{ij}^{-1}(z_1)$, $B(z_1)=\sum_{i\neq j}^{p}\big(\beta_{ij}(z_1)-b_1(z_1)\big)\bbH_n(z_1)\bbr_i\bbr_i^{*}\bbD_{ij}^{-1}(z_1)$ and $$C(z_1)=n^{-1}b_1(z_1)\bbH_n(z_1)\sum_{i\neq j}^{p}\big(\bbD_{ij}^{-1}(z_1)-\bbD_j^{-1}(z_1)\big).$$

For any real $t$, $\Big|1-\frac{t}{z(1+n^{-1}\mathbb{E}tr\bbD_{12}^{-1}(z))}\Big|^{-1} \leq\frac{\big|z\big(1+n^{-1}\mathbb{E}tr\bbD_{12}^{-1}(z)\big)\big|}{\Im\big(z(1+n^{-1}\mathbb{E}tr\bbD_{12}^{-1}(z))\big)}\leq \frac{|z|(1+1/v_0)}{v_0}$.

Thus,
\begin{eqnarray}\label{yang6}
||\bbH_n(z_1)||\leq\frac{1+1/v_0}{v_0}.
\end{eqnarray}

For any random matrix $\bbM$, denote its nonrandom bound by $|||\bbM|||$.
From (\ref{0521(20)}), Lemma \ref{lem han} and (\ref{yang6}), we get
\begin{eqnarray}\label{yang8}
\mathbb{E}\big|tr\bbB(z_1)\bbM\big|&\leq&p\mathbb{E}^{1/2}\big(\big|\beta_{12}(z_1)-b_1(z_1)\big|^2\big)\cdot\mathbb{E}^{1/2}\big(\big|\bbr_i^{*}\bbD_{ij}^{-1}(z_1)\bbM H_n(z_1)\bbr_i\big|^2\big)\non
&\leq&K|||\bbM|||\frac{|z_1|^2(1+1/v_0)}{v_0^5}n^{1/2}.
\end{eqnarray}

For any $n\times n$ matrix $\bbA$, we have
\begin{eqnarray}\label{yang7}
\big|tr\big(\bbD^{-1}(z)-\bbD_j^{-1}(z)\big)\bbA\big|\leq\frac{||\bbA||}{\Im(z)},
\end{eqnarray}

With (\ref{yang7}), we obtain
\begin{eqnarray}\label{1002.5}
\big|tr\bbC(z_1)\bbM\big|\leq |||\bbM|||\frac{|z_1|(1+1/v_0)}{v_0^3}.
\end{eqnarray}

For any nonrandom $\bbM$, it follows from Lemma \ref{lem han} and (\ref{yang6}) that
\begin{eqnarray}\label{1002.4}
\mathbb{E}\big|tr\bbA(z_1)\bbM\big|&\leq&K\mathbb{E}^{1/2}\big(tr\bbD_{ij}^{-1}(z_1)\bbM\bbH_n(z_1)\\
&&\ \ \ \ \ \ \ \ \ \ \ \ \cdot\bbH_n(\bar z_1)\bbM^{*}\bbD_{ij}^{-1}(\bar z_1)\big)\non
&\leq&K||\bbM||\frac{1+1/v_0}{v_0^2}n^{1/2}.\nonumber
\end{eqnarray}

By using (\ref{0626(2)}), we can write $tr\mathbb{E}_j\big(\bbA(z_1)\big)\bbD_j^{-1}(z_2)=A_1(z_1,z_2)+A_2(z_1,z_2)+A_3(z_1,z_2)+R(z_1,z_2)$, where
\medskip

\begin{eqnarray*}
A_1(z_1,z_2)&=&-tr\sum_{i<j}^{p}\bbH_n(z_1)\bbr_i\bbr_i^{*}\mathbb{E}_j\big(\bbD_{ij}^{-1}(z_1)\big)\beta_{ij}(z_2)\bbD_{ij}^{-1}(z_2)\bbr_i\bbr_i^{*}\bbD_{ij}^{-1}(z_2)\non
&=&-\sum_{i<j}^{p}\beta_{ij}(z_2)\bbr_i^{*}\mathbb{E}_j\big(\bbD_{ij}^{-1}(z_1)\big)\bbD_{ij}^{-1}(z_2)\bbr_i\bbr_i^{*}\bbD_{ij}^{-1}(z_2)\bbH_n(z_1)\bbr_i;\non
A_2(z_1,z_2)&=&-tr\sum_{i<j}^{p}\bbH_n(z_1)n^{-1}\mathbb{E}_j\big(\bbD_{ij}^{-1}(z_1)\big)\big(\bbD_j^{-1}(z_2)-\bbD_{ij}^{-1}(z_2)\big);\non
A_3(z_1,z_2)&=&tr\sum_{i<j}^{p}\bbH_n(z_1)\big(\bbr_i\bbr_i^{*}-n^{-1}\bbI_n\big)\mathbb{E}_j\big(\bbD_{ij}^{-1}(z_1)\big)\bbD_{ij}^{-1}(z_2);\non
R(z_1 , z_2)&=&tr\mathbb{E}_{j}\sum_{i>j}\bbH_n(z_1)\big(-\frac{1}{n(n-1)}\bbe\bbe^{*}+\frac{1}{n(n-1)}\bbI_n\big)\bbD_{ij}^{-1}(z_1)\bbD_{j}^{-1}(z_2),
\end{eqnarray*}
where $\bbe$ is an $n$-dimensional vector with all elements being $1$ and $\mathbb{E}\bar{\bbr}_{ik}\bbr_{ij}=-\frac{1}{n(n-1)}$, $k\neq j$(see (\ref{y0})).

It is easy to see that $R(z_1 , z_2)=O(1)$.
We get from (\ref{yang7}) and (\ref{yang6}) that $|A_2(z_1,z_2)|\leq\frac{1+1/v_0}{v_0^2}$. Similar to (\ref{yang8}), we have $\mathbb{E}|A_3(z_1,z_2)|\leq\frac{1+1/v_0}{v_0^3}n^{1/2}$. Using Lemma \ref{lem han} and (\ref{0521(20)}), we have, for $i<j$,
\begin{eqnarray}\label{0521(22)}
&&\mathbb{E}\big|\beta_{ij}(z_2)\bbr_i^{*}\mathbb{E}_j\big(\bbD_{ij}^{-1}(z_1)\big)\bbD_{ij}^{-1}(z_2)\bbr_i\bbr_i^{*}\bbD_{ij}^{-1}(z_2)\bbH_n(z_1)\bbr_i\non
&&\ \ \ \ \ -b_1(z_2)n^{-2}tr\big(\mathbb{E}_j\big(\bbD_{ij}^{-1}(z_1)\big)\bbD_{ij}^{-1}(z_2)\big)tr\big(\bbD_{ij}^{-1}(z_2)\bbH_n(z_1)\big)\big|\leq Kn^{-1/2}.
\end{eqnarray}
By (\ref{yang7}), we have
\begin{eqnarray}\label{0521(23)}
&&\big|tr\big(\mathbb{E}_j\big(\bbD_{ij}^{-1}(z_1)\big)\bbD_{ij}^{-1}(z_2)\big)tr\big(\bbD_{ij}^{-1}(z_2)\bbH_n(z_1)\big)\non
&&\ \ -tr\big(\mathbb{E}_j\big(\bbD_{j}^{-1}(z_1)\big)\bbD_{j}^{-1}(z_2)\big)tr\big(\bbD_{j}^{-1}(z_2)\bbH_n(z_1)\big)\big|\leq Kn.
\end{eqnarray}
It follows from (\ref{0521(22)}) and (\ref{0521(23)}) that
\begin{eqnarray*}
\mathbb{E}\big|A_1(z_1,z_2)+\frac{j-1}{n^2}b_1(z_2)tr\big(\mathbb{E}_j\big(\bbD_j^{-1}(z_1)\big)\bbD_j^{-1}(z_2)\big)tr\big(\bbD_j^{-1}(z_2)\bbH_n(z_1)\big)\big|\leq Kn^{1/2}.
\end{eqnarray*}

Therefore, by (\ref{yang10})--(\ref{0521(23)}), we obtain that
\begin{eqnarray*}
&&tr\big(\mathbb{E}_j\big(\bbD_j^{-1}(z_1)\big)\bbD_j^{-1}(z_2)\big)
\big[1+\frac{j-1}{n^2}b_1(z_1)b_1(z_2)tr\big(\bbD_j^{-1}(z_2)\bbH_n(z_1)\big]\non
&=&-tr\big(\bbH_n(z_1)\bbD_j^{-1}(z_2)\big)+A_4(z_1,z_2),
\end{eqnarray*}
where $\mathbb{E}|A_4(z_1,z_2)|\leq Kn^{1/2}$.

By (\ref{yang10}) for $\bbD_j^{-1}(z_2)$ and (\ref{yang8})-(\ref{1002.4}), we have
\begin{eqnarray}\label{0521(32)}
&&tr\big(\mathbb{E}_j\big(\bbD_j^{-1}(z_1)\big)\bbD_j^{-1}(z_2)\big)
\big[1-\frac{j-1}{n^2}b_1(z_1)b_1(z_2)tr\big(\bbH_n(z_2)\bbH_n(z_1)\big)\big]\\
&=&tr\Big(\bbH_n(z_2)\bbH_n(z_1)\Big)+A_5(z_1,z_2),\nonumber
\end{eqnarray}
where $\mathbb{E}|A_5(z_1,z_2)|\leq Kn^{1/2}$.

From (\ref{yang7}), we have $|b_1(z)-b_n(z)|\leq Kn^{-1}$. Using $\mathbb{E}\big|\frac{1}{n}tr\bbD^{-1}(z)-\mathbb{E}\frac{1}{n}tr\bbD^{-1}(z)\big|^{k}\leq C_kn^{-k/2}$,
we have
\begin{eqnarray}\label{0521(31)}
|b_n(z)-\mathbb{E}\beta_1(z)|\leq Kn^{-1/2}.
\end{eqnarray}

As in $(2.2)$ of \cite{S1995}, one may verify that
\begin{eqnarray}\label{0521(28)}
m_n(z)=-\frac{1}{pz}\sum^{p}_{j=1}\beta_j(z).
\end{eqnarray}
It follows from (\ref{0521(28)}) that
\begin{eqnarray}\label{0521(29)}
\mathbb{E}\beta_1(z)=-z\mathbb{E}m_n(z).
\end{eqnarray}
From (\ref{0521(31)}), (\ref{0521(29)}) and Lemma \ref{lem han}, we have
\begin{eqnarray}\label{0524(1)}
|b_1(z)+zm_{c_n}(z)|\leq Kn^{-1/2}.
\end{eqnarray}
Let $\bbQ_n(z)=\big(\bbI_n+\frac{p-1}{n}m_{c_n}(z)\bbI_n\big)^{-1}.$ So by (\ref{0521(32)}), we get
\begin{eqnarray}\label{yang13}
&&tr\big(\mathbb{E}_j\big(\bbD_j^{-1}(z_1)\big)\bbD_j^{-1}(z_2)\big)
\big[1-\frac{j-1}{n^2}m_{c_n}(z_1)m_{c_n}(z_2)tr\bbQ_n(z_2)\bbQ_n(z_1)\big]\\
&&=\frac{1}{z_1z_2}tr\big(\bbQ_n(z_2)\bbQ_n(z_1)\big)+A_6(z_1,z_2),\nonumber
\end{eqnarray}
where $\mathbb{E}|A_6(z_1,z_2)|\leq Kn^{1/2}$.

Rewrite (\ref{yang13}) as
\begin{eqnarray}\label{0521(35)}
&&tr\big(\mathbb{E}_j\big(\bbD_j^{-1}(z_1)\big)\bbD_j^{-1}(z_2)\big)
\big[1-\frac{j-1}{n}\frac{m_{c_n}(z_1)m_{c_n}(z_2)}{\big(1+\frac{p-1}{n}m_{c_n}(z_2)\big)\big(1+\frac{p-1}{n}m_{c_n}(z_1)\big)}\big]\\
&&=\frac{n}{z_1z_2}\frac{1}{\big(1+\frac{p-1}{n}m_{c_n}(z_1)\big)\big(1+\frac{p-1}{n}m_{c_n}(z_2)\big)}+A_6(z_1,z_2).\nonumber
\end{eqnarray}
Then $J_2$ can be written as $J_2=a_n(z_1,z_2)\frac{1}{p}\sum^{p}_{j=1}\frac{1}{1-\frac{j-1}{p}a_n(z_1,z_2)}+A_7(z_1,z_2)$, where
$a_n(z_1,z_2)=\frac{c_n m_{c_n}(z_1)m_{c_n}(z_2)}{\big(1+\frac{p-1}{n}m_{c_n}(z_1)\big)\big(1+\frac{p-1}{n}m_{c_n}(z_2)\big)}$ and $\mathbb{E}|A_7(z_1,z_2)|\leq Kn^{-1/2}$.

Note that the limit of $a_n(z_1, z_2)$ is $a(z_1,z_2)=\frac{c m(z_1)m(z_2)}{\big(1+cm(z_1)\big)\big(1+cm(z_2)\big)}$. Thus by (\ref{0521(35)}), the i.p. limit of $\frac{\partial^2}{\partial z_2\partial z_1}J_2$ is
\begin{eqnarray*}
&&\frac{\partial^2}{\partial z_2\partial z_1}\int^{a(z_1,z_2)}_{0}\frac{1}{1-z}dz
=\frac{\partial}{\partial z_2}\big(\frac{\partial a(z_1,z_2)/\partial z_1}{1-a(z_1,z_2)}\big)\non
&=&\frac{\partial}{\partial z_2}\Big(\frac{cm(z_2)m^{'}(z_1)}{\Big(1+cm(z_1)\Big)\Big(1+c\big(m(z_1)+m(z_2)\big)+c(c-1)m(z_1)m(z_2)\Big)}\Big)\non
&=&\frac{cm^{'}(z_1)m^{'}(z_2)}{\Big(1+c\big(m(z_1)+m(z_2)\big)+c(c-1)m(z_1)m(z_2)\Big)^2}.
\end{eqnarray*}

For $\frac{\partial^2}{\partial z_1\partial z_2}J_1$ in (\ref{yang4}), similar to (\ref{0524(51)}), by (\ref{0524(1)}) and (\ref{0524(6)}),
we have
\begin{eqnarray*}
&&\mathbb{E}\left|\frac{1}{n}\sum_{j=1}^p\mathbb{E}|X_{j1}|^4\frac{1}{n}tr\mathbb{E}_j\bbD_j^{-1}(z_1)\frac{1}{n}tr\mathbb{E}_j\bbD_j^{-1}(z_2)-\frac{\underline m_{c_n}(z_1)\underline m_{c_n}(z_2)}{n}\sum_{j=1}^p\mathbb{E}|X_{j1}|^4\right|=o(1)
\end{eqnarray*}
So we can conclude that
\begin{eqnarray*}
J_1\stackrel{i.p.}{\rightarrow}\Big(1-\kappa\Big)cz_1m(z_1)z_2m(z_2)\underline{m}(z_1)\underline{m}(z_2)
=\Big(1-\kappa\Big)\frac{c\underline{m}(z_1)\underline{m}(z_2)}{\big(1+\underline{m}(z_1)\big)\big(1+\underline{m}(z_2)\big)},
\end{eqnarray*}
where the equality above uses the relation between $m(z)$ and $\underline{m}(z)$: $m(z)=-\frac{1}{z\big(1+\underline{m}(z)\big)}$. Then the second derivative of $J_1$ with respect to $z_1$ and $z_2$ is
\begin{eqnarray*}
\frac{\partial^2}{\partial z_1\partial z_2}J_1\stackrel{i.p.}{\rightarrow}\Big(1-\kappa\Big)\frac{c\underline{m}^{'}(z)\underline{m}^{'}(z_2)}{\big(1+\underline{m}(z_1)\big)^2\big(1+\underline{m}(z_2)\big)^2}.
\end{eqnarray*}

The next aim is to establish the limit of $\frac{\partial^2}{\partial z_1\partial z_2}J_3$ in (\ref{yang4}).
It is enough to find the limit of $\frac{1}{n^2}\sum^{n}_{k=1}\sum^{p}_{j=1}(\mathbb{E}|X_{j1}|^4-2-|\mathbb{E}X_{j1}^2|^2)\mathbb{E}_j\big(\bbD_j^{-1}(z_1)\big)_{kk}\mathbb{E}_j\big(\bbD_j^{-1}(z_2)\big)_{kk}$.

First, we claim that
\begin{eqnarray}\label{yang40}
&&\\
\frac{1}{n^2}\sum^{n}_{k=1}\sum^{p}_{j=1}(\mathbb{E}|X_{j1}|^4-2-|\mathbb{E}X_{j1}^2|^2)\mathbb{E}_j\big(\bbD_j^{-1}(z_1)-\mathbb{E}\bbD_j^{-1}(z_1)\big)_{kk}\mathbb{E}_j\big(\bbD_j^{-1}(z_2)\big)_{kk}=O_p(n^{-1/2}).\nonumber
\end{eqnarray}

Actually,
\begin{eqnarray}\label{yang40}\\
&&\mathbb{E}|\frac{1}{n^2}\sum^{n}_{k=1}\sum^{p}_{j=1}(\mathbb{E}|X_{j1}|^4-2-|\mathbb{E}X_{j1}^2|^2)\mathbb{E}_j\big(\bbD_j^{-1}(z_1)-\mathbb{E}\bbD_j^{-1}(z_1)\big)_{kk}\mathbb{E}_j\big(\bbD_j^{-1}(z_2)\big)_{kk}|
\nonumber\\
&\le&\frac{pK}{n^2 v_0} \sum^{n}_{k=1}\mathbb{E}|\bbe'_{k}(\bbD_{1}^{-1}(z_1)-\mathbb{E}\bbD_{1}^{-1}(z_1))\bbe_{k}| \leq K n^{-1/2},\nonumber
\end{eqnarray}
where the last inequality follows from (\ref{1009.1}) (replacing $\bbD$ by $\bbD_1 $). 
With (\ref{yang40}), it remains to find the limit of
\begin{eqnarray}\label{0803.1}
\frac{1}{n}\sum^{n}_{k=1}\mathbb{E}\big(\bbD_j^{-1}(z_1)\big)_{kk}\mathbb{E}\big(\bbD_j^{-1}(z_2)\big)_{kk}.
\end{eqnarray}
It is easy to see that the sum of expectations in (\ref{0803.1}) is exactly the same for any j. Moreover, we have
\begin{eqnarray}\label{0524(5)}
\frac{1}{n}\sum^{n}_{k=1}\mathbb{E}\big(\bbD_j^{-1}(z_1)\big)_{kk}\mathbb{E}\big(\bbD_j^{-1}(z_2)\big)_{kk}\stackrel{i.p.}{\rightarrow}\underline{m}(z_1)\underline{m}(z_2).
\end{eqnarray}

By (\ref{*}), (\ref{0524(1)}) and (\ref{0524(5)}), we get $J_3\stackrel{i.p.}{\rightarrow}\Big(\kappa-2-|\mathbb{E}X_{11}^2|^2\Big)cz_1z_2m(z_1)\underline{m}(z_1)m(z_2)\underline{m}(z_2)$.

Thus the limit of $\frac{\partial^2}{\partial z_1\partial z_2}J_3$ is
\begin{eqnarray}\label{0521(60)}\\
&&\frac{\partial^2}{\partial z_1\partial z_2}J_3\stackrel{i.p.}{\rightarrow}\Big(\kappa-|\mathbb{E}X_{11}^2|^2-2\Big)
 \cdot c(m(z_1)\underline{m}(z_1)+z_1m(z_1)\underline{m}^{'}(z_1)+z_1m'(z_1)\underline{m}(z_1))
\nonumber\\
 && \times \ (m(z_2)\underline{m}(z_2)+z_2m(z_2)\underline{m}^{'}(z_2)+z_2m'(z_2)\underline{m}(z_2)).
 \nonumber
\end{eqnarray}

For the complex case, the limit of $\frac{\partial^2}{\partial z_1\partial z_2}J_4$ is derived in Lemma 9.  \\

\textbf{Step 2}:

The tightness of $M_n^{(1)}(z)$ is similar to that provided in \cite{BS2004}. It is sufficient to prove the moment condition (12.51) of \cite{Bill1968}, i.e.
\begin{eqnarray}\label{0718(0)}
\sup_{n; z_1, z_2\in\mathcal{C}_n}\frac{\mathbb{E}|M_n^{(1)}(z_1)-M_n^{(1)}(z_2)|^2}{|z_1-z_2|^2}
\end{eqnarray}
is finite.

Before proceeding, we provide some results needed in the proof later. First, moments of $||\bbD^{-1}(z)||$, $||\bbD^{-1}_j(z)||$ and $||\bbD_{ij}^{-1}(z)||$ are bounded in $n$ and $z\in\mathcal{C}_n$. It is easy to see that it is true for $z\in\mathcal{C}_u$ and for $z\in\mathcal{C}_{\ell}$ if $x_{\ell}<0$. For $z\in\mathcal{C}_r$ or, if $x_{\ell}>0$, $z\in\mathcal{C}_{\ell}$, we have from Lemma \ref{lem eig} that
\begin{eqnarray*}
\mathbb{E}||\bbD^{-1}_j(z)||^{m}&\leq&K_1+v^{-m}P(||\bbB_{j}||\geq\eta_r \ or \ \lambda_{\min}^{\bbB_{j}}\leq\eta_{\ell})\non
&\leq&K_1+K_2n^{m}\varepsilon^{-m}n^{-\ell}\leq K
\end{eqnarray*}
for large $\ell$. Here $\eta_r$ is any number between $(1+\sqrt{c})^2$ and $x_r$; if $x_{\ell}>0$, $\eta_{\ell}$ is any number between $x_{\ell}$ and $(1-\sqrt{c})^2$ and if $x_{\ell}<0$, $\eta_{\ell}$ can be any negative number. So for any positive integer $m$,
\begin{eqnarray}
\max\Big(\mathbb{E}||\bbD^{-1}(z)||^{m}, \mathbb{E}||\bbD_j^{-1}(z)||^{m}, \mathbb{E}||\bbD_{ij}^{-1}(z)||^{m}\Big)\leq K.
\end{eqnarray}
By the argument above, we can extend (\ref{1019.2}) in the remark of Lemma \ref{lem han} and get
\begin{eqnarray}\label{0718(1)}
\Big|\mathbb{E}\Big(a(v)\prod^{q}_{\ell=1}\big(\bbr_1^{*}\bbB_{\ell}(v)\bbr_1-n^{-1}tr\bbB_{\ell}(v)\big)\Big)\Big|\leq Kn^{-1}\delta_n^{(2q-4)},
\end{eqnarray}
where $\bbB_{\ell}(v)$ is independent of $\bbr_1$ and
\begin{eqnarray*}
\max(|a(v)|, ||\bbB_{\ell}(v)||)\leq K(1+n^{s}I(||\bbB_n||\geq\eta_r \ or \ \lambda_{\min}^{\widetilde{\bbB}}\leq\eta_{\ell})),
\end{eqnarray*}
with $\widetilde{\bbB}$ being $\underline{\bbB}_j^{(n)}$ or $\underline{\bbB}_n$.

By (\ref{0718(1)}), we have
\begin{eqnarray}\label{0718(3)}
\mathbb{E}|\varepsilon_j(z)|^{m}\leq K_mn^{-1}\delta_n^{2m-4}.
\end{eqnarray}

Let $\gamma_j(z)=\bbr_j^{*}\bbD_j^{-1}(z)\bbr_j-n^{-1}\mathbb{E}tr\bbD^{-1}_j(z)$.
By Lemma \ref{lem0527}, (\ref{0718(1)}) and H$\ddot{o}$lder's inequality, with similar derivation on page 580 of \cite{BS2004}, we have
\begin{eqnarray}\label{0718(2)}
\mathbb{E}|\gamma_j(z)-\varepsilon_j(z)|^{m}\leq\frac{K_m}{n^{m/2}}.
\end{eqnarray}
It follows from (\ref{0718(3)}) and (\ref{0718(2)}) that
\begin{eqnarray}\label{0718(4)}
\mathbb{E}|\gamma_j|^{m}\leq K_mn^{-1}\delta_n^{2m-4}, \ \ m\geq 2.
\end{eqnarray}

Next, we prove that $b_n(z)$ is bounded. With (\ref{0718(1)}), we have for any $m\geq 1$,
\begin{eqnarray}\label{0718(5)}
\mathbb{E}|\beta_1(z)|^{m}\leq K_m.
\end{eqnarray}
Since $b_n(z)=\beta_1(z)+\beta_1(z)b_n(z)\gamma_1(z)$, it is derived from (\ref{0718(4)}) and (\ref{0718(5)}) that $|b_n(z)|\leq K_1+K_2|b_n(z)|n^{-1/2}$.

Then
\begin{eqnarray}\label{0718(6)}
|b_n(z)|\leq\frac{K_1}{1-K_2n^{-1/2}}\leq K.
\end{eqnarray}

With (\ref{0718(1)})--(\ref{0718(6)}) and the same approach on page 581-583 of \cite{BS2004}, we can obtain that (\ref{0718(0)}) is finite.

\textbf{Steps 3 and 4}:

First, we list some results which are used later in this part. The derivations of these results are similar to those for sample covariance matrices in \cite{BS2004}:
\begin{eqnarray*}
\sup_{z\in\mathcal{C}_n}|\underline{m}_n(z)-\underline{m}(z)|\rightarrow 0, \ \ as \ \ n\rightarrow\infty.
\end{eqnarray*}
\begin{eqnarray}\label{yang20}
\sup_{n;z\in\mathcal{C}_n}\big|\big|\big(c_n\mathbb{E}\underline{m}_n(z)\bbI_n+\bbI_n\big)^{-1}\big|\big|<\infty.
\end{eqnarray}
\begin{eqnarray*}
\sup_{z\in\mathcal{C}_n}\big|\frac{\mathbb{E}\big(\underline{m}_n^2(z)\big)}{\big(1+c_n\mathbb{E}\underline{m}_n(z)\big)^2}\big|<\xi, \ \ \xi\in(0,1).
\end{eqnarray*}
\begin{eqnarray}\label{yang24}
\mathbb{E}\big|tr\bbD^{-1}_1(z)\bbM-\mathbb{E}tr\bbD^{-1}_1(z)\bbM\big|^2\leq K||\bbM||^2.
\end{eqnarray}

Next, we derive an identity.
Let
\begin{eqnarray}\label{gn}
\bbG_n(z)=c_n\mathbb{E}m_n(z)\bbI_n+\bbI_n.
\end{eqnarray}
Write $\underline{\bbB}_n-z\bbI_n-\big(-z\bbG_n(z)\big)$ as $\sum^{p}_{j=1}\bbr_j\bbr_j^{*}-\big(-zc_n\mathbb{E}m_n(z)\big)\bbI_n$. Taking first inverse and then expected value, we get
\begin{eqnarray*}
&&\big(-z\bbG_n(z)\big)^{-1}-\mathbb{E}\big(\underline{\bbB}_n-z\bbI_n\big)^{-1}\non
&=&\big(-z\bbG_n(z)\big)^{-1}\mathbb{E}\big[\big(\sum^{p}_{j=1}\bbr_j\bbr_j^{*}-(-zc_{n}\mathbb{E}m_n(z)\bbI_n)\big)\big(\underline{\bbB}_n-z\bbI_n\big)^{-1}\big]\non
&=&-z^{-1}\sum^{p}_{j=1}\mathbb{E}\beta_j(z)\big[\bbG^{-1}_n(z)\bbr_j\bbr_j^{*}\big(\underline{\bbB}_{(j)}^{n}-z\bbI_n\big)^{-1}\big]\non
&&+z^{-1}\mathbb{E}\big[\bbG^{-1}_n(z)\big(-zc_{n}\mathbb{E}m_n(z)\big)\bbI_n\big(\underline{\bbB}_n-z\bbI_n\big)^{-1}\big]\non
&=&-z^{-1}\sum^{p}_{j=1}\mathbb{E}\beta_j(z)\big[\bbG^{-1}_n(z)\bbr_j\bbr_j^{*}\big(\underline{\bbB}_{(j)}^{n}-z\bbI_n\big)^{-1}-\frac{1}{n}\bbG^{-1}_n(z)\mathbb{E}\big(\underline{\bbB}_n-z\bbI_n\big)^{-1}\big]\non
&=&-z^{-1}\sum^{p}_{j=1}\mathbb{E}\beta_j(z)\big[\bbG^{-1}_n(z)\bbr_j\bbr_j^{*}\bbD_j^{-1}(z)-\frac{1}{n}\bbG^{-1}_n(z)\mathbb{E}\bbD^{-1}(z)\big],
\end{eqnarray*}
where the second equality uses (\ref{0521(40)}) and the third equality uses (\ref{0521(29)}).

Taking trace on both sides and dividing by $-\frac{p}{z}$, we get
\begin{eqnarray*}
&&\frac{1}{c_n\big(1+c_{n}\mathbb{E}m_n(z)\big)}+\frac{z}{c_n}\mathbb{E}\underline{m}_n(z)\non
&=&\frac{1}{p}\sum^{p}_{j=1}\mathbb{E}\beta_j(z)\big[\bbr_j^{*}\bbD_j^{-1}(z)\bbG^{-1}_n(z)\bbr_j-\frac{1}{n}\mathbb{E}tr\left(  \bbG^{-1}_n(z)\bbD^{-1}(z)\right)\big],
\end{eqnarray*}
Next, we investigate the limit of
\begin{eqnarray}\label{0521(46)}
&&n\Big(\frac{1}{c_n\big(1+c_{n}\mathbb{E}m_n(z)\big)}+\frac{z}{c_n}\mathbb{E}\underline{m}_n(z)\Big)\non
&=&\frac{n}{p}\sum^{p}_{j=1}\mathbb{E}\beta_j(z)\big[\bbr_j^{*}\bbD_j^{-1}(z)\bbG^{-1}_n(z)\bbr_j-\frac{1}{n}\mathbb{E}tr\left(  \bbG^{-1}_n(z)\bbD^{-1}(z)\right)\big].
\end{eqnarray}

We need only to calculate the limit of $\mathbb{E}\beta_1(z)\big[\bbr_1^{*}\bbD_1^{-1}(z)\bbG^{-1}_n(z)\bbr_1-\frac{1}{n}\mathbb{E}tr  \big(\bbG^{-1}_n(z)\bbD^{-1}(z)\big)\big]$. By similar arguments to Steps 1 and 2, we can get the limit of (\ref{0521(46)}).
Let $\gamma_j(z)=\bbr_j^{*}\bbD_j^{-1}(z)\bbr_j-\frac{1}{n}\mathbb{E}tr\bbD_j^{-1}(z)$.
By (\ref{0521(40)}) and the fact that $\beta_1(z)=b_n(z)\big(1-\beta_1(z)\gamma_1(z)\big)$,
we have
\begin{eqnarray}\label{yang21}
&&\mathbb{E}tr\big(\bbG^{-1}_n(z)\bbD_1^{-1}(z)\big)-\mathbb{E}tr\big(\bbG^{-1}_n(z)\bbD^{-1}(z)\big)\non
&=&\mathbb{E}\beta_1(z)tr\bbG^{-1}_n(z)\bbD_1^{-1}(z)\bbr_1\bbr_1^{*}\bbD_1^{-1}(z)\non
&=&b_n(z)\mathbb{E}\left[(1-\beta_1(z)\gamma_1(z))\bbr_1^{*}\bbD_1^{-1}(z)\bbG^{-1}_n(z)\bbD_1^{-1}(z)\bbr_1\right].
\end{eqnarray}
From Lemma \ref{lem han} and (\ref{yang20}), we get
\begin{eqnarray*}
\big|\mathbb{E}\beta_1(z)\gamma_1(z)\bbr_1^{*}\bbD_1^{-1}(z)\bbG^{-1}_n(z)\bbD_1^{-1}(z)\bbr_1\big|\leq Kn^{-1}.
\end{eqnarray*}
Therefore,
\begin{eqnarray*}
\big|(\ref{yang21})-n^{-1}b_n(z)\mathbb{E}tr\bbD^{-1}_1(z)\bbG^{-1}_n(z)\bbD_1^{-1}(z)\big|\leq Kn^{-1}.
\end{eqnarray*}

Since $\beta_1(z)=b_n(z)-b_n^2(z)\gamma_1(z)+\beta_1(z)b_n^2(z)\gamma_1^2(z)$, we have
\begin{eqnarray}\label{0521(45)}
&&n\mathbb{E}\beta_1(z)\bbr_1^{*}\bbD_1^{-1}(z)\bbG^{-1}_n(z)\bbr_1-\mathbb{E}\beta_1(z)\mathbb{E}tr\bbG^{-1}_n(z)\bbD_1^{-1}(z)\non
&=&\frac{nb_{n}(z)}{c_n\mathbb{E}m_{n}(z)+1}\mathbb{E}\gamma_{1}(z)+b_{n}^{2}(z)\frac{tr\mathbb{E}\bbD_{1}^{-1}(z)}{c_n\mathbb{E}m_{n}(z)+1}\mathbb{E}\gamma_{1}(z)-b_n^2(z)n\mathbb{E}\gamma_1(z)\bbr_1^{*}\bbD_1^{-1}(z)\bbG^{-1}_n(z)\bbr_1\non
&&+b_n^2(z)\Big(n\mathbb{E}\beta_1(z)\gamma_1^2(z)\bbr_1^{*}\bbD^{-1}_1(z)\bbG^{-1}_n(z)\bbr_1 -\big(\mathbb{E}\beta_1(z)\gamma_1^2(z)\big)\mathbb{E}tr\bbG^{-1}_n(z)\bbD^{-1}_1(z)\Big)\non
&&+nb_1(z)\Big(\bbr_1^{*}\bbD_1^{-1}(z)\big[\bbG_n(z)\big]^{-1}\bbr_1-\mathbb{E}tr\bbG^{-1}_n(z)\bbD_1^{-1}(z)\Big)+o(1)\non
&=&\frac{nb_{n}(z)}{c_n\mathbb{E}m_{n}(z)+1}\mathbb{E}\gamma_{1}(z)+b_{n}^{2}(z)\frac{tr\mathbb{E}\bbD_{1}^{-1}(z)}{c_n\mathbb{E}m_{n}(z)+1}\mathbb{E}\gamma_{1}(z)-b_n^2(z)n\mathbb{E}\gamma_1(z)\bbr_1^{*}\bbD_1^{-1}(z)\bbG^{-1}_n(z)\bbr_1\non
&&+b_n^2\Big(\mathbb{E}\big[n\beta_1(z)\gamma_1^{2}(z)\bbr_1^{*}\bbD_1^{-1}(z)\bbG^{-1}_n(z)\bbr_1  -\beta_1(z)\gamma_1^2(z)tr\bbD_1^{-1}(z)\bbG^{-1}_n(z)\big]\Big)\non
&&+b_n^2(z)Cov\big(\beta_1(z)\gamma_1^2(z), tr\bbD^{-1}_1(z)\bbG^{-1}_n(z)\big)+o(1).\non
\end{eqnarray}

By (\ref{0718(4)})
and (\ref{yang21}), we have
\begin{eqnarray*}
\Big|\mathbb{E}\big[n\beta_1(z)\gamma_1^2(z)\bbr_1^{*}\bbD_1^{-1}(z)\bbG^{-1}_n(z)\bbr_1-\beta_1(z)\gamma_1^2(z)tr\bbD_1^{-1}(z)\bbG^{-1}_n(z)\big]\Big|\leq K\delta_n^2.
\end{eqnarray*}

Using (\ref{0718(4)}), (\ref{yang21}), (\ref{yang24}) and (\ref{0718(5)}),
we have
\begin{eqnarray*}
&&\Big|Cov\big(\beta_1(z)\gamma_1^2(z), tr\bbD_1^{-1}(z)\bbG^{-1}_n(z)\big)\Big|\non
&\leq&\big(\mathbb{E}|\beta_1(z)|^4\big)^{1/4}\big(c_{n}\mathbb{E}|\gamma_1(z)|^{8}\big)^{1/4}
\Big(\mathbb{E}\big|tr\bbD_1^{-1}(z)\bbG^{-1}_n(z)-\mathbb{E}tr\bbD_1^{-1}(z)\bbG^{-1}_n(z)\big|^2\Big)^{1/2}\non
&\leq&K\delta_n^{3}n^{-1/4}.
\end{eqnarray*}

Since $\beta_1(z)=b_n(z)-b_n(z)\beta_1(z)\gamma_1(z)$, from (\ref{0718(4)}) and (\ref{0718(5)}), it follows that $\mathbb{E}\beta_1(z)=b_n(z)+O(n^{-1/2})$. Write
\begin{eqnarray*}
&&\mathbb{E}n\gamma_1(z)\bbr_1^{*}\bbD_1^{-1}(z)\bbG^{-1}_n(z)\bbr_1\non
&=&n\mathbb{E}\big[\big(\bbr_1^{*}\bbD_1^{-1}\bbr_1-n^{-1}tr\bbD_1^{-1}(z)\big)
\big(\bbr_1^{*}\bbD_1^{-1}(z)\bbG^{-1}_n(z)\bbr_1-n^{-1}tr\bbD_1^{-1}(z)\bbG^{-1}_n(z)\big)\big]\non
&&+n^{-1}Cov\big(tr\bbD_1^{-1}(z), tr\bbD_1^{-1}(z)\bbG^{-1}_n(z)\big).
\end{eqnarray*}
From (\ref{yang24}), we see that the second term above is $O(n^{-1})$. For the other term $\mathbb{E}\beta_j(z)\big[\bbr_j^{*}\bbD_j^{-1}(z)\bbG^{-1}_n(z)\bbr_j-\frac{1}{n}\mathbb{E}tr  \bbG^{-1}_n(z)\bbD_j^{-1}(z)$, repeat the same steps above, we can get a similar result by replacing the subscript 1 by j. By (\ref{0521(46)}) and (\ref{0521(45)}), we arrive at
\begin{eqnarray}\label{yang26}
&&n\Big(\frac{1}{c_n\big(1+c_{n}\mathbb{E}m_n(z)\big)}+\frac{z}{c_n}\mathbb{E}m_n(z)\Big)\\
&=&\frac{1}{p}\sum_{j=1}^p(W^{(j)}_1+W^{(j)}_2+W^{(j)}_3)+o(1),\nonumber
\end{eqnarray}
where $$W^{(j)}_1=b_n^2(z)n^{-1}\mathbb{E}tr\bbD_j^{-1}(z)\bbG^{-1}_n(z)\bbD_j^{-1}(z), $$ $$W^{(j)}_2=-b_n^2(z)n\mathbb{E}\big[\big(\bbr_1^{*}\bbD_j^{-1}(z)\bbr_1-n^{-1}tr\bbD_j^{-1}(z)\big) \big(\bbr_1^{*}\bbD_j^{-1}(z)\bbG_n^{-1}(z)\bbr_1-n^{-1}tr\bbD_j^{-1}(z)\bbG_n^{-1}(z)\big)\big],$$ and $$W^{(j)}_3=\frac{nb_{n}(z)}{c_n\mathbb{E}m_{n}(z)+1}\mathbb{E}\gamma_{j}(z)+b_{n}^{2}(z)\frac{tr\mathbb{E}\bbD_{j}^{-1}(z)}{c_n\mathbb{E}m_{n}(z)+1}\mathbb{E}\gamma_{j}(z).$$

To calculate the limit of $W^{(j)}_1$, we need to expand $\bbD_j^{-1}(z)$ to the form like (\ref{yang10}). Similar to \cite{BS2004}, we recalculate (\ref{yang8}) and (\ref{1002.3}) by (\ref{0718(3)})--(\ref{0718(4)}). We omit the details here. After these steps, we have
\begin{eqnarray}
\lim_{n\rightarrow\infty}W^{(j)}_1=\frac{z^2m^2(z)\underline{m}^{'}(z)}{cm(z)+1}.
\end{eqnarray}
For $W_2$,
using Lemma \ref{lem yang} on $W_2$, we have
\begin{eqnarray}\label{0521(50)}
&&\\
W^{(j)}_2=\left\{\begin{array}{cc}
W^{(j)}_{2,1}+W^{(j)}_{2,2}+W^{(j)}_{2,3}+W^{(j)}_{2,4}+W^{(j)}_{2,5}, & under \ the \ complex \ case;  \\
W^{(j)}_{2,1}+2W^{(j)}_{2,2}+W^{(j)}_{2,4}+W^{(j)}_{2,5}, & under \ the \ real \ case,  \\
\end{array}
\right.\nonumber
\end{eqnarray}
where $$W^{(j)}_{2,1}=-\frac{1-\mathbb{E}|X_{j1}|^4}{n^2}b_n^2(z)\mathbb{E}tr\bbD_j^{-1}(z)\mathbb{E}tr\big[\bbD_j^{-1}(z)\bbG^{-1}_n(z)\big],$$
$$W^{(j)}_{2,2}=-\frac{1}{n-1}b_n^2(z)\mathbb{E}tr\Big[\bbD_j^{-2}(z)\bbG_n^{-1}(z)\Big],$$ $$W^{(j)}_{2,3}=-\frac{|\mathbb{E}X_{j1}^2|^2}{n}b^2_n(z)\mathbb{E}tr\bbD_j^{-1}(z)\bbG_n^{-1}(z)\big(\bbD_j^{-1}(z)\big)^{*},$$ $$W^{(j)}_{2,4} = -\Big[\frac{\mathbb{E}|X_{j1}|^4}{n}-\frac{2}{n}-n\mathbb{E}\big(\frac{(X_{j1}^{*})^2X_{j2}^2}{n^2}\big)\Big]b^2_n(z)\sum^{n}_{k=1}\mathbb{E}\big(\bbD_j^{-1}(z)\bbG^{-1}_n(z)\big)_{kk}\big(\bbD_j^{-1}(z)\big)_{kk},$$ and $W^{(j)}_{2,5} = o_{L_1}(1)$ uniformly for $j$.

The limits of $W^{(j)}_{2,1}$, $W^{(j)}_{2,2}$ and $W^{(j)}_{2,3}$ can be easily obtained as
\begin{eqnarray*}
&&\lim_{n\rightarrow\infty}W^{(j)}_{2,1}=\Big(\mathbb{E}|X_{j1}|^4-1\Big)\frac{z\underline{m}^2(z)}{\big(1+\underline{m}(z)\big)\big(z(1+\underline{m}(z))-c\big)},
\non
&&\lim_{n\rightarrow\infty}W^{(j)}_{2,2}=-\frac{z\underline{m}^{'}(z)}{\big(1+\underline{m}(z)\big)\big(z+z\underline{m}(z)-c\big)},
\non
&&\lim_{n\rightarrow\infty}W^{(j)}_{2,3}=-|\mathbb{E}X_{11}^2|^2\frac{cm^2(z)}{1+cm(z)}\cdot\frac{1}{(1+cm(z))^2-c|\mathbb{E}X_{11}^2|^2m^2(z)},
\end{eqnarray*}
where the last limit is derived similarly to Lemma 9.

For $W^{(j)}_{2,4}$, similar to deriving the limit of $\frac{\partial^2}{\partial z_1\partial z_2}J_3$ in (\ref{0521(60)}), we have
\begin{eqnarray}\label{0524(20)}
\frac{1}{n}\sum^{n}_{k=1}\mathbb{E}\big[\big(\bbD_{j}^{-1}(z)\big)_{kk}\big]\mathbb{E}\big[\big(\bbD_j^{-1}(z)\bbG_n^{-1}(z)\big)_{kk}\big]
\stackrel{i.p.}{\rightarrow}\frac{\underline{m}^2(z)}{cm(z)+1}.
\end{eqnarray}
By (\ref{0524(20)}) and a simple calculation, we get $\lim_{n\rightarrow\infty}W^{(j)}_{2,4}=-\Big(\mathbb{E}|X_{j1}|^4-3\Big)\frac{z^2m^2(z)\underline{m}^2(z)}{cm(z)+1}$.

For $W_3$, by lemma 8, we have
\begin{eqnarray}\label{0524(50)}
n\mathbb{E}\gamma_j(z)&=&n\mathbb{E}\Big(\frac{1}{n(n-1)}\sum_{k\neq\ell}\big(\bbD_j^{-1}(z)\big)_{k\ell}\Big)\non
&=&\frac{1}{n-1}\mathbb{E}\textbf{1}^{T}\bbD_j^{-1}(z)\textbf{1}-\frac{1}{n-1}\mathbb{E}tr\bbD_j^{-1}(z)\stackrel{i.p.}{\rightarrow}-\frac{1}{z}-\underline{m}(z).
\end{eqnarray}
Then it follows from (\ref{0524(50)}) and (\ref{0524(1)}) that
\begin{eqnarray*}
W^{(j)}_3\stackrel{i.p.}{\rightarrow}&\frac{m(z)+zm(z)\underline{m}(z)-zm^2(z)\underline{m}(z)-z^2m^2(z)\underline{m}^2(z)}{1+cm(z)}.
\end{eqnarray*}

Therefore, it follows from (\ref{yang26}) and calculations above, we can obtain
\begin{eqnarray}\label{0521(61)}
&& \lim_{n\rightarrow\infty}n\Big(\frac{1}{c_n\big(1+c_{n}\mathbb{E}m_n(z)\big)}+\frac{z}{c_n}\mathbb{E}m_n(z)\Big)
\nonumber\\
&&= \ \lim_{n\rightarrow\infty}\frac{1}{p}\sum^p_{j=1}(W^{(j)}_{1}+W^{(j)}_{2,1}+W^{(j)}_{2,2}+W^{(j)}_{2,3}+W^{(j)}_{2,4}+W^{(j)}_3),
\nonumber\\
&& \ \mbox{under the complex case, or};
\nonumber\\
&& = \lim_{n\rightarrow\infty}\frac{1}{p}\sum_{j=1}^p(W^{(j)}_{1}+W^{(j)}_{2,1}+2W^{(j)}_{2,2}+W^{(j)}_{2,4}+W^{(j)}_3),
\nonumber\\
 && \mbox{under the real case}.
\end{eqnarray}

The goal is to find the limit of $M_n^{(2)}(z)=n\Big(\mathbb{E}\underline{m}_n(z)-\underline{m}_{c_n}(z)\Big)$. It has a relation to the limit in (\ref{0521(61)}). We illustrate this point below.

Recall that $m_{c_n}(z)$ and $\underline{m}_{c_n}(z)$ satisfy the following equations
\begin{eqnarray}\label{0521(63)}
m_{c_n}(z)=\frac{1}{1-c_n-c_nzm_{c_n}(z)-z},
\end{eqnarray}
\begin{eqnarray}\label{yang29}
\underline{m}_{c_n}(z)=-\big(z-\frac{c_n}{1+\underline{m}_{c_n}(z)}\big)^{-1}.
\end{eqnarray}
Let $A_n(z)=\frac{1}{c_n\big(1+c_{n}\mathbb{E}m_n(z)\big)}+\frac{z}{c_n}\mathbb{E}\underline{m}_n(z)$.
Since
\begin{eqnarray}\label{0521(67)}
\mathbb{E}\underline{m}_{n}(z)=-\frac{1-c_n}{z}+c_n\mathbb{E}m_n(z),
\end{eqnarray}
With (\ref{0521(67)}), we have
\begin{eqnarray*}
&& A_n(z) = \frac{1}{c_n+c_{n}\mathbb{E}\underline{m}_n(z)+c_{n}\frac{1-c_n}{z}}+\frac{z}{c_n}\mathbb{E}\underline{m}_n(z)\non
&& = \ \mathbb{E}\underline{m}_n(z)\Big[\frac{z}{c_n}+\frac{1}{c_n+c_n(1-c_n)/z+(c_n-1)\mathbb{E}\underline{m}_n(z)}\Big(\frac{1}{\mathbb{E}\underline{m}_n(z)}-\frac{1}{c_n(1+\mathbb{E}\underline{m}_n(z)+\frac{1-c_n}{z})}\Big)\Big].
\end{eqnarray*}
Then it follows that
\begin{eqnarray}\label{yang30}
&&\\
&&\mathbb{E}\underline{m}_n(z)=\Big[-\frac{z}{c_n}\Big(c_n+\frac{c_n(1-c_n)}{z}+(c_n-1)\mathbb{E}\underline{m}_n(z)\Big)\non
&&\ \ \ \ \ \ \ \ \ \ \ \ \ \ \ \ \ \ \ +\frac{1}{c_n(1+\mathbb{E}\underline{m}_n(z)+\frac{1-c_n}{z})}+\Big(c_n+\frac{c_n(1-c_n)}{z}+(c_n-1)\mathbb{E}\underline{m}_n(z)\Big)\frac{A_n(z)}{\mathbb{E}\underline{m}_n(z)}\Big]^{-1}.\nonumber
\end{eqnarray}

By (\ref{yang29}) and (\ref{yang30}), we have
\begin{eqnarray}\label{0524(32)}
\mathbb{E}\underline{m}_n(z)-\underline{m}_{c_n}(z)=A^{-1}B^{-1}(B+A),
\end{eqnarray}
where $A=z-\frac{c_n}{1+\underline{m}_{c_n}(z)}$ and
\begin{eqnarray*}
&&B=-\frac{z}{c_n}\Big(c_n+\frac{c_n(1-c_n)}{z}+(c_n-1)\mathbb{E}\underline{m}_n(z)\Big)
+\frac{1}{c_n(1+\mathbb{E}\underline{m}_n(z)+\frac{1-c_n}{z})}\non
&& + \ \Big(c_n+\frac{c_n(1-c_n)}{z}+(c_n-1)\mathbb{E}\underline{m}_n(z)\Big)\frac{A_n(z)}{\mathbb{E}\underline{m}_n(z)}
\end{eqnarray*}

By the definition of $A_n(z)$ and (\ref{0521(67)}), we know
\begin{eqnarray}\label{0524(30)}
&&\\
&&\frac{1}{c_n\big(1+\mathbb{E}\underline{m}_n(z)+\frac{1-c_n}{z}\big)}=\frac{1}{c_n(1+c_n\mathbb{E}m_n(z))}=A_n(z)-\frac{z}{c_n}\mathbb{E}\underline{m}_n(z).\nonumber
\end{eqnarray}

Then it follows from (\ref{0524(30)}) that
\begin{eqnarray}\label{0524(31)}
&& B+A = -\frac{1}{\underline{m}_{c_n}(z)}-z-(1-c_n)-z\mathbb{E}\underline{m}_n(z)
\nonumber\\
&& + \ A_n(z)\Big(1+\frac{1}{\mathbb{E}\underline{m}_n(z)}\big(c_n+\frac{c_n(1-c_n)}{z}\big)+(c_n-1)\Big)\non
&& = \ c_nz\Big(m_{c_n}(z)-\mathbb{E}m_n(z)\Big)+A_n(z)\Big(1+\frac{1}{\mathbb{E}\underline{m}_n(z)}\big(c_n+\frac{c_n(1-c_n)}{z}\big)+(c_n-1)\Big),\non
\end{eqnarray}
where the last equality uses (\ref{0521(63)}) and (\ref{0521(67)}).

From (\ref{0524(32)}), (\ref{0524(31)}) and the fact that $n\Big(\mathbb{E}\underline{m}_n(z)-\underline{m}_{c_n}(z)\Big) =p\Big(\mathbb{E}m_n(z)-m_{c_n}(z)\Big)$,
we have
\begin{eqnarray}\label{0524(35)}
n\Big(\mathbb{E}\underline{m}_n(z)-\underline{m}_{c_n}(z)\Big)
=\frac{nA_n(z)\Big(1+\frac{U_n}{\mathbb{E}\underline{m}_n(z)}A^{-1}B^{-1}\Big)}{1+A^{-1}B^{-1}z},
\end{eqnarray}
where $U_n=c_n+\frac{c_n(1-c_n)}{z}+(c_n-1)\mathbb{E}\underline{m}_n(z)$.

With tedious but simple calculations, we obtain the limit of each part on the right hand side of (\ref{0524(35)}) as follows.
\begin{eqnarray}\label{0524(36)}
&&\lim_{n\rightarrow\infty}A=z-\frac{c}{1+\underline{m}(z)}, \ \
\lim_{n\rightarrow\infty}B=-z-(1-c)-z\underline{m}(z),\non
&&\lim_{n\rightarrow}\Big(1+\frac{U_n}{\mathbb{E}\underline{m}_n(z)}\Big)=c+\frac{zc+c(1-c)}{z\underline{m}(z)}.
\end{eqnarray}
It follows from (\ref{0524(36)}) that
\begin{eqnarray}\label{0524(37)}
\lim_{n\rightarrow}\frac{1+\frac{U_n}{\mathbb{E}\underline{m}_n(z)}A^{-1}B^{-1}}{1+A^{-1}B^{-1}z}
=-\frac{c(1+\underline{m}(z))\big(z(1+\underline{m}(z))+1-c\big)}{z\underline{m}(z)\Big(\big(z(1+\underline{m}(z))-c\big)^2-c\Big)}.
\end{eqnarray}

Thus, it follows from (\ref{0521(61)}) and (\ref{0524(36)}) that
\begin{eqnarray*}
&&n\Big(\mathbb{E}\underline{m}_n(z)-\underline{m}_{c_n}(z)\Big)\non
&&= \lim_{p\rightarrow\infty}\frac{1}{p}\sum_{j=1}^p\Big(W^{(j)}_1+W^{(j)}_{2,1}+W^{(j)}_{2,2}+W^{(j)}_{2,3}+W^{(j)}_{2,4}+W^{(j)}_{3}\Big)
\nonumber\\
&& \times \ \left(-\frac{c(1+\underline{m}(z))\big(z(1+\underline{m}(z))+1-c\big)}{z\underline{m}(z)\left(\big(z(1+\underline{m}(z))-c\big)^2-c\right)}\right),
\nonumber\\
&& \ \mbox{under the complex random variable case, or}
\nonumber\\
&& = \ \lim_{n\rightarrow\infty}\frac{1}{p}\sum_{j=1}^p \left(W^{(j)}_1+W^{(j)}_{2,1}+2W^{(j)}_{2,2}+W^{(j)}_{2,4}+W^{(j)}_{3}\right)
\nonumber\\
&& \times \ \left(-\frac{c(1+\underline{m}(z))\big(z(1+\underline{m}(z))+1-c\big)}{z\underline{m}(z)\left(\big(z(1+\underline{m}(z))-c\big)^2-c\right)}\right),
\nonumber\\
&& \ \mbox{under the real random variable case.}
\nonumber
\end{eqnarray*}

}

\section{Supplementary: Some Lemmas and Proofs}

In this section, we provide the detailed proofs of some necessary lemmas
used in the proof of Proposition 1, including Lemmas 5--7.
Moreover, three additional lemmas, i.e. Lemmas 8--10 are given and proved
in this section.
\subsection{\textbf{Proof of Lemma 5}}

\begin{proof}
To avoid confusion, denote by $K_q$ a positive constant large enough depending on $q$ only.
By inequality
\begin{eqnarray}\label{z10}
\left|\frac{\|\bby-\bar{\bby}\|^2}{n}-1\right| \le \left| \frac{\sum_{i=1}^{n}X_{i}^2}{n}-1\right|+ \left(\frac{\sum_{i=1}^{n}X_i}{n}\right)^2,
\end{eqnarray}
we have $P\left(B_{n}^{c}(\ep)\right) \le P\left(\left| \frac{\sum_{i=1}^{n}X_{i}^2}{n}-1\right|\ge \frac{\ep }{2}\right)+P\left(\left|\frac{\sum_{i=1}^{n}X_i}{n}\right|\ge \sqrt{\frac{\ep }{2}}\right)$.

By Markov's inequality and Lemma \ref{lem0527}, we have
\begin{eqnarray}\label{z8}
&&\\
P\left(\left|\frac{\sum_{i=1}^{n}X_{i}^2}{n}-1\right|\ge \frac{\ep }{2}\right) &\leq& \frac{2^q\mathbb{E}|\sum_{i=1}^{n}(X_{i}^2-1)|^q}{(n\ep)^q}\non
&\le& \frac{K_q(\sum_{i=1}^{n}\mathbb{E}|X_{i}^2-1|^2)^{q/2}+\mathbb{E}\sum_{i=1}^{n}|X_{i}^2-1|^{q}}{(n\ep)^q}\non
&=&O( n^{-q/2}v^{q/2}_4+ n^{-q+1}v_{2q}).\nonumber
\end{eqnarray}

Similarly, we have
\begin{eqnarray}\label{z9}
\\
&&P\left(\left|\frac{\sum_{i=1}^{n}X_i}{n}\right|\ge \sqrt{\frac{\ep }{2}}\right)\le \frac{K_q ((\sum_{i=1}^{n}\mathbb{E}X^2_{i})^{q}+\mathbb{E}\sum_{i=1}^{n}|X_{i}|^{2q})}{(n\ep)^{2q}}
\nonumber\\
&& = \ O( n^{-q}+ n^{-2q+1}v_{2q}).\nonumber
\end{eqnarray}
which implies $P(B_{n}^{c}(\ep)) =O( n^{-q/2}v^{q/2}_4+ n^{-q+1}v_{2q})$.

Note that
$$\pmb\alpha^{*}\bbA\pmb\alpha-\frac{1}{n}tr\bbA=\left[\pmb\alpha^{*}\bbA\pmb\alpha-\frac{1}{n}tr\bbA\right]I(B_{n}(\ep))+\left[\pmb\alpha^{*}\bbA\pmb\alpha-\frac{1}{n}tr\bbA\right]I(B_{n}^{c}(\ep)).$$

There exists a positive constant $K_q$ such that
$$\mathbb{E}\left|\pmb\alpha^{*}\bbA\pmb\alpha-\frac{1}{n}tr\bbA \right|^q \le K_q \left(\mathbb{E}\left|(\pmb\alpha^{*}\bbA\pmb\alpha-\frac{1}{n}tr\bbA)I(B_n(\ep))\right|^q+\mathbb{E}\left|(\pmb\alpha^{*}\bbA\pmb\alpha-\frac{1}{n}tr\bbA)I(B^{c}_n(\ep))\right|^q\right).$$

By the fact that $|\pmb\alpha^{*}\bbA\pmb\alpha-\frac{1}{n}tr\bbA| \le 2 \|\bbA\|$, we obtain
\begin{eqnarray}\label{han1}
&&\mathbb{E}\left|(\pmb\alpha^{*}\bbA\pmb\alpha-\frac{1}{n}tr\bbA)I(B^{c}_n(\ep))\right|^q\\
 &&\le 2^{q}\|\bbA\|^q P(B^{c}_n(\ep))=O( n^{-q/2}v^{q/2}_4+ n^{-q+1}v_{2q})\|\bbA\|^{q}.
 \nonumber
\end{eqnarray}

Observe that
$$\pmb\alpha^{*}\bbA\pmb\alpha-\frac{1}{n}tr\bbA=\frac{1}{\|\bby-\bar{\bby}\|^2}((\bby-\bar{\bby})^{T}\bbA(\bby-\bar{\bby})-tr\bbA)+\left(\frac{trA}{\|\bby-\bar{\bby}\|^2}-\frac{1}{n}tr\bbA\right)=a_1+a_2.$$

For $0 < \ep <1/2$, there exists a positive constant $K_q$ such that
\begin{eqnarray}
\\
&&\mathbb{E}\left|\left[\pmb\alpha^{*}\bbA\pmb\alpha-\frac{1}{n}tr\bbA\right]I(B_n(\ep))\right|^q \le K_q (\mathbb{E}|a_{1}I(B_n(\ep))|^q +\mathbb{E}|a_{2}I(B_n(\ep))|^q).\nonumber
\end{eqnarray}

Consider $a_2$ first. It is easy to see $\frac{\|\bby-\bar{\bby}\|^2}{n} \ge 1-\ep $ on the event $B_n$, so that
\begin{eqnarray}\label{z11}
&&\mathbb{E}|a_{2}I(B_n(\ep))|^q \leq K_q \left(\frac{tr\bbA}{n}\right)^q \ \mathbb{E}\left|1-\frac{\|\bby-\bar{\bby}\|^2}{n}\right|^q\\
& \le & K_q \left(\frac{tr\bbA}{n}\right)^q \cdot \left\{\mathbb{E} \left|\frac{\sum_{i=1}^{n}X_{i}^2}{n}-1 \right|^q +\mathbb{E}\left|\left(\frac{\sum_{i=1}^{n}X_i}{n}\right)^2\right|^q\right\}\non
&\le & \ K_q \ \left(\frac{tr\bbA}{n}\right)^q \  (n^{-q/2}v^{q/2}_4+ n^{-q+1}v_{2q}),\nonumber
\end{eqnarray}
where the second inequality follows from (\ref{z10}) and the last inequality follows from (\ref{z8}) and (\ref{z9}).

Therefore, we have
\begin{eqnarray}\label{han2}
\mathbb{E}|a_{2}I(B_n(\ep))|^q  \le  K_q \  \left(\frac{tr\bbA}{n}\right)^q \  (n^{-q/2}v^{q/2}_4+ n^{-q+1}v_{2q}).
\end{eqnarray}

Similarly, for $a_1$, by writing $\overline{y}=\bbe\bbe^T\bby$, we have
\begin{eqnarray}
\\
&&\mathbb{E}|a_{1}I(B_n(\ep))|^q\le K_q \frac{1}{n^q}\mathbb{E}|(\bby-\bar{\bby})^{T}\bbA(\bby-\bar{\bby})-tr\bbA|^q\non
&& \le \ \frac{K_q}{n^q} \ \left(\mathbb{E}|\bby^{T}\bbA\bby-tr\bbA|^{q}+E\left|\frac{1}{n}\bby^{T}\bbe\bbe^{T}\bbA\bby+\frac{1}{n}\bby^{T}\bbA\bbe\bbe^{T}\bby\right|^{q}+E\left|\frac{1}{n^2}\bby^{T}\bbe\bbe^{T}\bbA\bbe\bbe^{T}\bby\right|^{q}\right).
\nonumber
\end{eqnarray}
Noting that $\frac{1}{n}tr\bbe\bbe^{T}\bbA = \frac{1}{n}tr\bbA\bbe\bbe^{T}=\frac{1}{n^2}tr\bbe\bbe^{T}\bbA\bbe\bbe^{T} \le \|\bbA\| $ and Lemma 2.2 in \cite{BS2004}, we have
\begin{eqnarray}
\\
&&\mathbb{E}|\bby^{T}\bbA\bby-tr\bbA|^{q}+\mathbb{E}|\frac{1}{n}\bby^{T}\bbe\bbe^{T}\bbA\bby+\frac{1}{n}\bby^{T}\bbA\bbe\bbe^{T}\bby|^{q}+\mathbb{E}|\frac{1}{n^2}\bby^{T}\bbe\bbe^{T}\bbA\bbe\bbe^{T}\bby|^{q}\non
&& \le \ K_q \ (v_{2q}tr(\bbA\bbA^*)^{q}+(v_4 tr(\bbA\bbA^*))^{q/2}).\nonumber
\end{eqnarray}

Hence, we obtain
\begin{eqnarray}\label{han3}
\mathbb{E}|a_{1}I(B_n(\ep))|^q \le K_q  n^{-q} (v_{2q}tr(\bbA\bbA^*)^{q}+(v_4 tr(\bbA\bbA^*))^{q/2}).\end{eqnarray}

Combining (\ref{han1}),(\ref{han2}) and (\ref{han3}) together , we can conclude that
\bea
&&\mathbb{E}\left|\pmb\alpha^{*}\bbA\pmb\alpha-\frac{1}{n}tr \bbA\right|^q
\nonumber\\
&& \leq \ K_q \ \left\{n^{-q} (v_{q}tr(\bbA\bbA^*)^{q}+(v_4 tr(\bbA\bbA^*))^{q/2})+(n^{-q/2}v^{q/2}_4+ n^{-q+1}v_{2q})\|\bbA\|^{q}\right\},
\nonumber
\eea
where $K_q$ is a positive constant depending on $q$ only.

\end{proof}

\subsection{\textbf{Proof of Lemma 6}}

    \begin{proof}
    At first, we evaluate some expectations. Set $\mathbb{E}|\alpha_1|^4=\mu_4$ and $\mathbb{E}(\bar{\alpha}_1\alpha_2)^2=\mu_{12}$ for convenience. Note that
    \begin{eqnarray}\label{1001.1}
    \sum_{i=1}^{n}\alpha_i=0\ \ \sum_{i=1}^{n}(\bar{\alpha}_i \alpha_i)=1\ \ \text{and}\ \ \mathbb{E}(\bar{\alpha}_1 \alpha_1)=\frac{1}{n}.\end{eqnarray}

  It follows that for $i\neq j$
    \begin{eqnarray}\label{y0}
    \mathbb{E}(\bar{\alpha}_i \alpha_j)=\mathbb{E}(\bar{\alpha}_1 \alpha_2)&=&\frac{1}{n-1}\left[\mathbb{E}(\bar{\alpha}_1 \sum_{i=1}^{n} \alpha_i)-\mathbb{E}(\bar{\alpha}_1 \alpha_1)\right]\\
    &=&-\frac{1}{n(n-1)},\nonumber
    \end{eqnarray}
    \begin{eqnarray}\label{y1}
    \mathbb{E}(\bar{\alpha}_1 \alpha_1\bar{\alpha}_2 \alpha_2)&=&\frac{1}{n-1}\mathbb{E}\left(\bar{\alpha}_1 \alpha_1[\sum_{i=1}^{n}(\bar{\alpha}_i \alpha_i)-\bar{\alpha}_1 \alpha_1]\right)\\
    &=&\frac{1}{n(n-1)}-\frac{1}{n-1}\mu_4,\nonumber
    \end{eqnarray}
    \begin{eqnarray}\label{y2}
    \mathbb{E}(\bar{\alpha}_1 \alpha_1\bar{\alpha}_1 \alpha_2) & = & \frac{1}{n-1}\mathbb{E}(\bar{\alpha}_1 \alpha_1\bar{\alpha}_1 (\sum_{i=1}^{n}\alpha_i -\alpha_1))
     \nonumber\\
     & = & -\frac{1}{n-1}\mu_4,
    \end{eqnarray}
    and via (\ref{y1}),
    \begin{eqnarray}\label{y3}
    &&\mathbb{E}(\bar{\alpha}_1 \alpha_1\bar{\alpha}_2 \alpha_3)\\
    &=&\frac{1}{n-2}\mathbb{E}[\bar{\alpha}_1 \alpha_1 \bar{\alpha}_2 \sum_{i=1}^{n}(\bar{\alpha}_i-\alpha_1 -\alpha_2)]\non
    &=&-\frac{1}{n-2}\mathbb{E}(\bar{\alpha}_1 \alpha_1 \bar{\alpha}_2 \alpha_1)-\frac{1}{n-2}\mathbb{E}(\bar{\alpha}_1 \alpha_1 \bar{\alpha}_2 \alpha_2)
    \nonumber
    \end{eqnarray}
    \begin{eqnarray}
    &=&-\frac{1}{(n-1)(n-2)}\mathbb{E}[\bar{\alpha}_1 \alpha_1 \alpha_1 (\sum_{i=1}^{n}\bar{\alpha}_i-\alpha_1)]-\frac{1}{n-2}\mathbb{E}(\bar{\alpha}_1 \alpha_1 \bar{\alpha}_2 \alpha_2)\non
    &=&\frac{1}{(n-1)(n-2)}\mu_4-\frac{1}{n(n-1)(n-2)}+\frac{1}{(n-1)(n-2)}\mu_4\non
    &=&\frac{2}{(n-1)(n-2)}\mu_4-\frac{1}{n(n-1)(n-2)}.\nonumber
    \end{eqnarray}

    Analogously, we can get
    \begin{eqnarray}\label{y7}
    \mathbb{E}(\bar{\alpha}_1 \bar{\alpha}_1 \alpha_2 \alpha_3)&=&\frac{1}{(n-1)(n-2)}\mu_4-\frac{1}{n-2}\mu_{12},
    \end{eqnarray}
    \begin{eqnarray}\label{y4}
    &&\mathbb{E}(\bar{\alpha}_1 \alpha_2\bar{\alpha}_3 \alpha_4)\\
    &=&\frac{1}{n-3}\mathbb{E}[\bar{\alpha}_1 \alpha_2 \bar{\alpha}_3 \sum_{i=1}^{n}(\alpha_i-\alpha_1 -\alpha_2-\alpha_3)]\non
    &=&-\frac{1}{n-3}\mathbb{E}(\bar{\alpha}_1 \alpha_1 \alpha_2 \bar{\alpha}_3)-\frac{1}{n-3}E(\bar{\alpha}_1 \alpha_2 \alpha_2 \bar{\alpha}_3)-\frac{1}{n-3}E(\bar{\alpha}_1 \alpha_2  \bar{\alpha}_3 \alpha_3)\non
    &=&\frac{1}{(n-3)}\left(\frac{2}{n(n-1)(n-2)}-\frac{5}{(n-1)(n-2)}\mu_4+\frac{1}{n-2}\mu_{12}\right).\nonumber
    \end{eqnarray}

Let's calculate $\mu_4$. We claim that $\mu_4=\frac{\mathbb{E}|X_1|^4}{n^2}+o(n^{-2})$. To prove it, we consider the real case only, the complex case can be proved similarly. We below use the same notation $B_n(\ep)$ as Lemma 5. Suppose that $\ep$ is a positive constant such that $\mathbb{P}(B_{n}^c(\ep))=o(n^{-1})$ (by the convergence rate in the law of large numbers that $n\mathbb{P}\left(\left|\frac{\sum_{i=1}^n (X_i^2-1)}{n}\right|\ge \ep \right)\rightarrow 0$, which can be referred to Theorem 28 of \cite{VV1975}). Then, we have
\begin{eqnarray}
&&\\
    &&\left|\mathbb{E}(\frac{|X_1 -\bar{x}|^4}{\|\bby-\bar{\bby}\|^4}-\frac{(X_1 -\bar{x})^4}{n^2})\right|\non
    &\leq&\left|\mathbb{E}\frac{|X_1 -\bar{x}|^4}{n^{2}\|\bby-\bar{\bby}\|^4}(n^2 -\|\bby-\bar{\bby}\|^4)I(B_n(\ep))\right|\non
    &&+\left|\mathbb{E}\Big[\mathbb{E}\Big(\frac{|X_1 -\bar{x}|^4}{\|\bby-\bar{\bby}\|^4}|B^{c}_n(\ep)\Big)I(B^{c}_n(\ep))\Big]\right|+\left|\mathbb{E}\frac{|X_1 -\bar{x}|^4}{n^2}I(B^{c}_n(\ep))\right|\non
    &\leq&\left|\mathbb{E}\frac{|X_1 -\bar{x}|^4}{n^{2}\|\bby-\bar{\bby}\|^4}(n -\|\bby-\bar{\bby}\|^2)(n +\|\bby-\bar{\bby}\|^2)I(B_n(\ep))\right|+\frac{\mathbb{P}(B_n^c(\ep))}{n}+o(n^{-2})\non
    &\leq& \frac{K \ep}{n^2}\mathbb{E}|X_1 -\bar{x}|^4+o(n^{-2}) \leq \  K \ep \frac{1}{n^2}+o(n^{-2})\nonumber,
    \end{eqnarray}
    where the second part of the second inequality follows from
\begin{eqnarray*}
\mathbb{E}\Big[\frac{|X_1 -\bar{x}|^4}{\|\bby-\bar{\bby}\|^4}|B^{c}_n(\ep)\Big]
=\frac{1}{n}\mathbb{E}\Big[\frac{\sum^{n}_{i=1}|X_i-\bar{x}|^4}{\|\bby-\bar{\bby}\|^4}|B^{c}_n(\ep)\Big]
\leq\frac{1}{n},
\end{eqnarray*}
using $\frac{\sum^{n}_{i=1}|X_i-\bar{x}|^4}{\|\bby-\bar{\bby}\|^4}\leq 1$; the third part of the second inequality uses the fact that
$\Big|\mathbb{E}|X_1-\bar x|^4I(B_n^{c}(\varepsilon))\Big|\rightarrow 0$.

  It means the  inequality holds for any $\ep >0$ and $n$ large enough, so we have proved
    $$\mathbb{E}\left(\frac{|X_1 -\bar{x}|^4}{\|\bby-\bar{\bby}\|^4}-\frac{|X_1 -\bar{x}|^4}{n^2}\right)=o(n^{-2}).$$

   In a similar way, we can obtain
        \begin{eqnarray}\label{y5}
    \mu_{12}=\mathbb{E}(\bar{\alpha}_1\bar{\alpha}_1\alpha_2\alpha_2)=\frac{|\mathbb{E}X_{1}^{2}|^2}{n^2}+o(n^{-2}).
    \end{eqnarray}
    It is easy to get
    \begin{eqnarray}\label{y6}
    &&\mathbb{E}(\pmb\alpha^{*}\bbA\pmb\alpha-\frac{1}{n}tr\bbA)(\pmb\alpha^{*}\bbB\pmb\alpha-\frac{1}{n}tr\bbB)\\
    &=&\mathbb{E}(\pmb\alpha^{*}\bbA\pmb\alpha\pmb\alpha^{*}\bbB\pmb\alpha)-\frac{1}{n}(tr\bbA)\mathbb{E}(\pmb\alpha^{*}\bbB\pmb\alpha-\frac{1}{n}tr\bbB)\non
    &&-\frac{1}{n}(tr\bbB)\mathbb{E}(\pmb\alpha^{*}\bbA\pmb\alpha-\frac{1}{n}tr\bbA)-\frac{1}{n^2}tr\bbA tr\bbB\nonumber
    \end{eqnarray}
    and $\mathbb{E}(\pmb\alpha^{*}\bbB\pmb\alpha-\frac{1}{n}tr\bbB)=\mathbb{E}(\bar{\alpha}_1\alpha_2)\sum_{k\neq l}\bbB_{kl}=\frac{1}{n-1}\sum_{k\neq l}\bbB_{kl}E[\bar{\alpha}_1 (\sum_{i=1}^{n}\alpha_i-\alpha_1)]$.

    By (\ref{1001.1}),
    we further have
    \begin{eqnarray}\label{1001.2}
    \mathbb{E}(\pmb\alpha^{*}\bbB\pmb\alpha-\frac{1}{n}tr\bbB)&=&-\frac{1}{n-1}\sum_{k\neq l}\bbB_{kl}\mathbb{E}[\bar{\alpha}_1 \alpha_1]\\
    &=&-\frac{1}{n(n-1)}\sum_{k\neq l}\bbB_{kl},\nonumber
    \end{eqnarray}
by which we can conclude
    $$\frac{1}{n}(tr\bbA)\mathbb{E}(\pmb\alpha^{*}\bbB\pmb\alpha-\frac{1}{n}tr\bbB)=-\frac{tr\bbA}{n^{2}(n-1)}\sum_{k\neq l}\bbB_{kl}.$$
    In the same way, we can get
    $$\frac{1}{n}(tr\bbB)\mathbb{E}(\pmb\alpha^{*}\bbA\pmb\alpha-\frac{1}{n}tr\bbA)=-\frac{tr\bbB}{n^{2}(n-1)}\sum_{k\neq l}\bbA_{kl}.$$
    To calculate $E(\pmb\alpha^{*}\bbA\pmb\alpha\pmb\alpha^{*}\bbB\pmb\alpha)$, we expand the expression as
    \begin{eqnarray}\label{y8}
    &&\mathbb{E}(\pmb\alpha^{*}\bbA\pmb\alpha\pmb\alpha^{*}\bbB\pmb\alpha)\\
    &=&\mathbb{E}(\sum_{i,j}\bar{\alpha}_i\bbA_{ij}\alpha_j \sum_{k,l}\bar{\alpha}_k\bbA_{kl}\alpha_l )=\sum_{i,j,k,l}\mathbb{E}\bar{\alpha}_i\alpha_j\bar{\alpha}_k\alpha_l\bbA_{ij}\bbB_{kl}.\nonumber
    \end{eqnarray}

    To calculate (\ref{y8}), we split it into the following cases:

    1. i=j=k=l, $\sum_{i}(\bar{\alpha}_i\alpha_i)^2 \bbA_{ii} \bbB_{ii}$;

    2. i=j, k=l, i$\neq$k, $\sum_{i,k \atop i\neq k}(\bar{\alpha}_i\alpha_i)(\bar{\alpha}_k\alpha_k)\bbA_{ii}\bbB_{kk}$;

    3. i=j, k$\neq$ l, $\sum_{i,k,l \atop k\neq l}(\bar{\alpha}_i\alpha_i)(\bar{\alpha}_k\alpha_l)\bbA_{ii}\bbB_{kl}$;

    4. i$\neq$j, k=l, $\sum_{i,j,k \atop i\neq j}(\bar{\alpha}_i\alpha_j)(\bar{\alpha}_k\alpha_k)\bbA_{ij}\bbB_{kk}$;

    5. i$\neq$j, k$\neq$l, i=k,j=l $\sum_{i,j \atop i\neq j}(\bar{\alpha}_i\bar{\alpha}_i)(\alpha_j\alpha_j)\bbA_{ij}\bbB_{ij}$;

    6. i$\neq$j, k$\neq$l, i=l,j=k $\sum_{i,j \atop i\neq j}(\bar{\alpha}_i\alpha_i)(\bar{\alpha}_j\alpha_j)\bbA_{ij}\bbB_{ji}$;

    7. i$\neq$j, k$\neq$l, i=k, l$\neq$ j, $\sum_{i,j,l \atop i\neq j \neq l}(\bar{\alpha}_i\alpha_j)(\bar{\alpha}_i\alpha_j)\bbA_{ij}\bbB_{il}$;

    8. i$\neq$j, k$\neq$l, l=j, i$\neq$ k, $\sum_{i,j,k \atop i\neq j \neq k}(\bar{\alpha}_i\alpha_j)(\bar{\alpha}_k\alpha_j)\bbA_{ij}\bbB_{kj}$;

    9. i$\neq$j, k$\neq$l, k=j, i$\neq$ l, $\sum_{i,j,l \atop i\neq j \neq l}(\bar{\alpha}_i\alpha_j)(\bar{\alpha}_j\alpha_l)\bbA_{ij}\bbB_{jl}$;

    10. i$\neq$j, k$\neq$l, i=l, k$\neq$ j, $\sum_{i,j,k \atop i\neq j \neq k}(\bar{\alpha}_i\alpha_j)(\bar{\alpha}_k\alpha_i)\bbA_{ij}\bbB_{ki}$;

    11. i$\neq$j, k$\neq$l, l$\neq $j, i$\neq$ k, $\sum_{i,j,k,l \atop i\neq j \neq k\neq l}(\bar{\alpha}_i\alpha_j)(\bar{\alpha}_k\alpha_l)\bbA_{ij}\bbB_{kl}$.

    For ease of presentation, we still keep $\mu_4$ in the expectations although we have evaluated the value.

    Case 1: \ $\mathbb{E}\left(\sum_{i}(\bar{\alpha}_i\alpha_i)^2 \bbA_{ii} \bbB_{ii}\right)=\mathbb{E}(\bar{\alpha}_1\alpha_1)^2 \sum_{i}\bbA_{ii}\bbB_{ii}=\mu_4 \sum_{i}\bbA_{ii}\bbB_{ii}$.

    Case 2: \ By (\ref{y1}), we have
    \begin{eqnarray*}
    \mathbb{E}\left(\sum_{i,k \atop i\neq k}(\bar{\alpha}_i\alpha_i)(\bar{\alpha}_k\alpha_k)\bbA_{ii}\bbB_{kk}\right)&=&E(\bar{\alpha}_1 \alpha_1\bar{\alpha}_2 \alpha_2)\sum_{i,k \atop i\neq k}\bbA_{ii}\bbB_{kk}\non
    &=&(\frac{1}{n(n-1)}-\frac{1}{n-1}\mu_4)\sum_{i,k \atop i\neq k}\bbA_{ii}\bbB_{kk}.
    \end{eqnarray*}

    Case 3: \ By (\ref{y2}) and (\ref{y3}), we have
    \begin{eqnarray}
    &&\mathbb{E}\sum_{i,k,l \atop k\neq l}(\bar{\alpha}_i\alpha_i)(\bar{\alpha}_k\alpha_l)\bbA_{ii}\bbB_{kl}\\
    &=&\mathbb{E}(\bar{\alpha}_1 \alpha_1\bar{\alpha}_2 \alpha_3)\sum_{i,k,l \atop i\neq k \neq l}\bbA_{ii}\bbB_{kl}+\mathbb{E}(\bar{\alpha}_1 \alpha_1\bar{\alpha}_1 \alpha_2)\sum_{i,l \atop l\neq i}\bbA_{ii}\bbB_{il}\non
    &&+E(\bar{\alpha}_1 \alpha_1\bar{\alpha}_2 \alpha_1)\sum_{i,k \atop k\neq i}\bbA_{ii}\bbB_{ki}\non
    &=&(\frac{2}{(n-1)(n-2)}\mu_4-\frac{1}{n(n-1)(n-2)})\sum_{i,k,l \atop i\neq k \neq l}\bbA_{ii}\bbB_{kl}\non
    &&-\frac{1}{n-1}\mu_4 \sum_{i,l \atop l\neq i}\bbA_{ii}(\bbB_{il}+\bbB_{li}).\nonumber
    \end{eqnarray}

    Case 4: \  Similarly to Case 3, we obtain
\begin{eqnarray*}
 \mathbb{E}\sum_{i,j,k \atop i\neq j}(\bar{\alpha}_i\alpha_j)(\bar{\alpha}_k\alpha_k)\bbA_{ij}\bbB_{kk}&=&(\frac{2}{(n-1)(n-2)}\mu_4-\frac{1}{n(n-1)(n-2)})\sum_{i,k,l \atop l\neq k \neq l}\bbB_{ii}\bbA_{kl}\non
&&-\frac{1}{n-1}\mu_4 \sum_{i,l \atop l\neq i}\bbB_{ii}(\bbA_{il}+\bbA_{li}).
\end{eqnarray*}

    Case 5: \ It follows from (\ref{y5}) that
   \begin{eqnarray*}
  \mathbb{E}\sum_{i,j \atop i\neq j}(\bar{\alpha}_i\bar{\alpha}_i\alpha_j\alpha_j)\bbA_{ij}\bbB_{ij}&=&\mathbb{E}(\bar{\alpha}_1\bar{\alpha}_1\alpha_2\alpha_2)\sum_{i,j \atop i\neq j}\bbA_{ij}\bbB_{ij}=\mu_{12}\sum_{i,j \atop i\neq j}\bbA_{ij}\bbB_{ij}\non
   &=&(\frac{|\mathbb{E}X_1^2|^2}{n^2}+o(n^{-2}))\sum_{i,j \atop i\neq j}\bbA_{ij}\bbB_{ij}.
   \end{eqnarray*}

    Case 6: \ By (\ref{y1}), we have
\bea
&& \mathbb{E}\left[\sum_{i,j \atop i\neq j}(\bar{\alpha}_i\alpha_i)(\bar{\alpha}_j\alpha_j)\bbA_{ij}\bbB_{ji}\right]=\mathbb{E}(\bar{\alpha}_1\alpha_1\bar{\alpha}_2\alpha_2)\sum_{i,j \atop i\neq j}\bbA_{ij}\bbB_{ji}
    \nonumber\\
    && = \ (\frac{1}{n(n-1)}-\frac{1}{n-1}\mu_4)\sum_{i,j \atop i\neq j}\bbA_{ij}\bbB_{ji}.
    \nonumber
    \eea

    Case 7: \ In view of (\ref{y7}), we have
\begin{eqnarray*}
  &&\mathbb{E}\left[\sum_{i,j,l \atop i\neq j\neq l}(\bar{\alpha}_i\alpha_j)(\bar{\alpha}_i\alpha_l)\bbA_{ij}\bbB_{il}\right]
=\mathbb{E}(\bar{\alpha}_1\alpha_2\bar{\alpha}_1\alpha_3)\sum_{i,j,l \atop i\neq j\neq l}\bbA_{ij}\bbB_{il}\non
&=&(\frac{1}{(n-1)(n-2)}\mu_4-\frac{1}{n}\mu_{12})\sum_{i,j,l \atop i\neq j\neq l}\bbA_{ij}\bbB_{il}.
\end{eqnarray*}

    Case 8: \ Similarly to Case 7, we have
\begin{eqnarray*}
 &&\mathbb{E}\left[\sum_{i,j,k \atop i\neq j \neq k}(\bar{\alpha}_i\alpha_j)(\bar{\alpha}_k\alpha_j)\bbA_{ij}\bbB_{kj}\right]=\mathbb{E}(\bar{\alpha}_1\alpha_2\bar{\alpha}_1\alpha_3)\sum_{i,j,k \atop i\neq j \neq k}\bbA_{ij}\bbB_{kj}\non
&=&(\frac{1}{(n-1)(n-2)}\mu_4-\frac{1}{n}\mu_{12})\sum_{i,j,k \atop i\neq j \neq k}\bbA_{ij}\bbB_{kj}.
\end{eqnarray*}

    Case 9: \ By (\ref{y3}), we have
\begin{eqnarray*}
&&\mathbb{E}\left[\sum_{i,j,l \atop i\neq j \neq l}(\bar{\alpha}_i\alpha_j)(\bar{\alpha}_j\alpha_l)\bbA_{ij}\bbB_{jl}\right]=\mathbb{E}(\bar{\alpha}_1 \alpha_1\bar{\alpha}_2 \alpha_3)\sum_{i,j,l \atop i\neq j \neq j}\bbA_{ij}\bbB_{kl}\non
&=&(\frac{2}{(n-1)(n-2)}\mu_4-\frac{1}{n(n-1)(n-2)})\sum_{i,j,l \atop i\neq j \neq k}\bbA_{ij}\bbB_{jl}.
\end{eqnarray*}

    Case 10: \ Similarly to Case 9, we have
\begin{eqnarray*}
&&\mathbb{E}\left[\sum_{i,j,k \atop i\neq j \neq k}(\bar{\alpha}_i\alpha_j)(\bar{\alpha}_k\alpha_i)\bbA_{ij}\bbB_{ki}\right]=\mathbb{E}(\bar{\alpha}_1 \alpha_1\bar{\alpha}_2 \alpha_3)\sum_{i,j,k \atop i\neq j\neq k}\bbA_{ij}\bbB_{ki}\non
&=&(\frac{2}{(n-1)(n-2)}\mu_4-\frac{1}{n(n-1)(n-2)})\sum_{i,j,k \atop i\neq j\neq l}\bbA_{ij}\bbB_{ik}.
\end{eqnarray*}

    Case 11: \ We conclude from (\ref{y4}) that

    \begin{eqnarray}
    &&\\
    &&\mathbb{E}\left[\sum_{i,j,k,l \atop i\neq j \neq k\neq l}(\bar{\alpha}_i\alpha_j)(\bar{\alpha}_k\alpha_l)\bbA_{ij}\bbB_{kl}\right]=\mathbb{E}(\bar{\alpha}_1\alpha_2\bar{\alpha}_3\alpha_4)\sum_{i,j,k,l \atop i\neq j \neq k\neq l}\bbA_{ij}\bbB_{kl}\non
    &=&\frac{1}{(n-3)}\left(\frac{2}{n(n-1)(n-2)}-\frac{5}{(n-1)(n-2)}\mu_4+\frac{1}{n-2}\mu_{12}\right)\sum_{i,j,k,l \atop i\neq j \neq k\neq l}\bbA_{ij}\bbB_{kl}.\nonumber
    \end{eqnarray}

    Summarizing the terms above, we conclude that
    \begin{equation}
    \begin{aligned}
    &\mathbb{E}(\pmb\alpha^{*}\bbA\pmb\alpha-\frac{1}{n}tr\bbA)(\pmb\alpha^{*}\bbB\pmb\alpha-\frac{1}{n}tr\bbB)\\
    &=\sum_{i=1}^{n}\frac{1}{n^2}(\mathbb{E}|X_{1}|^{4}-|E(X_{1}^{2})|^2 -2)\bbA_{ii}\bbB_{ii}+\frac{|\mathbb{E}X_{1}^{2}|^2}{n^2}tr(\bbA\bbB^{T})\non
    &+\frac{1}{n^2}tr(\bbA\bbB)+\frac{1-\mathbb{E}|X_{1}|^{4}}{n^{3}}tr \bbA tr\bbB+\Omega_n,
    \end{aligned}
    \end{equation}
where
    \begin{eqnarray}\label{1001.3}
    &\Omega_n&=\frac{2\mathbb{E}|X_{1}|^{4}-2}{n^4}\sum_{i,k,l \atop i\neq k \neq l}\bbA_{ii}\bbB_{kl}-\frac{\mathbb{E}|X_{1}|^{4}-1}{n^3} \sum_{i,l \atop l\neq i}\bbA_{ii}(\bbB_{il}+\bbB_{li})\\
    &&+\frac{2\mathbb{E}|X_{1}|^{4}-2}{n^4}\sum_{i,k,l \atop i\neq k \neq l}\bbB_{ii}\bbA_{kl}-\frac{\mathbb{E}|X_{1}|^{4}-1}{n^3} \sum_{i,l \atop l\neq i}\bbB_{ii}(\bbA_{il}+\bbA_{li})\non
    &&-\frac{1+|\mathbb{E}x_{1}^2|^2}{n(n-1)(n-2)}\sum_{i,j,l \atop i\neq j\neq l}\bbA_{ij}(\bbB_{il}+\bbB_{li})\non
    &&+\frac{2+|\mathbb{E}X_1^2|^2}{n(n-1)(n-2)(n-3)}\sum_{i,j,k,l \atop i\neq j \neq k\neq l}\bbA_{ij}\bbB_{kl}+o\left(\frac{1}{n}\right).\nonumber
    \end{eqnarray}

We next prove that $\Omega_n=o(1/n)$. Actually, it is straightforward to get the following derivations:
\begin{eqnarray}\label{x1}
&&\sum_{i,k,l \atop i\neq k \neq l}\bbA_{ii}\bbB_{kl}=\sum_{i,k,l}\bbA_{ii}\bbB_{kl}-2\sum_{i,l}\bbA_{ii}\bbB_{il}-\sum_{i, l}\bbA_{ii}\bbB_{ll}+2\sum_{i}\bbA_{ii}B_{ii},
\end{eqnarray}
\begin{eqnarray}\label{x2}
\mathbb{E}|\sum_{i,k,l}\bbA_{ii}\bbB_{kl}|&\le& \mathbb{E}|tr\bbA\bbe^{T}\bbB\bbe|=O(n^2),
\end{eqnarray}
\begin{eqnarray}\label{x4}
&&  \mathbb{E}|\sum_{i, l}\bbA_{ii}\bbB_{il}|\le \mathbb{E}|\sum_{i}\bbA_{ii}\bbe_{i}^{T}\bbB\bbe|\\
&\le& (\sum_{i}\bbA_{ii}\overline{\bbA}_{ii})^{1/2}\left( \mathbb{E}\left[\sum_{i}\bbe^{T}\bbB\bbe_i \overline{\bbe_{i}^T \bbB\bbe}\right]\right)^{1/2}=O(n),
\nonumber
\end{eqnarray}
\begin{eqnarray}\label{x5}
\sum_{i, l}\bbA_{ii}\bbB_{ll}=tr\bbA tr\bbB=O(n^2)\ \ \text{and} \ \ \sum_{i}\bbA_{ii}B_{ii}=O(n).
\end{eqnarray}

We conclude from (\ref{x1})-(\ref{x5}) that $\sum_{i,k,l \atop i\neq k \neq l}\bbA_{ii}\bbB_{kl}=O(n^2)$,
\begin{eqnarray}
&&\sum_{i,k,l \atop i\neq k\neq l}\bbA_{ik}\bbB_{il}=\bbe^{T}\bbA\bbB\bbe-tr\bbA\bbB-2\sum_{i\neq l}\bbB_{ii}\bbB_{il}=O(n),
\end{eqnarray}
and
\begin{eqnarray}
\frac{1}{n^2}\sum_{i,t,k,l \atop i\neq t \neq k\neq l}\bbA_{it}\bbA_{kl}&=&\Big[\frac{1}{n}(\bbe^{T}\bbA\bbe-tr\bbA)\Big]\Big[\frac{1}{n}(\bbe^{T}\bbB\bbe-tr\bbA)\Big]\\
&&-\frac{2}{n^2}\sum_{i,k,l \atop i\neq k\neq l}\bbA_{ik}\bbB_{il}
=O(1).\nonumber
\end{eqnarray}

We have proved the lemma.
\end{proof}

\subsection{\textbf{Proof of Lemma 7}}

\begin{proof}
Set $\bbX_n=(\bby_1, \bby_2, \ldots, \bby_p)^*$ . Denote  the i-th largest eigenvalue of $\bbB_n$ by $\lambda_{i}$ and the i-th largest eigenvalue of $C_n =\frac{1}{n}(\bbX_n-\bar{\bbX}_n)(\bbX_n-\bar{\bbX}_n)^{*}$ by $\nu_i$.
Noticing the trivial inequalities for any positive constant $\sigma$ small enough such that $\mu_1-\sigma >(1+\sqrt{c})^2$ and $\mu_2+\sigma <\mathbf{I}_{(0,1)}(c)(1-\sqrt{c})^2$, we have
\begin{eqnarray}
&&P(\lambda_1 \geq \mu_1)\\
&=& P(\lambda_1 \geq \mu_1, \nu_1 \ge \mu_1 -\sigma)+P(\lambda_1 \geq \mu_1, \nu_1 < \mu_1 -\sigma)\non
&\le&P(\nu_1 \ge \mu_1 -\sigma)+P(|\lambda_1-\nu_1|\ge \sigma)\nonumber
\end{eqnarray}
and
\begin{eqnarray}
&&P(\lambda_p \leq \mu_2)\\
&=& P(\lambda_p \leq \mu_2, \nu_p \le \mu_2 -\sigma)+P(\lambda_p \leq \mu_p, \nu_p > \mu_p -\sigma)\non
&\le&P(\nu_p \le \mu_2 -\sigma)+P(|\lambda_p-\nu_p|\ge \sigma).\nonumber
\end{eqnarray}

For the moment, we assume that

\begin{eqnarray}\label{z1}
P(\|\bbC_n\|=\nu_1 \geq \mu_1-\sigma)=o(n^{-l})
\end{eqnarray}
and
\begin{eqnarray}\label{z2}
P(\lambda_{\min}^{\bbC_n}=\nu_p \leq \mu_2-\sigma)=o(n^{-l}).
\end{eqnarray}
It then suffices to bound $\max_{1\leq i\leq p}P(|\lambda_i-\nu_i|\ge \sigma)$. By Lemma \ref{lem0525} , we have
\begin{eqnarray}\label{z12}
&&\max_{1 \leq i \leq p}|\sqrt{\lambda_{i}}-\sqrt{\nu_{i}}|\leq \|n^{-1/2}(\bbX_n -\bar{\bbX}_n)\|\|\sqrt{n}\bbD_n-\mathbf{I}_n\|\non
&&= \|n^{-1/2}(\bbX_n -\bar{\bbX}_n)\|.\max_{1 \leq i \leq p}\left|\frac{n^{1/2}}{\|\bby_i -\bar{\bby}_i\|}-1\right|.
\end{eqnarray}

In view of the above inequality, it is enough to show that for any fixed $\ep$, we have $P\left(\max_{1 \leq i \leq p}|\frac{n^{1/2}}{\|\bby_i -\bar{\bby}_i\|}-1|\ge \ep\right)=o(n^{-l})$, which can be guaranteed by
$$P\left(\max_{1 \leq i \leq p}|\frac{\|\bby_i -\bar{\bby}_i\|^2}{n}-1|\ge \ep\right)=o(n^{-l}).$$

By the inequality
$$\max_{1 \leq i \leq p} \left|\frac{\|\bby_i -\bar{\bby}_i\|^2}{n}-1\right| \le \max_{1 \leq i \leq p}\left|\frac{\sum_{j=1}^{n}|X_{ij}|^2}{n}-1\right|+\max_{1 \leq i \leq p}
\left|\frac{1}{n}\sum_{j=1}^{n}X_{ij}\right|^2,$$
it suffices to show the following two inequalities:
$$P\left(\max_{1 \leq i \leq p}|\frac{\sum_{j=1}^{n}|X_{ij}|^2}{n}-1|\ge \ep\right)=o(n^{-l})$$
and
\begin{eqnarray}\label{z7}
P\left(\max_{1 \leq i \leq p}|\frac{1}{n}\sum_{j=1}^{n}X_{ij}|\ge \ep\right)=o(n^{-l})
\end{eqnarray}
such that we can obtain (\ref{gg1}) and (\ref{gg2}).

To prove these two inequalities, one can refer to the proof of inequality (9) in \cite{BG2012} for details, we omit them here (one should note that p and n here are of the same order, which is different from \cite{BG2012}, but the proof is almost the same).

To finish the proof, we need to show that (\ref{z1}) and (\ref{z2}). Denoting the i-th largest eigenvalue of $\bbS_n  =\frac{1}{n}\bbX_n\bbX_{n}^{*}$ by $\tau_i$, referring to \cite{BS2004}, we know that
\begin{eqnarray}\label{z4}
P(\|\bbS_n\|=\tau_1 \geq \mu_1-\sigma/2)=o(n^{-l})
\end{eqnarray}
and
\begin{eqnarray}\label{z5}
P(\lambda_{min}^{\bbS_n}=\tau_p \leq \mu_2-\sigma/2)=o(n^{-l}).
\end{eqnarray}

Similarly to (\ref{z12}), we have
\begin{eqnarray}\label{z6}
&&\max_{1 \leq i \leq p}|\sqrt{\tau_{i}}-\sqrt{\nu_{i}}|\leq \|n^{-1/2}\overline{X}_n\|\leq \sqrt{\frac{1}{n}\sum_{i=1}^{p}|\frac{1}{n}\sum_{j=1}^{n}X_{ij}|^2}.
\end{eqnarray}

Combining (\ref{z7}), (\ref{z4}), (\ref{z5}), (\ref{z6}) together, we have (\ref{z1}) and  (\ref{z2}).
\end{proof}

\subsection{\textbf{Lemma 8}}
Suppose that $x_{n}=\bbe/\sqrt{n}=\textbf{1}/\sqrt{n}$ is a unit vector, then for the truncated random variable satisfying (\ref{*}), we have
$\mathbb{E}|x_{n}^{*}\bbD^{-1}(z)x_{n} +\frac{1}{z}|^2\rightarrow 0$.

\begin{proof} By Lemma 5, we obtain for any $2\le r\in \mathbb{N}^+$
$$\mathbb{E}|\bbr_{j}^{*}\bbD_{j}^{-1}(z)x_{n}x_{n}^{T}\bbD_{j}^{-1}(z)\bbr_{j}|^{r}=O(n^{-2}\delta_n^{2r-4}).$$ 
Rewrite it as a martingale
\begin{eqnarray}\label{mart}
&&x_{n}^{*}\bbD^{-1}(z)x_{n}-x_{n}^{*}\mathbb{E}\bbD^{-1}(z)x_{n}\non
&=&\sum_{j=1}^{p}x_{n}^{*}\mathbb{E}_{j}\bbD^{-1}(z)x_{n}-x_{n}^{*}\mathbb{E}_{j-1}\bbD^{-1}(z)x_{n}\non
&=&\sum_{j=1}^{p}x_{n}^{*}\mathbb{E}_{j}(\bbD^{-1}(z)-\bbD_{j}^{-1}(z))x_{n}-x_{n}^{*}\mathbb{E}_{j-1}(\bbD^{-1}(z)-\bbD_{j}^{-1}(z))x_{n}\non
&=&-\sum_{j=1}^{p}(\mathbb{E}_{j} -\mathbb{E}_{j-1})\beta_{j}(z)\bbr_{j}^{*}\bbD_{j}^{-1}(z)x_{n}x_{n}^{*}\bbD_{j}^{-1}(z)\bbr_{j}.
\nonumber
\end{eqnarray}
By Burkholder's inequality and (\ref{0718(5)}), we have
\begin{eqnarray}\label{xxx}
&&\mathbb{E}|x_{n}^{*}\bbD^{-1}(z)x_{n}-x_{n}^{*}\mathbb{E}\bbD^{-1}(z)x_{n}|^{2}\\
&\leq& K\sum_{j=1}^{p}\mathbb{E}|(\mathbb{E}_{j} -\mathbb{E}_{j-1})\beta_{j}(z)\bbr_{j}^{*}\bbD_{j}^{-1}(z)x_{n}x_{n}^{*}\bbD_{j}^{-1}(z)\bbr_{j}|^{2}\non
&\leq& K\sum_{j=1}^{p}(\mathbb{E}|\beta_j(z)|^4)^{1/2}(\mathbb{E}|\bbr_{j}^{*}\bbD_{j}^{-1}(z)x_{n}x_{n}^{*}\bbD_{j}^{-1}(z)\bbr_{j}|^4)^{1/2}
=O(\delta_n^2).
\nonumber
\end{eqnarray}

Thus, we have $\mathbb{E}|x_{n}^{*}\bbD^{-1}(z)x_{n}-x_{n}^{*}\mathbb{E}\bbD^{-1}(z)x_{n}|^2\longrightarrow 0$.
If $\Im z \ge v_0>0$, then $|\beta_j(z)|\le \frac{|z|}{v_0}$, so (\ref{xxx}) can get a sharper bound
\begin{eqnarray}\label{1009.1}\mathbb{E}|x_{n}^{*}\bbD^{-1}(z)x_{n}-x_{n}^{*}\mathbb{E}\bbD^{-1}(z)x_{n}|^{2}=O\left(\frac{1}{n}\right).
\end{eqnarray}

From the proof above, one should note that (\ref{xxx}) and (\ref{1009.1}) hold for $\bbD_j(z)$ and any unit vector.

Note that $\bbD(z)+(c_n z \mathbb{E}m_n (z)+z)\bbI =\sum_{j=1}^{p}r_{j}r_{j}^{*}+c_n z\mathbb{E}m_{n}(z)\bbI$.

Recalling $m_n(z)=-\frac{1}{pz}\sum_{j=1}^p\beta_j(z), \ \ \bbG_n(z)=c_n\mathbb{E}m_n(z)\bbI_n+\bbI_n$, and using the identity
$\bbr_{j}^{*}\bbD^{-1}(z)=\beta_{j}(z)\bbr_{j}^* \bbD_{j}^{-1}(z)$,  we obtain
\begin{eqnarray*}
&&(-z\bbG_n(z))^{-1}-\mathbb{E}\bbD^{-1}(z)\non
&=&-z^{-1}\bbG^{-1}_n(z)\mathbb{E}\big[\big(\sum^{p}_{j=1}\bbr_j\bbr_j^{*}-(-z c_n \mathbb{E}m_n(z)\bbI_n)\big)\bbD^{-1}(z)\big]\non
&=&-z^{-1}\sum^{p}_{j=1}\mathbb{E}\beta_j\big[\bbG^{-1}_n(z)\bbr_j\bbr_j^{*}\bbD_j^{-1}(z)\big]\non
&&-z^{-1}\mathbb{E}\big[\bbG^{-1}_n(z)\big(-c_n z\mathbb{E}m_n(z)\big)\bbI_n\bbD^{-1}(z)\big]\non
&=&-z^{-1}\sum^{p}_{n=1}\mathbb{E}\beta_j\big[\bbG^{-1}_n(z)\bbr_j\bbr_j^{*}\big(\underline{\bbB}_{(j)}^{n}-z\bbI_n\big)^{-1}-\frac{1}{n}\bbG^{-1}_n(z)\mathbb{E}\bbD^{-1}(z)\big]\non
&=&-z^{-1}p\mathbb{E}\beta_1\big[\bbG^{-1}_n(z)\bbr_1\bbr_1^{*}\bbD^{-1}_1(z)
-\frac{1}{n}\bbG^{-1}_n(z)\mathbb{E}\bbD^{-1}(z)\big].
\end{eqnarray*}

Multiplying by $(-x_{n}^*)$ on the left and $x_n$ on the right, we have
\begin{eqnarray*}
&&x_{n}^* \bbD^{-1}(z)x_{n}-x_{n}^* (-z\bbG_n(z))^{-1}x_{n}\non
&=&z^{-1}p\mathbb{E}\beta_1\big[x_{n}^* \bbG^{-1}_n(z)\bbr_1\bbr_1^{*}\bbD^{-1}_1(z)x_{n}
-\frac{1}{n}x_{n}^* \bbG^{-1}_n(z)\mathbb{E}\bbD^{-1}(z)x_{n} \big]\non
&\triangleq&\delta_1 +\delta_2 +\delta_3 ,
\end{eqnarray*}
where $\delta_1 =\frac{p}{z}\mathbb{E}(\beta_1 (z) \alpha_1 (z))$, $\alpha_1 (z)=x_{n}^* \bbG^{-1}_n(z)\bbr_1\bbr_1^{*}\bbD^{-1}_1(z)x_{n}
-\frac{1}{n}x_{n}^* \bbG^{-1}_n(z)\bbD^{-1}_1(z)x_{n}$,
$$\delta_2 =\frac{1}{z}\mathbb{E}\beta_1 (z)x_{n}^* \bbG^{-1}_n(z)(\bbD_{1}^{-1}(z)-\bbD^{-1}(z))x_n,$$
and $\delta_3 =\frac{1}{z}\mathbb{E}\beta_1 (z)x_{n}^* \bbG^{-1}_n(z)(\bbD^{-1}(z)-\mathbb{E}\bbD^{-1}(z))x_n$.

Recalling the notations defined above (\ref{1002.2}) and by the following equalities:
$\delta_1 =\frac{p}{z}\mathbb{E}\tilde \beta_1 (z) \alpha_1 (z)-\frac{p}{z}\mathbb{E}\left[\beta_1 (z) \tilde \beta_1 (z) \varepsilon_1 (z)\alpha_1 (z)\right]$,
$$\tilde \beta_1 (z)=b_n(z)-\frac{1}{n}b_n(z)\tilde \beta_{1}(z)tr(\bbD_1^{-1}(z)-\mathbb{E}\bbD_{1}^{-1}(z)),$$
and $\mathbb{E}\alpha_1 =-(c_n \mathbb{E}m_n(z)+1\big)^{-1}\frac{1}{(n-1)}[\mathbb{E}x^{*}_{n}\bbD_{1}^{-1}(z)x_n+o(1)]$,
it is easy to see $p\mathbb{E}\beta_1 (z) \alpha_1 (z)=[\frac{1}{1+\mathbb{E}\underline{m}_{n}(z)}+o(1)]p\mathbb{E} \alpha_1 (z)$.

Therefore, $\delta_1 =c_n\frac{zm_n (z)}{(c_n zm_n (z)+z)}x_{n}^{*}\mathbb{E}(\bbD_{1}^{-1}(z))x_n +o(1)$.

Similarly to Bai, Miao and Pan(2007), one may have $\delta_2 =o(1)$ and $\delta_3 =o(1)$.
Hence, we obtain $\Big(1- \frac{c_n z m_n (z)}{c_n zm_n (z)+z}\Big)x_{n}^{*}\mathbb{E}(\bbD^{-1}(z))x_n+\frac{1}{c_n zm_n (z)+z} \longrightarrow 0$,
which implies $x_{n}^{*}\mathbb{E}(\bbD^{-1}(z))x_n \longrightarrow -\frac{1}{z}$.

\end{proof}

\begin{rmk}
This is an interesting result that the limit of $\frac{1}{n}\bbe^{T}\mathbb{E}(\bbD^{-1}(z))\bbe$ is independent of the corresponding Stieltjes transform $\underline{m}(z)$. Meanwhile, the limit of  $x_{n}^{*}\mathbb{E}(\bbD^{-1}(z))x_n$ depends on the limit of $x_n$, one can check this by the fact:
$$\frac{1}{n}tr\mathbb{E}(\bbD^{-1}(z))\longrightarrow \underline{m}(z),$$ which depends on the Stieltjes transform $\underline{m}(z)$  different from the result of Lemma 8.
\end{rmk}

\subsection{\textbf{Lemma 9}}
For $z_1$, $z_2 \in \mathcal{C}_u$, we have
\begin{eqnarray}
&&\\
&&\frac{\partial^2}{\partial z_1\partial z_2}J_4\stackrel{i.p.}{\rightarrow}\frac{|\mathbb{E}X_{11}^2|^2cm^{'}(z_1)m^{'}(z_2)}{[(1+c_1m(z_1))(1+cm(z_2))-c|\mathbb{E}X_{11}^2|^2m(z_1)m(z_2)]^2}.\nonumber
\end{eqnarray}

\begin{proof}
From (\ref{yang10}) and bounds, we have
\begin{eqnarray}\label{han10}
\bbD_j^{-1}(z_1)=-\bbH_n(z_1)+b_1(z_1)\bbA(z_1)+\bbB(z_1)+\bbC(z_1).
\end{eqnarray}
Therefore, recalling (\ref{yang8})--(\ref{1002.5}), we have
\begin{eqnarray}
&&\frac{1}{n}tr\big[\mathbb{E}_j\big(\bbD_j^{-1}(z_1)\big)\mathbb{E}_j\big(\bbD_j^{-1}(z_2)^{T}\big)\big]\\
&=&-\bbH_n(z_1)tr\mathbb{E}_j(\bbD_j^{-1}(z_2))^{T}+\frac{1}{n}b_1(z_1)tr\mathbb{E}_j\bbA(z_1)(\bbD_j^{-1}(z_2))^{T}+o(1).\nonumber
\end{eqnarray}

We can write
\begin{eqnarray*}
tr\bbE_j\big(\bbA(z_1)\big)\big(\bbD_j^{-1}(z_2)\big)^{T}=B_1(z_1,z_2)+B_2(z_1,z_2)+B_3(z_1,z_2)+N(z_1 , z_2),
\end{eqnarray*}
where
\begin{eqnarray*}
B_1(z_1,z_2)&=&-tr\sum_{i<j}\bbH_n(z_1)\bbr_i\bbr_i^{*}\mathbb{E}_j\big(\bbD_{ij}^{-1}(z_1)\big)(\beta_{ij}(z_2)\bbD_{ij}^{-1}(z_2)\bbr_i\bbr_i^{*}\bbD_{ij}^{-1}(z_2))^{T}\non
&=&-\sum_{i<j}^{p}\beta_{ij}(z_2)\bbr_i^{*}\mathbb{E}_j\big(\bbD_{ij}^{-1}(z_1)\big)(\bbD_{ij}^{-1}(z_2))^{T}\bar{\bbr}_i\bbr_i^{'}(\bbD_{ij}^{-1}(z_2))^{T}\bbH_n(z_1)\bbr_i;\non
B_2(z_1,z_2)&=&-tr\sum_{i<j}\bbH_n(z_1)n^{-1}\mathbb{E}_j\big(\bbD_{ij}^{-1}(z_1)\big)\big(\bbD_j^{-1}(z_2)-\bbD_{ij}^{-1}(z_2)\big)^T;\non
B_3(z_1,z_2)&=&tr\sum_{i<j}\bbH_n(z_1)\big(\bbr_i\bbr_i^{*}-n^{-1}\bbI_n\big)\mathbb{E}_j\big(\bbD_{ij}^{-1}(z_1)\big)(\bbD_{ij}^{-1}(z_2))^{T};\non
N(z_1 , z_2)&=&tr\mathbb{E}_{j}\sum_{i>j}\bbH_n(z_1)\big(-\frac{1}{n(n-1)}\bbe\bbe^{*}+\frac{1}{n(n-1)}\bbI_n\big)\bbD_{ij}^{-1}(z_1)(\bbD_{j}^{-1}(z_2))^{T}.
\end{eqnarray*}

It is easy to see $N(z_1 , z_2)=O(1)$. We get from (\ref{yang6}) and (\ref{yang7}) that $|B_2(z_1,z_2)|\leq\frac{1+1/v_0}{v_0^2}$.
Similarly to (\ref{yang8}) , we have $\mathbb{E}|B_3(z_1,z_2)|\leq\frac{1+1/v_0}{v_0^3}n^{1/2}$.

Using Lemma 5 and (\ref{0521(20)}), we have, for $i<j$,
\begin{eqnarray}\label{0525(22)}
&&\\
&&\mathbb{E}\big|\beta_{ij}(z_2)\bbr_i^{*}\mathbb{E}_j\big(\bbD_{ij}^{-1}(z_1)\big)(\bbD_{ij}^{-1}(z_2))^{T}\bar{\bbr}_i\bbr_i^{'}(\bbD_{ij}^{-1}(z_2))^{T}\bbH_n(z_1)\bbr_i\non
&&\ \ \ \ \ -b_1(z_2)n^{-2}|\mathbb{E}X_{11}^2|^2tr\big(\mathbb{E}_j\big(\bbD_{ij}^{-1}(z_1)\big)(\bbD_{ij}^{-1}(z_2))^{T}\big)tr\big((\bbD_{ij}^{-1}(z_2))^{T}\bbH_n(z_1)\big)\big|\non
&\leq&Kn^{-1/2}.\nonumber
\end{eqnarray}

By (\ref{yang7}), we have
\begin{eqnarray}\label{0525(23)}
&&\big|tr\big(\mathbb{E}_j\big(\bbD_{ij}^{-1}(z_1)\big)(\bbD_{ij}^{-1}(z_2))^{T}\big)tr\big((\bbD_{ij}^{-1}(z_2))^{T}\bbH_n(z_1)\big)\\
&&\ \ -tr\big(\mathbb{E}_j\big(\bbD_{j}^{-1}(z_1)\big)(\bbD_{j}^{-1}(z_2))^{T}\big)tr\big((\bbD_{j}^{-1}(z_2))^{T}\bbH_n(z_1)\big)\big|\leq Kn.\nonumber
\end{eqnarray}

It follows from (\ref{0525(22)}) and (\ref{0525(23)}) that
\begin{eqnarray*}
&&\mathbb{E}\big|B_1(z_1,z_2)+\frac{j-1}{n^2}b_1(z_2)|\mathbb{E}X_{11}^2|^2 tr\big(\mathbb{E}_j\big(\bbD_{j}^{-1}(z_1)\big)(\bbD_{j}^{-1}(z_2))^{T}\big)tr\big((\bbD_{j}^{-1}(z_2))^{T}\bbH_n(z_1)\big)\big|\non
&&\leq Kn^{1/2}.
\end{eqnarray*}

Analogously, recalling (7.48), we may obtain
\begin{eqnarray}\label{han13}
&&tr\big(\mathbb{E}_j\big(\bbD_j^{-1}(z_1)\big)\big(\bbD_j^{-1}(z_2)\big)^{T}\big)\\
&\times&
\big[1-\frac{j-1}{n^2}|\mathbb{E}X_{11}^2|^2 m_{c_n}(z_1)m_{c_n}(z_2)tr(\bbQ_n(z_2)\bbQ_n(z_1))\big]\non
&&=\frac{1}{z_1z_2}tr\big(\bbQ_n(z_2)\bbQ_n(z_1)\big)+B_6(z_1,z_2),\nonumber
\end{eqnarray}
where $\mathbb{E}|B_6(z_1,z_2)|\leq Kn^{1/2}$.

Rewrite (\ref{han13}) as
\begin{eqnarray}\label{0525(35)}
&&\frac{1}{n}tr\big(\mathbb{E}_j\big(\bbD_j^{-1}(z_1)\big)\big(\bbD_j^{-1}(z_2)\big)^{T}\big)
\big[1-\frac{j-1}{n}|\mathbb{E}X_{11}^2|^2\frac{m_{c_n}(z_1)m^{0}_n(z_2)}{\big(1+\frac{p-1}{n}m_{c_n}(z_2)\big)\big(1+\frac{p-1}{n}m_{c_n}(z_1)\big)}\big]\non
&=&\frac{1}{z_1z_2}\frac{1}{\big(1+\frac{p-1}{n}m_{c_n}(z_1)\big)\big(1+\frac{p-1}{n}m_{c_n}(z_2)\big)}+\frac{1}{n}B_6(z_1,z_2).\non
\end{eqnarray}

Therefore, we have
\begin{eqnarray}
&&\frac{1}{n}tr\big(\mathbb{E}_j\big(\bbD^{-1}_j(z_1)\big)\big(\bbD_j^{-1}(z_2)\big)^{T}\big)\\
&=&\frac{a_n(z_1,z_2)}{z_1z_2 m_{c_n}(z_1)m_{c_n}(z_2)\big[1-\frac{j-1}{p}|\mathbb{E}X_{11}^2|^2a_n(z_1,z_2)\big]}+o(1),\nonumber
\end{eqnarray}
where $a_n(z_1,z_2)=\frac{p}{n}\frac{m_{c_n}(z_1)m_{c_n}(z_2)}{\big(1+\frac{p-1}{n}m_{c_n}(z_1)\big(1+\frac{p-1}{n}m_{c_n}(z_2)\big)}$.

Because the limit of $a_n(z_1, z_2)$ is $a(z_1,z_2)=\frac{c m(z_1)m(z_2)}{\big(1+cm(z_1)\big)\big(1+cm(z_2)\big)}$, we have
\begin{eqnarray}
&&\frac{1}{n}tr\big(\mathbb{E}_j\big(\bbD_j^{-1}(z_1)\big)\big(\bbD_j^{-1}(z_2)\big)^{T}\big)\\
&=&\frac{a_n(z_1,z_2)}{z_1z_2 m(z_1)m(z_2)\big[1-\frac{j-1}{p}|\mathbb{E}X_{11}^2|^2a_n(z_1,z_2)\big]}+o(1).\nonumber
\end{eqnarray}
Therefore, $J_4$ can be written as
\begin{eqnarray*}
J_4=|\mathbb{E}X_{11}^2|^2a_n(z_1,z_2)\frac{1}{p}\sum^{p}_{j=1}\frac{1}{1-\frac{j-1}{p}|\mathbb{E}X_{11}^2|^2 a_n(z_1,z_2)}+B_7(z_1,z_2),
\end{eqnarray*}
where $\mathbb{E}|B_7(z_1,z_2)|\leq Kn^{-1/2}$.

Thus, by (\ref{0525(35)}), the i.p. limit of $J_4$ is $|\mathbb{E}X_{11}^2|^2\int^{a(z_1,z_2)}_{0}\frac{1}{1-|\mathbb{E}X_{11}^2|^2z}dz$. We can then directly get the limit of $\frac{\partial^2}{\partial z_1\partial z_2}J_4$ and we omit the details here.
\end{proof}

\subsection{\textbf{Lemma 10}}
When $v_0= \Im z>0$ is bounded, we have
$$|\mathbb{E}\underline{m}_n(z)-\underline{m}_{c_n}(z)| \le K n^{-1/2}.$$

\begin{proof}
First, by (\ref{0521(6)}), one can prove that
\begin{eqnarray}\label{bixu}
\mathbb{E}|\frac{1}{n}tr\bbD^{-1}(z)-\frac{1}{n}\mathbb{E}tr\bbD^{-1}(z)|^{q} \le K n^{-q/2}, \ \ \forall q \in N_+.
\end{eqnarray}
By $\bbD(z)=\sum_{i=1}^p \bbr_i \bbr_{i}^* -z\mathbf{I},$
we have
$$\mathbf{I}=\sum_{i=1}^p \bbr_i \bbr_{i}^* \bbD^{-1}(z) -z\bbD^{-1}(z)=\sum_{i=1}^p\frac{\bbr_i \bbr_{i}^* \bbD^{-1}_i(z) }{1+\bbr_{i}^* \bbD^{-1}_i(z)\bbr_i}  -z\bbD^{-1}(z).$$

Taking trace and expectation on both sides, then divided by n, we have
\begin{eqnarray}\label{xx1}
1=c_n-c_n\mathbb{E}\frac{1}{1+\bbr_{1}^* \bbD^{-1}_1(z)\bbr_1}-z\mathbb{E}\underline{m}_n(z).
\end{eqnarray}
Denote
\begin{eqnarray}\label{xx2}\rho_n(z)=c_n(\mathbb{E}\frac{1}{1+\bbr_{1}^* \bbD^{-1}_1(z)\bbr_1}-\frac{1}{1+\frac{1}{n}\mathbb{E}tr\bbD^{-1}_1(z)})/\mathbb{E}\underline{m}_n(z).\end{eqnarray}
Combining (\ref{xx1}) and (\ref{xx2}) together, we obtain
$$\mathbb{E}\underline{m}_n(z)=\frac{1}{c_n\frac{1}{1+\mathbb{E}\underline{m}_n(z)}-z+\rho_n(z)}.$$

As we know that $\underline{m}_{c_n}(z)$ satisfies the following equation
$$\underline{m}_{c_n}(z)=\frac{1}{c_n\frac{1}{1+\underline{m}_{c_n}(z)}-z}.$$

Then we have
\bea
&& \mathbb{E}\underline{m}_n(z)-\underline{m}_{c_n}(z)=\frac{c_n(\mathbb{E}\underline{m}_n(z)-\underline{m}_{c_n}(z))\frac{1}{(1+\mathbb{E}\underline{m}_n(z))(1+\underline{m}_{c_n}(z))}}{(c_n\frac{1}{1+\mathbb{E}\underline{m}_n(z)}-z+\rho_n(z))(c_n\frac{1}{1+\underline{m}_{c_n}(z)}-z)}
\nonumber\\
&& + \ (\mathbb{E}\underline{m}_n(z))\underline{m}_{c_n}(z)\rho_n(z).
\nonumber
\eea

Rewrite it as
$$(\mathbb{E}\underline{m}_n(z)-\underline{m}_{c_n}(z))(1-c_n\frac{\mathbb{E}\underline{m}_n(z)\underline{m}_{c_n}(z))}{(1+\mathbb{E}\underline{m}_n(z))(1+\underline{m}_{c_n}(z))})=(\mathbb{E}\underline{m}_n(z))\underline{m}_{c_n}(z)\rho_n(z).$$
Because $\Im z$ is bounded, it is straightforward to obtain $\underline{m}_{c_n}(z)=O(1)$. By the definition of $\rho_n(z)$, there exists a constant C such that
\begin{eqnarray}
&&|\mathbb{E}\underline{m}_n(z)\rho_n(z)| \le C \mathbb{E}|\bbr_{1}^* \bbD^{-1}_1(z)\bbr_1-\frac{1}{n}\mathbb{E}tr\bbD^{-1}_1(z)| \non
&&\le C (\mathbb{E}|\bbr_{1}^* \bbD^{-1}_1(z)\bbr_1-\frac{1}{n}tr\bbD^{-1}_1(z)|+\mathbb{E}|\mathbb{E}\frac{1}{n}tr\bbD^{-1}(z)-\frac{1}{n}tr\bbD^{-1}(z)|\non
&&+\mathbb{E}|\frac{1}{n}tr\bbD^{-1}(z)-\frac{1}{n}tr\bbD^{-1}_1(z)| ) \le C n^{-1/2},\non
\end{eqnarray}
where the last inequality follows from Lemma 5 and (\ref{bixu}).
Similar to (2.19) of \cite{BS2004}, combining with $|\mathbb{E}\underline{m}_n(z)\rho_n(z)|\le C n^{-1/2}$ and $\underline{m}_{c_n}(z)=O(1)$, we have
$$|c_n\frac{\mathbb{E}\underline{m}_n(z)\underline{m}_{c_n}(z))}{(1+\mathbb{E}\underline{m}_n(z))(1+\underline{m}_{c_n}(z))}|<1.$$

Thus, we have proved it, and therefore completed all the proofs in this supplementary document.
\end{proof}


\end{document}